\documentclass[reqno,a4paper,10pt]{amsart}

\usepackage[utf8]{inputenc}
\usepackage{graphicx}
\usepackage{subfigure}
\usepackage{amsmath, amssymb, amsthm}
\usepackage{enumitem}

\usepackage{graphicx}
\usepackage{array}
\usepackage{tikz}
\usepackage{mathrsfs}
\usetikzlibrary{arrows}
\usetikzlibrary{shapes}
\usepackage{cleveref}
\usepackage[english]{babel}

\renewcommand{\theenumi}{(\roman{enumi})}

\newcommand{\ol}[1]{\overline{#1}}
\newcommand{\br}[1]{\left(#1\right)}
\newcommand{\abs}[1]{\left\vert#1\right\vert}
\newcommand{\torus}[1]{\mathbb{S}_{#1}} %surface of genus #1
\newcommand{\Sg}{\mathbb{S}_g} %surface of genus g
\newcommand{\Z}{\mathbb{Z}} %integers
\newcommand{\N}{\mathbb{N}} %natural numbers
\newcommand{\R}{\mathbb{R}} %real numbers
\newcommand{\class}[1]{\mathcal{#1}}
\newcommand{\whp}{whp}%with high probability}
\def\kernel{K}        %kernels of graphs
\def\core{C}        % cores of graphs
\def\complex{Q}       %GF for connected complex graphs
\def\general{S}       %general graphs
\def\outside{U}   %Part outside the complex / non-\comppart\
\def\firstsub{\textsc{1Sub}}
\def\firstcrit{\textsc{1Crit}}
\def\firstsup{\textsc{1Sup}}
\def\intermediate{\textsc{Int}}
\def\secondsub{\textsc{2Sub}}
\def\secondcrit{\textsc{2Crit}}
\def\secondsup{\textsc{2Sup}}
\newcommand{\comppart}{complex part}
\def\component{H}

\DeclareMathOperator{\edgesfirst}{\lambda} %deviation from n/2 edges in first phase transition
\DeclareMathOperator{\edgessecond}{\zeta} %deviation from n edges in second phase transition
\DeclareMathOperator{\intdegree}{\alpha} %average degree in intermediate regime
\DeclareMathOperator{\pl}{pl}   %(size of) planar part
\newcommand{\cf}{w}    %compensation factor
\DeclareMathOperator{\ex}{ex}   %Excess of a graph
\DeclareMathOperator{\deficiency}{d}  %deficiency of a graph
\newcommand{\prob}[1]{\mathbb{P}\left[#1\right]}    %Probability
\newcommand{\expec}[1]{\mathbb{E}\left[#1\right]}    %Expectation
\DeclareMathOperator{\Sub}{s} %Number of cores with the same kernel

\def\dd#1{\,\mathrm{d}{#1}}
\newcommand{\ra}{\rightarrow}

\def\Gimenez{Gim{\'e}nez}
\def\Erdos{Erd\H{o}s}
\def\Renyi{R\'enyi}
\def\Luczak{\L{}uczak}
\def\ER{\Erdos--\Renyi}

\newtheorem{thm}{Theorem}[section]
\newtheorem{prop}[thm]{Proposition}
\newtheorem{coro}[thm]{Corollary}
\newtheorem{lem}[thm]{Lemma}

\newtheorem{question}[thm]{Question}
\newtheorem{problem}[thm]{Problem}
\newtheorem{conjecture}[thm]{Conjecture}
\newtheorem{remark}[thm]{Remark}
\theoremstyle{definition}
\newtheorem{definition}[thm]{Definition}

\newtheorem*{construction*}{Construction}
\newcommand{\proofof}[1]{\subsection{Proof of \Cref{#1}}}

\crefname{thm}{theorem}{theorems}
\crefname{prop}{proposition}{propositions}
\crefname{coro}{corollary}{corollaries}
\crefname{lem}{lemma}{lemmas}
\crefname{definition}{definition}{definitions}
\crefname{question}{question}{questions}
\crefname{problem}{problem}{problems}
\crefname{conjecture}{conjecture}{conjectures}

\newtheoremstyle{claim}%name
{}%Space above
{}%Space below
{\itshape}%Body font
{}%Indent amount
{\bf}%Theorem head font
{.}%Punctuation after theorem head
{.5em}%Space after theorem head
{}%Theorem head spec(can be left empty, meaning ‘normal’)
\theoremstyle{claim}

\crefname{claim}{claim}{claims}

\title{Phase transitions in graphs on orientable surfaces}
\thanks{An extended
  abstract of this paper has been published in  the Proceedings of the
  European Conference on Combinatorics, Graph Theory and Applications
  (EuroComb17), Electronic Notes in Discrete Mathematics 61:687--693, 2017.}
\author{M. Kang, M. Mo{\ss}hammer, P. Spr{\" u}ssel}
\address{Graz University of Technology, Institute of Discrete Mathematics, Steyrergasse 30, 8010 Graz, Austria}
\email{\{kang,mosshammer,spruessel\}@math.tugraz.at}
\thanks{Supported by Austrian Science Fund (FWF): P27290 and W1230 II}


\begin{document}

\begin{abstract}
  Let $\mathbb{S}_g$ be the orientable surface of genus $g$.
  We prove that the component structure of a graph chosen uniformly at random from the class $\mathcal{S}_g(n,m)$ of all graphs on vertex set $[n]=\{1,\dotsc,n\}$
  with $m$ edges embeddable on $\mathbb{S}_g$ features two phase transitions. The first phase 
  transition mirrors the classical phase transition in the Erd\H{o}s--R\'enyi random graph $G(n,m)$ chosen 
  uniformly at random from all graphs with vertex set $[n]$ and $m$ edges. It takes place at $m=\frac{n}{2}+O(n^{2/3})$,
  when a unique largest component, the so-called \emph{giant component}, emerges. The second phase transition occurs at
  $m = n+O(n^{3/5})$, when the giant component covers almost all vertices of the graph. This kind of phenomenon is 
  strikingly different from $G(n,m)$ and has only been observed for graphs on surfaces.
  Moreover, we derive an asymptotic estimation of the number of graphs in $\mathcal{S}_g(n,m)$ throughout the
  regimes of these two phase transitions.
\end{abstract}

\maketitle

\section{Introduction and results}\label{sec:intro}

\subsection{Background and motivation}\label{sec:background}

In their series of seminal papers~\cite{ErdosRenyi59,ErdosRenyi60,ErdosRenyi64,ErdosRenyi66},
\Erdos\ and \Renyi\ studied asymptotic stochastic properties of graphs chosen
according to a certain probability distribution---an approach that laid the foundations for the
classical theory of random graphs. The main questions considered by \Erdos, \Renyi, and many
others are of the following type. Consider the so-called \emph{\ER\ random graph} $G(n,m)$ chosen
uniformly at random from the class $\class{G}(n,m)$ of all graphs on vertex set $[n]:=\{1,\dotsc,n\}$
with $m=m(n)$ edges. What structural properties does $G(n,m)$ have \emph{with
high probability} (commonly abbreviated as \whp), that is, with probability tending to one as $n$
tends to infinity?

One of the most extensively studied properties of random graphs is the component structure. \Erdos\ and
\Renyi~\cite{ErdosRenyi60} proved that the \emph{order} (that is, the number of vertices) of the 
components of $G(n,m)$ changes drastically when $m$ is around $\frac{n}2$; this kind of behaviour 
is widely known as a \emph{phase transition}. The result of \Erdos\ and \Renyi\ states that \whp\ a) if 
the average degree $t := 2\frac{m}{n}$ of $G(n,m)$ is smaller than one, then all components have at most
logarithmic order; b) if $t=1$, the largest component has order $n^{2/3}$; c) if $t\to c>1$, then there 
is a unique component of linear order, while all other components are at most logarithmic. This phenomenon 
became known as the \emph{emergence of the giant component} and was considered by \Erdos\ and \Renyi\ to 
be `one of the most striking facts concerning random graphs'.

While the result of \Erdos\ and \Renyi\ seems to indicate a `double jump' in the order of the largest 
component from logarithmic to order $n^{2/3}$ to linear, Bollob{\'a}s~\cite{Bollobas84} proved that the phase transition is actually `smooth' when
we look more closely at the range of $t$ being around one, that is, when $s:=m-\frac{n}{2}$ is sublinear. 
Bollob{\'a}s' result, which was later improved by {\L}uczak~\cite{luczak1}, shows that the order of the 
largest component changes gradually, depending on whether $s$ has order at most $n^{2/3}$ (known as the \emph{critical 
regime}) or if $s$ has larger order and $s>0$ (the \emph{supercritical regime}) or $s<0$ (the \emph{subcritical
regime}). Subsequently, Aldous~\cite{Aldous97} further improved the result for the critical regime
using multiplicative coalescent processes and inhomogeneous Brownian motion.

In the supercritical regime and in the regime $t>1$, \emph{central limit theorems} and \emph{local limit theorems} 
provide stronger concentration results for the order and the size (that is, the number of edges) of the largest component. The methods used for these
results range from counting techniques~\cite{PittelWormald05,Stepanov70} over Fourier analysis~\cite{coja-oghlan2014}
to probabilistic methods such as Galton-Watson branching processes~\cite{BollobasRiordan2012B}, two-round
exposure~\cite{BehrCOKang10}, or random walks and martingales~\cite{BollobasRiordan2012A}.

Since the pioneering work of \Erdos\ and \Renyi, various random graph models have been introduced and studied. 
A particularly interesting model are random \emph{planar} graphs or, more generally, random graphs that are 
embeddable on a fixed two-dimensional surface. Here, a graph $G$ is called \emph{embeddable} on a surface
$\mathbb{S}$ if $G$ can be drawn on $\mathbb{S}$ without crossing edges.

Graphs embeddable on a surface and graphs \emph{embedded} on a surface---also known as \emph{maps}---have
been studied extensively since the pioneering work of Tutte (see e.g.~\cite{tutte1963-planar-maps})
in view of enumeration~\cite{Chapuy-recurrence,Chapuy-constellations,ChapuyFang-bipartite,Chapuy2011-enumeration-graphs-on-surfaces,MR1980342,gimenez2009,mcdiarmid2005},
random sampling~\cite{PanaSteger-sampling,Boltzmann11,Bodirsky07,Bodirsky08,fusy-sampling,Schaeffer-sampling}, and asymptotic 
properties~\cite{banderier-airy-phenomena,PanaSteger-outerplanar,Bodirsky08,bodirsky2007-cubic-graphs,Chapuy-diameter,Drmota-degree,Drmota-maxdegree,FountPana-3-core,GerkeMcDiarmid04,gerke-planar,gerke-avg-degree,gimenez2009,kang2012,McDiarmid-surfaces,maxdegree,mcdiarmid2005,Mcdiarmid2006,PanaSteger-biconnected,PanaSteger-degree}. 
Maps and embeddable graphs have also shown to have important applications in algebra and geometry 
(see e.g.~\cite{LandoZvonkin} for an overview) and statistical physics~\cite{Brezin78,KangLoebl09,Kasteleyn63}. In some of these applications
(e.g.~\cite{Kasteleyn63}) phase transitions play a crucial role, therefore it is an important
question whether random \emph{embeddable} graphs undergo similar phase transitions as \ER\ random
graphs and if they do, what the \emph{critical behaviour} close to the point of the phase
transition is.

For the order of the largest component of $G(n,m)$, the critical behaviour is described by the results 
of Bollob{\'a}s~\cite{Bollobas84} and {\L}uczak~\cite{luczak1} mentioned above. In order to formally 
state their results, we need to introduce some notation. A connected graph is called \emph{tree} if it
has no cycles, \emph{unicyclic} if it contains precisely one cycle, and \emph{complex} (or \emph{multicyclic}) otherwise.
Given a graph $G$, we enumerate its components as $\component_i=\component_i(G)$, $i=1,2,\dotsc$, in such a way that
they are ordered from large to small, that is, the orders $\abs{\component_1},\abs{\component_2},\dotsc$ of the components satisfy $\abs{\component_i}\ge\abs{\component_j}$
whenever $i<j$. We say that $\component_i$ is the \emph{$i$-th largest component} of $G$.

The results of Bollob{\'a}s and {\L}uczak can now be described as follows (for all order 
notation in the following, see~\Cref{def:landau}). If $m$ is smaller than $\frac{n}{2}$ and satisfies $m-\frac{n}{2}=\omega(n^{2/3})$,
then \whp\ all components of $G(n,m)$ have order $o(n^{2/3})$. Once $m-\frac{n}{2} = O(n^{2/3})$, several components of order
$\Theta_p(n^{2/3})$ appear simultaneously. Finally, if $m$ becomes even larger,
then the largest component $\component_1$ \whp\ has order $\omega(n^{2/3})$, while every
other component has order $o(n^{2/3})$ \whp. If we view this development as a process, this
means that all components of order $\Theta_p(n^{2/3})$ that appeared when $m-\frac{n}2 = O(n^{2/3})$
later merge into a single component that is then the unique component of order $\omega(n^{2/3})$.
This component is usually referred to as the \emph{giant component}.

\begin{thm}[\cite{Bollobas84,luczak1}]\label{thm:ER}
  Let $m=\br{1+\edgesfirst n^{-1/3}}\frac{n}{2}$, where $\edgesfirst=\edgesfirst(n)=o(n^{1/3})$, and let
  $\component_i=\component_i(G)$, $i=1,2,\dotsc$, be the $i$-th largest component of $G=G(n,m)$.
  \begin{enumerate}
  \item\label{Gnm:sub}
    If $\edgesfirst\to-\infty$, then for every $i\in\N\setminus\{0\}$ \whp\ $\component_i$ is a tree and has order
    \begin{equation*}
      \br{2+o\br{1}}\frac{n^{2/3}}{\edgesfirst^2}\log\br{-\edgesfirst^3}.
    \end{equation*}
  \item\label{Gnm:crit}
    If $\edgesfirst\to c$ for a constant $c\in\R$, then for every $i\in\N\setminus\{0\}$ the order of $\component_i$ is
    \begin{equation*}
      \Theta_p\br{n^{2/3}}.
    \end{equation*}
    Furthermore, the probability that $\component_i$ is complex is bounded away both from 0 and 1.
  \item\label{Gnm:super}
    If $\edgesfirst\to\infty$, then \whp\ the largest component $\component_1$ of $G$ is complex and has
    order
    \begin{equation*}
      (2+o(1))\edgesfirst n^{2/3}.
    \end{equation*}
    For $i\ge 2$, \whp\ $\component_i$ is a tree of order $o(n^{2/3})$.
  \end{enumerate}
\end{thm}

Returning to embeddable graphs, we call a graph \emph{planar} if it is embeddable on the sphere and
denote by $P(n,m)$ the graph chosen uniformly at random from the class $\class{P}(n,m)$ of all planar graphs with
vertex set $[n]$ and $m$ edges. Kang and {\L}uczak~\cite{kang2012} proved that $P(n,m)$
features a similar phase transition as $G(n,m)$, that is, the giant component emerges at
$m = \frac{n}2+O(n^{2/3})$.

\begin{thm}[\cite{kang2012}]\label{thm:planar:first}
  Let $m=\br{1+\edgesfirst n^{-1/3}}\frac{n}{2}$, where $\edgesfirst=\edgesfirst(n)=o(n^{1/3})$, and let
  $\component_i=\component_i(G)$, $i=1,2,\dotsc$, be the $i$-th largest component of $G=P(n,m)$. For every $i\in\N\setminus\{0\}$
  \whp
  \begin{equation*}
    \abs{\component_i} =
    \begin{cases}
      \br{2+o\br{1}}\frac{n^{2/3}}{\edgesfirst^2}\log\br{-\edgesfirst^3} & \text{if }\edgesfirst\to-\infty,\\
      \Theta\br{n^{2/3}} & \text{if }\edgesfirst\to c\in\R,\\
      (1+o(1))\edgesfirst n^{2/3} & \text{if }\edgesfirst\to\infty\text{ and }i=1,\\
      \Theta(n^{2/3}) & \text{if }\edgesfirst\to\infty\text{ and }i\ge 2.
    \end{cases}
  \end{equation*}
\end{thm}

The main difference to the \ER\ random graph lies in the case $\edgesfirst\to\infty$. In this regime,
the largest component of $P(n,m)$ is roughly half as large as the largest component of $G(n,m)$.
On the other hand, the order of the second largest component (or more generally, of the $i$-th
largest component for every fixed $i\ge2$) is much larger in $P(n,m)$ than in $G(n,m)$.

This behaviour, however, is not the most surprising feature of random planar graphs. Indeed, Kang
and {\L}uczak~\cite{kang2012} discovered that there is a second phase transition at
$m = n+O(n^{3/5})$, which is when the giant component covers almost all vertices. Such a behaviour
is not observed for \ER\ random graphs, where the number of vertices outside the giant component is
linear in $n$ as long as $m$ is linear.

\begin{thm}[\cite{kang2012}]\label{thm:planar:second}
  Let $m=\br{2+\edgessecond n^{-2/5}}\frac{n}{2}$, where $\edgessecond=\edgessecond(n)=o(n^{2/5})$.
  Then \whp\ the largest component $\component_1$ of $P(n,m)$ is complex and
  \begin{equation*}
    n-\abs{\component_1} =
    \begin{cases}
      \Theta\br{\abs{\edgessecond}n^{3/5}} & \text{if }\edgessecond\to-\infty,\\
      \Theta\br{n^{3/5}} & \text{if }\edgessecond\to c\in\R,\\
      \Theta\br{\edgessecond^{-3/2}n^{3/5}} & \text{if }\edgessecond\to\infty\text{ and }\edgessecond=o(n^{1/15}).
    \end{cases}
  \end{equation*}
\end{thm}

Given that this second phase transition has only been observed for random planar graphs, the
fundamental question that is raised by \Cref{thm:planar:second} is whether this is an intrinsic
phenomenon of planar graphs.

\begin{question}\label{question:main}
  Which other classes of graphs feature a phase transition analogous to \Cref{thm:planar:second}?
\end{question}

Canonical candidates for classes that lie `between' $\class{P}(n,m)$ and $\class{G}(n,m)$ are
graphs that are embeddable on a surface of fixed positive genus. In this paper, we consider graphs embeddable on the \emph{orientable} 
surface $\Sg$ with genus $g\in\N$. Let $\class{\general}_g(n,m)$ be the class of graphs 
with vertex set $[n]$ and $m$ edges that are embeddable on $\Sg$. (Of course, $\class{\general}_0(n,m)=\class{P}(n,m)$.) 
One of the main results of this paper is that for every fixed $g$,
the answer to \Cref{question:main} is positive for the class $\class{\general}_g(n,m)$.

For $m=\left\lfloor \mu n\right\rfloor$ with $\mu\in(1,3)$, \Gimenez\ and Noy~\cite{gimenez2009}
showed, among several other results, that \whp\ $P(n,m)$ has a component that covers all but
finitely many vertices. Observe that \Cref{thm:planar:second} leaves a gap of order $\Theta(n^{1/3})$ 
to the `dense' regime considered by \Gimenez\ and Noy. Subsequently, Chapuy, Fusy, \Gimenez, Mohar, and
Noy~\cite{Chapuy2011-enumeration-graphs-on-surfaces} proved analogous results in the dense regime for $\general_g(n,m)$.

\subsection{Main results}\label{sec:main}

This paper is the first to determine the component structure of $\general_g(n,m)$
for arbitrary $g\ge0$ in the `sparse' regime $m\le(1+o(1))n$. In terms of phase transitions,
the component structure of $\general_g(n,m)$ features particularly interesting phenomena in 
this regime, similar to $P(n,m)$. To derive these phenomena, we use a wide range of complementary 
methods from various fields (see \Cref{sec:techniques} for more details). 

With this paper, we strive to provide a deeper
understanding of the evolution of graphs embeddable on $\Sg$ for fixed $g$. Moreover, we pave
a way to better understand embeddability of random graphs, in particular a) the `typical' genus
of $G(n,m)$ when $m=m(n)$ is given and b) the evolution of graphs on a surface of non-constant 
genus $g=g(n)$.

The main contributions of this paper are fourfold. We determine the order and structure of the
largest components of a graph $\general_g(n,m)$ chosen uniformly at random from
$\class{\general}_g(n,m)$, where the number $m$ of edges is a) around $\frac{n}{2}$, b) around $n$,
or c) in between the previous two regimes. Moreover, we determine d) the asymptotic number of graphs 
in $\class{\general}_g(n,m)$ for all the aforementioned regimes.

Our first main result describes the appearance of the unique giant component in $\general_g(n,m)$. Similar
to various random graph models including \ER\ random graphs and random planar graphs (see
\Cref{thm:ER,thm:planar:first}), the critical range for the number of edges for the appearance of
the giant component is $m=\frac{n}2+O(n^{2/3})$. Below this range, the $i$-th largest component (for each $i\ge1$) of
$\general_g(n,m)$ \whp\ is a tree of order $o(n^{2/3})$. In the critical range, several components
of order $\Theta_p(n^{2/3})$ appear simultaneously. After the critical range, $\general_g(n,m)$
\whp\ has a unique component of order $\omega(n^{2/3})$ which in addition is complex and has genus $g$,
that is, it is embeddable on $\Sg$, but not on $\torus{g-1}$.

\begin{thm}\label{main1}
  Let $m=\br{1+\edgesfirst n^{-1/3}}\frac{n}{2}$, where $\edgesfirst=\edgesfirst(n)=o(n^{1/3})$, and denote by
  $\component_i=\component_i(G)$, $i=1,2,\dotsc$, the $i$-th largest component of $G=\general_g(n,m)$. For every $i\in\N\setminus\{0\}$
  \whp\ the following holds.
  \begin{enumerate}
  \item\label{firstsub}
    If $\edgesfirst\to-\infty$, then $\component_i$ is a tree of order
    \begin{equation*}
      \br{2+o\br{1}}\frac{n^{2/3}}{\edgesfirst^2}\log\br{-\edgesfirst^3}.
    \end{equation*}
  \item\label{firstcrit}
    If $\edgesfirst\to c$ for a constant $c\in\R$, then the probability that $G$ has
    complex components is bounded away both from 0 and 1.
    The $i$-th largest component has order
    \begin{equation*}
      \Theta_p\br{n^{2/3}}.
    \end{equation*}
  \item\label{firstsup}
    If $\edgesfirst\to\infty$, then $\component_1$ is complex and has order
    \begin{equation*}
      \edgesfirst n^{2/3}+O_p(n^{2/3}).
    \end{equation*}
    For $i\ge 2$, we have $\abs{\component_i} = \Theta_p(n^{2/3})$.

    Moreover, $G$ has $O_p(1)$ complex components. The probability that $G$ has at least $i$ complex 
    components is bounded away both from 0 and 1. If $G$ has at least $i$ complex
    components, then the $i$-th largest \emph{complex} component (by this we mean $\component_i(\complex_G)$, where
    $\complex_G$ is the union of all complex components of $G$) has order $\Theta_p(n^{2/3})$.

    Furthermore, if $g\ge 1$, then whp $\component_1$ is not embeddable on $\torus{g-1}$, while all
    other components of $G$ are planar.
 \end{enumerate}
\end{thm}

Comparing the special case of $g=0$ in \Cref{main1} with \Cref{thm:planar:first}, the following
discrepancies are apparent. Firstly, in the critical regime $\edgesfirst\to c\in\R$, \Cref{main1}\ref{firstcrit}
yields components of order $\Theta_p(n^{2/3})$ compared to $\Theta(n^{2/3})$ claimed by
\Cref{thm:planar:first}. The same holds for the orders of $\component_i$ for $i\ge 2$ in the supercritical
regime $\edgesfirst\to\infty$. Both points are due to minor mistakes in~\cite{kang2012}; the proofs
given there in fact yield order $\Theta_p(n^{2/3})$ instead of the claimed $\Theta(n^{2/3})$. Secondly, the
error term in the order of the giant component given in \Cref{main1}\ref{firstsup} is stronger than
the one from \Cref{thm:planar:first}. Finally, \Cref{main1}\ref{firstsup} tells us that for
positive genus, the giant component is not only the unique largest component but also the unique
\emph{non-planar} one.

Our second main result describes the time when the giant component covers almost all vertices. The
critical phase for the number of edges for this phenomenon is $m=n+O(n^{3/5})$. Here, the number of
vertices \emph{outside} the giant component changes from $\omega(n^{3/5})$ for $m$ below the
critical range to $\Theta(n^{3/5})$ within the critical range to $o(n^{3/5})$ beyond the critical
range.

\begin{thm}\label{main2}
  Let $m=\br{2+\edgessecond n^{-2/5}}\frac{n}{2}$, where $\edgessecond=\edgessecond(n)=o(n^{2/5})$.
  Then \whp\ the largest component $\component_1$ of $\general_g(n,m)$ is complex. Furthermore, for
  \begin{equation*}
    r(n) :=
    \begin{cases}
      \abs{\edgessecond}n^{3/5} & \text{if }\edgessecond\ra-\infty,\\
      n^{3/5} & \text{if }\edgessecond\ra c\in\R,\\
      \edgessecond^{-3/2}n^{3/5} & \text{if }\edgessecond\ra\infty\text{ and }\edgessecond=o((\log n)^{-2/3}n^{2/5}),
    \end{cases}
  \end{equation*}
  we have $n-\abs{\component_1}=O_p\br{r(n)}$ and \whp\ $n-\abs{\component_1}=\Omega(r(n))$.
\end{thm}

The main improvement of \Cref{main2} in comparison to \Cref{thm:planar:second} (the
corresponding result for $g=0$) is that \Cref{thm:planar:second}
only deals with the case $\zeta = o(n^{1/15})$ and therefore leaves a
gap to the dense regime $m=\lfloor\mu n\rfloor$ with $\mu\in(1,3)$ that has been
covered in~\cite{Chapuy2011-enumeration-graphs-on-surfaces,gimenez2009}.
\Cref{main2} closes this gap up to a factor $(\log n)^{2/3}$. Additionally, \Cref{main2}
provides a correction of the upper bound given in~\cite{kang2012} on the
number of vertices outside the giant component. In~\cite{kang2012}, the
upper bound was obtained with the help of an intermediate result (Theorem 2(iv) in~\cite{kang2012}) about
the structure of the complex part (see \Cref{sec:asymptotics} for
a definition). However, this intermediate result does not apply in the
regime $m\sim n$. \Cref{main2} provides a slightly weaker
upper bound that is of larger order than the lower bound (albeit the
orders differ by less than every \emph{fixed} growing function).

Our third main result covers the case when the number of edges is between the regimes of
the two phase transitions, that is, the average degree of the graph is between one and two. In this
`intermediate' regime, the largest component is complex, has genus $g$, and its order is
linear both in $n$ and in the average degree of the graph.

\begin{thm}\label{main3}
  Let $m=\intdegree\frac{n}{2}$, where $\intdegree=\intdegree(n)$ converges to a constant in
  $(1,2)$, and let $\component_i=\component_i(G)$, $i=1,2,\dotsc$, be the $i$-th largest component of
  $G=\general_g(n,m)$. Then
  \begin{equation*}
    \abs{\component_1} = \br{\intdegree-1}n+O_p\br{n^{2/3}}.
  \end{equation*}
  For $i\ge 2$, we have $\abs{\component_i} = \Theta_p(n^{2/3})$.

  Furthermore, if $g\ge 1$, then \whp\ $\component_1$ is not embeddable on $\torus{g-1}$, while all other
  components are planar.
\end{thm}

\begin{figure}[htbp]
  \begin{tikzpicture}[xscale=1.5,yscale=3]
    \draw[->] (0,0) -- (3.25,0) node[anchor=120] {$\intdegree$};
    \draw (0,0) -- (0,1) node[anchor=225] {$\abs{H_1}/n$};

    \draw (1.5,-.3) node {$G(n,m)$};

    \foreach \x in {0,1,2,3} \draw (\x,-.04) node[anchor=north] {$\x$} -- (\x,0);
    \foreach \y in {0,0.5,1} \draw (-.08,\y) node[anchor=east] {$\y$} -- (0,\y);

    \draw[very thick] (0,0) -- (1,0) -- (1.1,.176134) -- (1.2,.313698) -- (1.3,.42297) -- (1.4,.511011) -- (1.5,.582812) -- (1.6,.641981) -- (1.7,.691186) -- (1.8,.73243) -- (1.9,.767244) -- (2,.796812) -- (2.1,.822065) -- (2.2,.843739) -- (2.3,.862423) -- (2.4,.878596) -- (2.5,.892645) -- (2.6,.904889) -- (2.7,.915593) -- (2.8,.924975) -- (2.9,.933219) -- (3,.94048) -- (3.1,.946888) -- (3.2,.952555);

    \begin{scope}[shift={(4.25,0)}]
      \draw[->] (0,0) -- (3.25,0) node[anchor=120] {$\intdegree$};
      \draw (0,0) -- (0,1) node[anchor=225] {$\abs{H_1}/n$};

      \draw (1.5,-.3) node {$\general_g(n,m)$};

      \foreach \x in {0,1,2,3} \draw (\x,-.04) node[anchor=north] {$\x$} -- (\x,0);
      \foreach \y in {0,0.5,1} \draw (-.08,\y) node[anchor=east] {$\y$} -- (0,\y);

      \draw[very thick] (0,0) -- (1,0) -- (2,1) -- (3,1);
    \end{scope}
  \end{tikzpicture}
  \caption{Rescaled order of the largest component of $G(n,m)$ and of $\general_g(n,m)$.}
  \label{fig:largest}
\end{figure}
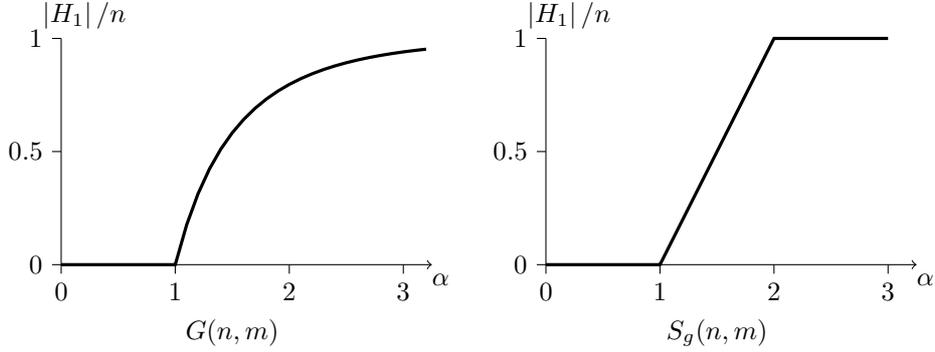

In the intermediate regime, or more generally, for $m=\intdegree\frac{n}{2}$ with $\intdegree>1$,
the classical \ER\ random graph $G(n,m)$ \whp\ has a largest component of order $(1+o(1))\beta n$,
where $\beta$ is the unique positive solution of the equation
\begin{equation*}
  1-\beta = e^{-\alpha\beta}.
\end{equation*}
In particular, as long as $\intdegree>1$ is a constant, the largest component of $G(n,m)$ will
leave a linear number of vertices uncovered, see \Cref{fig:largest}.
Indeed, Karp~\cite{Karp90} proved that the components of $G(n,m)$ can be explored via a
Galton-Watson branching process with offspring distribution $\text{Po}(\intdegree)$; the survival
property of such a process is given by $\beta$ above, yielding order $(1+o(1))\beta n$ of the
largest component. For graphs on surfaces, however, there is no such simple approach to explore 
components.

As our last main result, we derive the asymptotic number of graphs embeddable on $\Sg$.

\begin{thm}\label{main4}
  For $n\to\infty$, the number of graphs in $\class{\general}_g(n,m)$ is asymptotically
  given as follows.
  \begin{enumerate}
  \item\label{enum:first}
    If $m = \br{1+\edgesfirst n^{-1/3}}\frac{n}{2}$, where $\edgesfirst=\edgesfirst(n)=o(n^{1/3})$, then
    \begin{equation*}
      \abs{\class{\general}_g(n,m)} =
      \begin{cases}
        \frac{1+o(1)}{\pi^{1/2}e^{3/4}}\br{\frac{e}{1+\edgesfirst n^{-1/3}}}^{\frac{n}2+\frac{\edgesfirst}{2}n^{2/3}}n^{\frac{n}2+\frac{\edgesfirst}{2}n^{2/3}-\frac12} & \text{if }\edgesfirst\to-\infty,\\
        \Theta(1)e^{\frac{n}2-\frac{\edgesfirst^2}{4}n^{1/3}}n^{\frac{n}2+\frac{\edgesfirst}{2}n^{2/3}-\frac12} & \text{if }\edgesfirst\to c\in\R,\\
        \exp\br{O(\edgesfirst)}\br{\frac{e}{1-\edgesfirst n^{-1/3}}}^{\frac{n}2-\frac{\edgesfirst}{2}n^{2/3}}n^{\frac{n}2+\frac{\edgesfirst}{2}n^{2/3}-\frac12} & \text{if }\edgesfirst\to\infty.
      \end{cases}
    \end{equation*}
  \item\label{enum:intermediate}
    If $m=\intdegree\frac{n}{2}$, where $\intdegree=\intdegree(n)$ converges to a constant in $(1,2)$, then
    \begin{equation*}
      \abs{\class{\general}_g(n,m)} =
      \exp\br{O\big(n^{1/3}\big)}\br{\frac{e}{2-\intdegree}}^{(2-\intdegree)\frac{n}{2}}n^{\intdegree\frac{n}{2}}.
    \end{equation*}
  \item\label{enum:second}
    If $m = \br{2+\edgessecond n^{-2/5}}\frac{n}{2}$, where $\edgessecond=\edgessecond(n)=o(n^{2/5})$, then
    \begin{equation*}
      \abs{\class{\general}_g(n,m)} =
      \begin{cases}
        \exp\br{O\big(\abs{\edgessecond}^{-2/3}n^{3/5}\big)}\br{\frac{\abs{\edgessecond}}{e}}^{\frac{\edgessecond}{2}n^{3/5}}n^{n+\frac3{10}\edgessecond n^{3/5}} & \text{if }\edgessecond\to-\infty,\\
        \exp\br{O\big(n^{3/5}\big)}n^{n+\frac3{10}\edgessecond n^{3/5}} & \text{if }\edgessecond\to c\in\R,\\
        \exp\br{O\big(\edgessecond n^{3/5}\big)}\edgessecond^{-\frac34\edgessecond n^{3/5}} n^{n+\frac3{10}\edgessecond n^{3/5}} & \text{if }\edgessecond\to\infty \text{ and}\\
        & \hspace{-15pt}\edgessecond=o((\log n)^{-2/3}n^{2/5}).
      \end{cases}
    \end{equation*}
  \end{enumerate}
\end{thm}

\subsection{Proof techniques and outline}\label{sec:techniques}

The techniques used in this paper are novel in comparison to the vast majority of papers on 
random graphs. Classical random graph results are usually proved with the help of probabilistic arguments
such as first and second moment methods, independence of random variables, or martingales. On the other hand, 
papers about random graphs on surfaces, e.g.~\cite{Chapuy2011-enumeration-graphs-on-surfaces,gimenez2009}, 
use singularity analysis of generating functions. In contrast, we combine various complementary methods
to prove our results.

The starting point of our proofs are \emph{constructive decompositions} of graphs, a method mostly used
in enumerative combinatorics. Every graph in
$\class{\general}_g(n,m)$ can be decomposed into its complex components and non-complex components,
which then can further be decomposed into smaller parts. The most important structures occurring in
this decomposition are the so-called \emph{core} and \emph{kernel} of the graph. The decomposition
is \emph{constructive} in the sense that every graph can be constructed in a unique way starting
from its kernel via its core and complex components (see \Cref{sec:decomposition}).

We interpret the aforementioned constructive decomposition in terms of \emph{combinatorial counting}, 
in other words, we represent the number of graphs in the class $\class{\general}_g(n,m)$ as a sum of 
subclasses that are involved in the decomposition. We proceed by determining the main contributions 
to the sum using a combinatorial variant of \emph{Laplace's method} from complex analysis, a technique 
to derive asymptotic estimates of integrals that depend on a parameter $n$ tending to infinity. To 
illustrate how we apply this approach, assume that we want to analyse a sum of the form
\begin{equation*}
  A(n) = \sum_{i\in I}B(i)C(n-i),
\end{equation*}
where $i$ is a parameter related to one of the substructures occurring in the constructive decomposition,
e.g.\ the order of the core, say. We rewrite $A(n)$ as
\begin{equation*}
  A(n) = \sum_{i\in I}\exp\br{f(i)}
\end{equation*}
with $f(i)=\log(B(i)C(n-i))$ and then estimate the exponent $f(i)$ in order to determine the \emph{main 
contribution} to $A(n)$ in the following sense. We determine a set $J\subset I$ so that the partial sum 
over all $i\in I\setminus J$ (the \emph{tail} of the sum) is of smaller order than the total sum (see 
\Cref{def:sums} for a formal definition). The probabilistic interpretation of this main contribution
is that $\general_g(n,m)$ \whp\ has its corresponding parameter $i$ in the set $J$. In our example, this
will tell us the `typical' order of the core of $\general_g(n,m)$.

The exact method how we estimate the value of the tail and compare it to the total value of the sum will
differ from case to case. In some cases, rough bounds provided by maximising techniques will suffice; in
other cases, we need better bounds, which we derive by using Chernoff bounds or by bounding the sums via
integrals. Systematic applications of these techniques enable us to derive the exact ranges of the main 
contributions. From the main contributions, we deduce the orders of components, component structure, and other structural 
properties of $\general_g(n,m)$ by applying both combinatorial methods (e.g.\ double counting) and 
probabilistic techniques (e.g.\ Markov's and Chebyshev's inequalities).

This paper is organised as follows. After presenting the necessary notation and definitions
in \Cref{sec:prelim}, we give an overview of the proof strategy in \Cref{sec:strategy}; in particular,
we derive the aforementioned representation of $\abs{\class{\general}_g(n,m)}$ as a sum. In 
\Cref{sec:kernelcorecomplex}, we determine the main contributions to this sum using the techniques
mentioned above. From these results, we derive structural properties of $\general_g(n,m)$ in
\Cref{sec:structure}. \Cref{sec:proofs:main,sec:proofs:lemmas} are devoted to the proofs of our
main results and of the auxiliary results, respectively. Finally, we discuss 
various open questions in \Cref{sec:discussion}.

\subsection{Related work}\label{sec:related}

The order of the largest component of the \ER\ random graph $G(n,m)$ at the time of the phase
transition has been extensively studied
\cite{Bollobas84,BollobasRiordan2012A,luczak1,luczakpittel,PittelWormald05}. Most of the results
have been proved using purely probabilistic arguments (e.g.~random walks, martingales), leading to 
even stronger results than the ones stated in \Cref{thm:ER}, e.g.\ about the limiting distribution 
of the order and size of the largest component~\cite{BehrCOKang10,coja-oghlan2014,BollobasRiordan2012A,BollobasRiordan2017}. 
In the case of $\general_g(n,m)$, the additional condition of the graph being embeddable on $\Sg$ makes it 
virtually impossible to use the same techniques in order to derive such strong results.

Comparing \Cref{thm:ER,main1}, the main differences appear when the giant component arises in the
\emph{supercritical regime}, that is, when $\edgesfirst\to\infty$. Firstly, the order of the giant
component is only about half as large in $\general_g(n,m)$ as it is in $G(n,m)$. Secondly, the
$i$-th largest component $\component_i$ for fixed $i\ge 2$ is much larger in $\general_g(n,m)$ than in
$G(n,m)$. These two differences are closely related for the following reason. In $G(n,m)$, the
number $n'$ of vertices and $m'$ of edges \emph{outside} the giant component are such that $m' =
(1+\edgesfirst'n^{-1/3})\frac{n'}{2}$ with $\edgesfirst'\to-\infty$ and thus, $G(n',m')$ only has
small components. In $\general_g(n,m)$, the smaller order of the giant component enforces $m'$ to
be in the \emph{critical regime}, where $\edgesfirst'\to c\in\R$, thus resulting in larger orders
for $\component_i$ with $i\ge 2$. Lastly, while each such $\component_i$ is a tree \whp\ for the \ER\ random graph,
it has a positive probability to be complex for $\general_g(n,m)$.

Planar graphs and graphs embeddable on $\Sg$ have been investigated separately for the `sparse'
regime $m \le n+o(n)$ \cite{kang2012} and for the `dense' regime $m = \lfloor\mu n\rfloor$ with
$\mu\in(1,3)$ \cite{Chapuy2011-enumeration-graphs-on-surfaces,gimenez2009}. From a random graph
point of view, in particular when the giant component is considered, the sparse regime is the more
interesting regime. In the sparse regime, Kang and {\L}uczak \cite{kang2012} supplied new
resourceful proof methods---some of which we apply in a somewhat similar fashion in this paper---combining
probabilistic and graph theoretic methods with techniques from enumerative and analytic
combinatorics. On the other hand, minor mistakes in \cite{kang2012} led to results that featured
order terms that claimed to be stronger than what has actually been proved. One contribution of
this paper is to correct and strengthen these results from \cite{kang2012}.

In the dense regime, Gim{\'e}nez and Noy \cite{gimenez2009} and Chapuy, Fusy, \Gimenez, Mohar, and
Noy~\cite{Chapuy2011-enumeration-graphs-on-surfaces} use techniques from analytic combinatorics to
prove limit laws for graphs embeddable on $\Sg$. The advantage of their techniques is that one
method can be applied to derive a range of various limit laws, e.g.\ on the number of components,
the order of the largest component, and the chromatic and list-chromatic number. On the other hand,
the techniques are limited to a) the class $\class{\general}_g(n)$ of $n$-vertex graphs embeddable 
on $\Sg$, in other words, graphs with $n$ vertices and an arbitrary number of edges, or b) the class
$\class{\general}_g(n,\lfloor\mu n\rfloor)$, where $\mu$ is a constant. A random graph chosen from 
the class $\class{\general}_g(n)$ is averaged over all graphs with an arbitrary number of edges and
thus not appropriate when we look at a specific range of $m$.\footnote{In fact, the properties of a
  random graph chosen from $\class{\general}_g(n)$ are dominated by the graphs whose edge density 
  is quite large, more precisely, when $\mu\approx 2.21$~\cite{Chapuy2011-enumeration-graphs-on-surfaces,gimenez2009}.}
On the other hand, the class $\class{\general}_g(n,\lfloor\mu n\rfloor)$ scales the number $m$ of edges
as a linear function in $n$, which is not fine enough in order to capture the changes that take place
within the critical windows, which have length $\Theta(n^{2/3})$ for \Cref{main1} and $\Theta(n^{3/5})$
for \Cref{main2}. In terms of critical behaviour these techniques are therefore not applicable.

\section{Preliminaries}\label{sec:prelim}

\subsection{Asymptotic notations}\label{sec:asymptotics}

By $\N$ we denote the set of non-negative integers.
In order to express orders of components in a random graph when $n$ tends to infinity, we use the
following notation. Recall that an event holds \emph{with high probability}, or \whp\ for short,
if it holds with probability tending to one as $n$ tends to infinity.

\begin{definition}\label{def:landau}
  Let $X=(X_n)_{n\in\N}$ be a sequence of random variables and let $f\colon\N\to\R_{\ge 0}$ be a function. For
  $c\in\R^+$ and $n\in\N$, consider the inequalities
  \begin{align}
    \abs{X_n} &\le c f(n),\label{landau:upper}\\
    \abs{X_n} &\ge c f(n).\label{landau:lower}
  \end{align}
  We say that
  \begin{enumerate}
  \item $X_n = O(f)$ \whp, if there exists $c\in\R^+$ such that \eqref{landau:upper} holds \whp;
  \item $X_n = o(f)$ \whp, if for every $c\in\R^+$, \eqref{landau:upper} holds \whp;
  \item $X_n = \Omega(f)$ \whp, if there exists $c\in\R^+$ such that \eqref{landau:lower} holds \whp;
  \item $X_n = \omega(f)$ \whp, if for every $c\in\R^+$, \eqref{landau:lower} holds \whp;
  \item $X_n = \Theta(f)$ \whp, if both $X_n = O(f)$ and $X_n = \Omega(f)$ \whp;
  \item $X_n = O_p(f)$, if for every $\delta>0$, there exist $c_\delta\in\R^+$ and $N_\delta\in\N$
    such that \eqref{landau:upper} holds for $c=c_\delta$ and $n\ge N_\delta$ with probability at
    least $1-\delta$;
  \item $X_n = \Theta_p(f)$, if for every $\delta>0$, there exist $c_\delta^+,c_\delta^-\in\R^+$ and
    $N_\delta\in\N$ such that for $n\ge N_\delta$ with probability at least $1-\delta$, both
    \eqref{landau:upper} holds for $c=c_\delta^+$ and \eqref{landau:lower} holds for $c=c_\delta^-$.
  \end{enumerate}
\end{definition}
\noindent
The special case of $X=O_p(1)$ is also known as $X$ being \emph{bounded in probability}.

\subsection{Graphs on surfaces}\label{sec:surfaces}

Given a graph $G$, we denote its vertex set and its edge set by $V(G)$ and $E(G)$, respectively.
All graphs in this paper are vertex-labelled, that is, $V(G)=[n]$ for some $n\in\N$. Let $g\in\N$ be fixed.
An \emph{embedding} of a graph $G$ on $\Sg$, the orientable surface of genus $g$, is a drawing of
$G$ on $\Sg$ without crossing edges. If $G$ has an embedding on $\Sg$, we call $G$
\emph{embeddable} on $\Sg$. Clearly, embeddability is monotone in $g$, i.e.\ every graph that is
embeddable on $\Sg$ is also embeddable on $\torus{g+1}$. By the \emph{genus} of a given graph $G$ 
we denote the smallest $g\in\N$ for which $G$ is embeddable on $\Sg$. Graphs with genus zero are 
also called \emph{planar}.

Let $\component$ be a connected graph embeddable on $\Sg$. We say that $\component$ is \emph{unicyclic} 
if it contains precisely one cycle and we call $\component$ \emph{complex} (also known as \emph{multicyclic}) if it contains at least two 
cycles; the latter is the case if and only if $\component$ has more edges than vertices. If $\component$ 
is complex, we call
\begin{equation*}
  \ex(\component):=\abs{E(\component)}-\abs{V(\component)}
\end{equation*}
the \emph{excess} of $\component$. For a non-connected graph $G$, we define $\ex(G)$ to be the sum of the
excesses of its complex components (and set $\ex(G)=0$ as a convention if $G$ has no
complex components). $G$ is called \emph{complex} if all its components are complex.

\subsection{Complex part, core, and kernel}\label{sec:compcorekernel}

Let $G$ be any graph. The union $\complex_G$ of all complex components of $G$ is called the
\emph{\comppart} of $G$. The \emph{core} $\core_G$ of $G$ is defined as the maximal subgraph of
minimum degree at least two of $\complex_G$. The core can also be obtained from the \comppart\ by
recursively deleting vertices of degree one (in an arbitrary order). Vice versa, the \comppart\
can be constructed from the core by attaching trees to the vertices of the core. Finally, the
\emph{kernel} $\kernel_G$ of $G$ is constructed from the core $\core_G$ by replacing all vertices
of degree two in the following way. Every maximal path $P$ in $\core_G$ consisting of vertices of
degree two is replaced by an edge between the vertices of degree at least three that are adjacent
to the end vertices of $P$. By this construction, loops and multiple edges can occur. Reversing the construction, 
the core arises from the kernel by subdividing edges.

It is important to note that $\kernel_G$ is non-empty as soon as $\complex_G$ is, because each
component of the complex graph $\complex_G$ contains a non-empty core with at least one vertex 
of degree at least three. Furthermore, $\kernel_G$ has minimum degree at least three and might 
contain loops and multiple edges. Observe that $G$ is embeddable on $\Sg$ if and only if 
$\kernel_G$ is. In particular, $G$ and $\kernel_G$ have the same genus. Also observe that 
$\ex(G) = \ex(\complex_G)$ by definition and $\ex(\kernel_G) = \ex(\core_G) = \ex(\complex_G)$, 
because subdividing edges and attaching trees changes the number of vertices and edges by the same 
amount.

Given a graph $G$ with $n$ vertices, we denote the number of vertices of the \comppart\
$\complex_G$, the core $\core_G$, and the kernel $\kernel_G$ by $n_\complex$, $n_\core$, and
$n_\kernel$, respectively. The number of edges of $\complex_G$, $\core_G$, and $\kernel_G$ 
satisfy
\begin{equation*}
  \abs{E(\complex_G)} = n_\complex+\ex(G),
  \qquad
  \abs{E(\core_G)} = n_\core+\ex(G),
  \qquad
  \abs{E(\kernel_G)} = n_\kernel+\ex(G).
\end{equation*}
The kernel has minimum degree at least three by definition and thus has at least $\frac32n_\kernel$
edges. A kernel is called \emph{cubic} if all its vertices have degree three; in that case, it has
precisely $\frac32n_\kernel$ edges. The \emph{deficiency} of $G$ is defined as
\begin{equation*}
  \deficiency(G) := 2\abs{E(\kernel_G)}-3n_\kernel = 2\ex(G) - n_\kernel.
\end{equation*}
Clearly, the deficiency is always non-negative and $\deficiency(G)=0$ if and only if the kernel
$\kernel_G$ is either empty or cubic. The definition of the excess and deficiency of a graph
immediately implies the following relation between the deficiency, the excess, and the number of
vertices and edges of the kernel.

\begin{lem}\label{kernel-ld}
  Given a graph $G$, the numbers $n_\kernel$ of vertices and $m_\kernel$ of edges in the kernel $\kernel_G$ of $G$ are
  \begin{equation*}
    n_\kernel = 2\ex(G) - \deficiency(G)
    \qquad\text{and}\qquad
    m_\kernel = 3\ex(G) - \deficiency(G).
  \end{equation*}
\end{lem}

\subsection{Useful bounds}\label{sec:bounds}

We will
frequently use the following widely known formulas.
\begin{align}
  1+x&=\exp\br{x-\frac{x^2}{2}+\frac{x^3}{3}+O(x^4)}&
  &\text{if }x=o(1),&\label{explog}\\
  1+x&\le\exp\br{x},&&&\label{explogu}\\
  1+x&\ge\exp\br{x-\frac{x^2}{2}}&&\text{if }x\ge0,&\label{explogl}
\end{align}
To derive bounds for the factorial $n!$ and the falling factorial $(k)_i := k!/(k-i)!$ we shall use the inequalities
\begin{align}
 \sqrt{2\pi n}\br{\frac{n}{e}}^n&\leq n!\leq e\sqrt{n}\br{\frac{n}{e}}^n,\label{stirling:bound}\\
 k^i\exp\br{-\frac{i^2}{2(k-i)}}&\leq (k)_i\leq k^i\exp\br{-\frac{i(i-1)}{2k}}.\label{falling:bound}
\end{align}
For $1\le k\le n-1$ we will also use refined bounds for the binomial coefficient obtained by 
applying \eqref{stirling:bound} thrice.
\begin{align}
 \frac{\sqrt{2\pi}n^{n+1/2}}{e^2 k^{k+1/2}(n-k)^{n-k+1/2}}&\leq \binom{n}{k}\leq\frac{en^{n+1/2}}{2\pi k^{k+1/2}(n-k)^{n-k+1/2}}\,.\label{binomial:refined}
\end{align}
We shall also use the inequality
\begin{equation}\label{bound:frac}
  \frac{1}{a+b} \ge \frac{1}{a} - \frac{b}{a^2} \quad
  \text{if }a\not=0, a+b>0.
\end{equation}

Finally, we need some well known inequalities from probability theory. Given a random variable $X$, we
denote by $\expec{X}$ its expectation. The \emph{variance} of $X$ is then defined as $$\sigma^2:=\expec{\br{X-\expec{X}}^2} = \expec{X^2}-\expec{X}^2.$$
For a non-negative random variable $X$ and any $t>0$, Markov's inequality states that
\begin{equation}\label{eq:Markov}
  \prob{X\geq t}\leq \frac{\expec{X}}{t}.
\end{equation}
A stronger bound---which additionally holds for arbitrary random variables---is provided by Chebyshev's inequality.
For any random variable $X$ and any $t>0$, we have
\begin{equation}\label{eq:Chebyshev}
  \prob{\abs{X-\expec{X}}\geq t}\leq \frac{\sigma^2}{t^2}.
\end{equation}
In terms of Chernoff bounds, we shall need the two special cases of normal distributions
and binomial distributions. For a normally distributed random variable $X$, we have, for any given $t>0$,
\begin{equation}\label{eq:Chernoff:N}
  \prob{\Big.\!\abs{\big.X-\expec{X}}\ge t} \le 2\exp\br{-\frac{t^2}{2\sigma^2}}.
\end{equation}
If $X$ is a binomially distributed random variable, then
\begin{equation}\label{eq:Chernoff:B}
 \prob{\Big.\!\abs{\big.X-\expec{X}} \ge t} \le 2\exp\br{-\frac{t^2}{2\br{\expec{X}+\frac{t}{3}}}}.
\end{equation}

\section{Proof strategy}\label{sec:strategy}

\subsection{Decomposition and construction}\label{sec:decomposition}

Throughout the paper, let $g\in\N$ be fixed.
We have seen in \Cref{sec:compcorekernel} that any graph that is embeddable on $\Sg$ can be decomposed
into a) its \comppart\ and b) trees and unicyclic components. The \comppart\ can then further be
decomposed so as to obtain the core and the kernel. Vice versa, we can construct a graph on $\Sg$
by performing the reverse constructions.

\begin{construction*}
  The following steps construct every graph embeddable on $\Sg$.
  \renewcommand{\theenumi}{(\textrm{C}\arabic{enumi})}
  \begin{enumerate}
  \item\label{construction:kernel} Pick a kernel, i.e.\ a multigraph with minimum degree at least three
    that is embeddable on $\Sg$;
  \item\label{construction:kerneltocore} subdivide the edges of the kernel to obtain a core;
  \item\label{construction:coretocomplex} to every vertex $v$ of the core, attach a rooted tree $T_v$ (possibly only
    consisting of one vertex) by identifying $v$ with the root of $T_v$, so as to obtain a
    complex graph;
  \item\label{construction:complextogeneral} add trees and unicyclic components to obtain a general
    graph embeddable on $\Sg$.
  \end{enumerate}
\end{construction*}
\renewcommand{\theenumi}{(\roman{enumi})}

To avoid overcounting in \ref{construction:kerneltocore} if the kernel has loops or multiple
edges, multigraphs will always be weighted by the \emph{compensation factor} introduced by Janson,
Knuth, \Luczak, and Pittel~\cite{jansonea}, which is defined as follows. Given a multigraph $M$ and an integer $i\ge
1$, denote by $e_i(M)$ the number of (unordered) pairs $\{u,v\}$ of vertices for which there are
exactly $i$ edges between $u$ and $v$. Analogously, let $\ell_i(M)$ denote the number of vertices
$x$ for which there are precisely $i$ loops at $x$. Finally, let $\ell(M)=\sum_i i\ell_i(M)$ be the
number of loops of $M$. The compensation factor of $M$ is defined to be
\begin{equation}\label{eq:compensationfactor}
  \cf(M):=2^{-\ell(M)}\prod_{i=1}^{\infty}(i!)^{-e_i(M)-\ell_i(M)}.
\end{equation}
In \ref{construction:kerneltocore}, the compensation factor enables us
to distinguish multiple edges and loops at the same vertex (because of the factors $1/i!$) as well
as the different orientations of loops (because of the factor $2^{-\ell(M)}$). This fact ensures
that there is no overcounting in \ref{construction:kerneltocore}. Indeed, if a core $\core$
has kernel $\kernel$, then $\core$ can be constructed from $\kernel$ by subdividing edges in
precisely $\frac{1}{\cf(\kernel)}$ different ways; thus, assigning weight $\cf(\kernel)$ to
$\kernel$ prevents overcounting.

We denote by
\begin{itemize}
\item $\class{\general}_g$ the class of all graphs embeddable on $\Sg$;
\item $\class{\complex}_g$ the class of all \comppart s of graphs in $\class{\general}_g$;
\item $\class{\core}_g$ the class of all cores of graphs in $\class{\general}_g$;
\item $\class{\kernel}_g$ the class of all kernels of graphs in $\class{\general}_g$;
\item $\class{\outside}$ the class of all graphs without complex components.
\end{itemize}
In other words, $\class{\complex}_g$ is the class of all complex graphs embeddable on $\Sg$;
$\class{\core}_g$ consists of all complex graphs embeddable on $\Sg$ with minimum degree at least two;
and $\class{\kernel}_g$ comprises all (weighted) multigraphs embeddable on $\Sg$ with minimum degree at least three. 
The empty graph lies in all the classes above by convention.

If $n,m\in\N$ are fixed, we write $\class{\general}_g(n,m)$ for the subclass of
$\class{\general}_g$ containing all graphs with exactly $n$ vertices and $m$ edges. By
$\general_g(n,m)$ we denote a graph chosen uniformly at random from all graphs in
$\class{\general}_g(n,m)$. We use the corresponding notation also for the other classes defined
above.

The construction of graphs in $\class{\general}_g$ from their kernel via the core and \comppart\
as described in \ref{construction:kernel}--\ref{construction:complextogeneral} can be translated to relations between the numbers of graphs
in the previously defined classes. Starting from $\class{\general}_g(n,m)$,
\ref{construction:complextogeneral} immediately gives rise to the identity
\begin{equation}\label{eq:general}
  \abs{\class{\general}_g(n,m)} =
  \sum_{n_\complex,l}\binom{n}{n_\complex}\abs{\class{\complex}_g(n_\complex,n_\complex+l)}\cdot\abs{\class{\outside}(n_\outside,m_\outside)},
\end{equation}
where $n_\outside = n-n_\complex$ and $m_\outside = m-n_\complex-l$. Indeed, for each fixed number
$n_\complex$ of vertices in the \comppart\ and each fixed excess $l$
\begin{itemize}
\item the binomial coefficient counts the possibilities which vertices lie in the \comppart,
\item $\abs{\class{\complex}_g(n_\complex,n_\complex+l)}$ counts the \comppart s with $n_\complex$
  vertices and $n_\complex+l$ edges, and
\item $\abs{\class{\outside}(n_\outside,m_\outside)}$ counts all possible arrangements of non-complex
  components.
\end{itemize}

If $\abs{\class{\complex}_g(n_\complex,n_\complex+l)}$ and $\abs{\class{\outside}(n_\outside,m_\outside)}$
are known, then we can use \eqref{eq:general} to determine $\abs{\class{\general}_g(n,m)}$.
Determining $\abs{\class{\complex}_g(n_\complex,n_\complex+l)}$ turns out to be quite a challenging task,
to which we devote a substantial part of this paper. The number $\abs{\class{\outside}(n_\outside,m_\outside)}$,
on the other hand, can be determined using known results.

\subsection{Graphs without complex components}\label{sec:outside}

The class $\class{\outside}$ of graphs without complex components (i.e.\ each component is either a tree or unicyclic)
has been studied by Britikov~\cite{britikov2} and by Janson,
Knuth, \Luczak, and Pittel~\cite{jansonea}, who determined the number of graphs in $\class{\outside}(n,m)$ for different regimes of $m$.

\begin{lem}\label{britikov}
  Let $m=(1+\edgesfirst n^{-1/3})\frac{n}{2}$ with $\edgesfirst=\edgesfirst(n)<n^{1/3}$ and
  let $\rho(n,m)$ be such that
  \begin{equation*}
    \abs{\class{\outside}(n,m)}=\binom{\binom{n}{2}}{m}\rho\br{n,m}.
  \end{equation*}
  There exists a constant $c>0$ such that for
  \begin{equation*}
   f(n,m)=c\br{\frac{2}{e}}^{2m-n}\frac{m^{m+1/2}n^{n-2m+1/2}}{(n-m)^{n-m+1/2}},
  \end{equation*}
  we have
  \begin{enumerate}
  \item\label{britikov:sub}
    $\rho\br{n,m}=1+o(1)$, if $\edgesfirst\ra-\infty$;
  \item\label{britikov:crit}
    for each $a\in\R$, there exists a constant $b=b(a)>0$ such that $\rho\br{n,m}\ge b$ whenever $\edgesfirst\le a$;
  \item\label{britikov:sup}
    $\rho\br{n,m}\le n^{-1/2}f(n,m)$ if $\edgesfirst>0$ and $\edgesfirst=o(n^{1/12})$;
  \item\label{britikov:sup2}
    $\rho\br{n,m}\le f(n,m)$ if $\edgesfirst>0$.
  \end{enumerate}
\end{lem}

\Cref{britikov}\ref{britikov:sub}, \ref{britikov:crit}, and \ref{britikov:sup} are proven in 
\cite{britikov2} and \cite{jansonea}, but \ref{britikov:sup2} is a slight extension of the results in \cite{jansonea} 
which we prove in \Cref{sec:proofs:lemmas} along the following lines. Inspired by the proof of \ref{britikov:sup}
in \cite{jansonea}, we bound $\rho(n,m)$ by a contour integral and prove that this integral has value at most $f(n,m)$
for all $\edgesfirst>0$.

Clearly, every graph in $\class{\outside}$ is planar and thus also embeddable on $\Sg$. This fact,
together with \Cref{britikov} and \Cref{thm:ER}\ref{Gnm:sub} and \ref{Gnm:crit} will be enough to
prove \Cref{main1}\ref{firstsub} and \ref{firstcrit}. For all other regimes, \Cref{britikov} will
provide a very useful way to bound the number $\abs{\class{\outside}(n,m)}$ in \eqref{eq:general}.

\subsection{Complex parts}\label{sec:complex}

For the number $\abs{\class{\complex}_g(n_\complex,n_\complex+l)}$,
we analyse \ref{construction:kernel}--\ref{construction:coretocomplex}
in order to derive an identity similar to \eqref{eq:general}. Firstly, we need to sum over
all possible numbers $n_\core$ of vertices in the core; the number of edges in the core is then
given by $n_\core+l$. For fixed $n_\complex$, $n_\core$, and $l$, we have
\begin{itemize}
\item $\binom{n_\complex}{n_\core}$ choices for which vertices of the \comppart\ lie in the
  core,
\item $\abs{\class{\core}_g(n_\core,n_\core+l)}$ ways to choose a core, and
\item $n_\core n_\complex^{n_\complex-n_\core-1}$ possibilities to attach $n_\core$ rooted trees with
  $n_\complex$ vertices in total to the vertices of the core.
\end{itemize}
By \ref{construction:coretocomplex}, we thus deduce that
\begin{equation}\label{eq:complexcore}
  \abs{\class{\complex}_g(n_\complex,n_\complex+l)} = 
  \sum_{n_\core}\binom{n_\complex}{n_\core}\abs{\class{\core}_g(n_\core,n_\core+l)}n_\core n_\complex^{n_\complex-n_\core-1}\,.
\end{equation}

In order to determine $\abs{\class{\core}_g(n_\core,n_\core+l)}$, recall that by \Cref{kernel-ld},
the number of vertices and edges in the kernel depend only on the excess and deficiency of the
graph. Thus, we choose the deficiency $d$ as the summation index. The number of ways to construct a
core from a kernel according to \ref{construction:kerneltocore} cannot
be described in an easy fashion like the constructions in~\ref{construction:coretocomplex}
and~\ref{construction:complextogeneral}. We will investigate this construction step in more detail
in \Cref{lem:SigmaD}. For a kernel $\kernel \in \class{\kernel}_g(2l-d,3l-d)$,
consider the number of different ways to subdivide its edges that result in a core with $n_\core$
vertices and $n_\core+l$ edges. Denote by $\varphi_{n_\core,l,d}$ the average of this number, taken over
all kernels in $\class{\kernel}_g(2l-d,3l-d)$. With this notation, we deduce from \ref{construction:kernel} and
\ref{construction:kerneltocore} that
\begin{equation}\label{eq:core}
  \abs{\class{\core}_g(n_\core,n_\core+l)} =
  \sum_d\binom{n_\core}{2l-d}\abs{\class{\kernel}_g(2l-d,3l-d)} \varphi_{n_\core,l,d}.
\end{equation}
Recall that the multigraphs in $\class{\kernel}_g$ are weighted. Accordingly,
$\abs{\class{\kernel}_g(2l-d,3l-d)}$ does not denote the \emph{number} of these multigraphs, but the
\emph{sum of their weights}.

\subsection{Analysing the sums}\label{sec:maincontributions}

In each of \eqref{eq:general}, \eqref{eq:complexcore}, and \eqref{eq:core}, we may assume that the
parameters $n_\complex,n_\core,l,d$ of the sums only take those values for which the summands are
non-zero.
\begin{definition}\label{def:admissible}
  We call values for a parameter (or a set of parameters) \emph{admissible}, if there exists at least one graph 
  satisfying these values for the corresponding parameters.
\end{definition}
The definition of the parameters, together with \Cref{kernel-ld}, directly yield the following necessary
conditions for admissibility.
\renewcommand{\theenumi}{(\textrm{A}\arabic{enumi})}
\begin{enumerate}
\item $0\leq n_\complex\leq n$;
\item $0\leq n_\core\leq n_\complex$;
\item $0\leq l\leq m-n_\complex$;
\item $l=0$ if and only if $n_\complex=0$;
\item\label{ineq:euler}
  $l\leq 2n_\core + 6(g-1)$;
\item $0\leq d\leq 2l$.
\end{enumerate}
\renewcommand{\theenumi}{(\roman{enumi})}
Inequality \ref{ineq:euler} is due to Euler's formula applied to the core. These bounds will frequently be used; 
if we use other bounds, we will state them explicitly.

On the first glance, the sole application of \eqref{eq:general}, \eqref{eq:complexcore}, and \eqref{eq:core}
seems to be to determine the number of graphs with given numbers of vertices and edges in the classes
$\class{\general}_g$, $\class{\complex}_g$, and $\class{\core}_g$. However, we shall use these sums to
derive \emph{typical structural properties} of graphs chosen uniformly at random from one of these classes.

Our plan to derive such properties from the sums \eqref{eq:general}, \eqref{eq:complexcore}, and \eqref{eq:core} 
is as follows. Once we have determined the values $\abs{\class{\kernel}_g(2l-d,3l-d)}$
and $\varphi_{n_\core,l,d}$, we consider the parameters $n_\complex,n_\core,l,d$ of the sums
one after another. For each parameter $i$, we seek to determine which range for $i$ provides the `most important' 
summands. To make this more precise, let us introduce the following notation.

\begin{definition}\label{def:sums}
  For every $n\in\N$, let $I(n),I_0(n)\subset\N$ be finite index sets with $I_0(n)\subseteq I(n)$.
  For each $i\in I(n)$, let $A_i(n)\ge 0$. We say that the \emph{main contribution} to the sum
  \begin{equation*}
    \sum_{i\in I(n)}A_i(n)
  \end{equation*}
   is provided by $i\in I_0(n)$ if
  \begin{equation*}
    \sum_{i\in I(n)\setminus I_0(n)}A_i(n) = o\br{\sum_{i\in I(n)}A_i(n)},
  \end{equation*}
  where $n\to\infty$. The sum over $i\in I(n)\setminus I_0(n)$ is then called the \emph{tail} of
  $\sum A_i(n)$.
\end{definition}

We shall determine index sets $I_\complex(n)$, $I_\core(n)$, $I_l(n)$, $I_d(n)$ so that the main
contributions to the sums in \eqref{eq:general}, \eqref{eq:complexcore}, and \eqref{eq:core} are
provided by $n_\complex\in I_\complex(n)$, $n_\core\in I_\core(n)$, $l\in I_l(n)$, and $d\in
I_d(n)$, respectively. This will yield statements about the size of these values in the following way.
For fixed $m=m(n)$, the index set $I_\core(n)$, for example, will be of the type $[c_1f(n),c_2f(n)]$ 
for certain constants $0<c_1<c_2$ and a certain function $f=f(n)$. This implies
that if $G=\general_g(n,m)$, then \whp\ $n_\core\in I_\core(n)$ and thus $n_\core=\Theta(f)$.

The main challenge is to find the `optimal' intervals 
$I_\complex(n)$, $I_\core(n)$, $I_l(n)$, $I_d(n)$ in view of \Cref{def:sums} in the
sense that they should be a) large enough so as to provide the main contribution and b)
as small as possible so as to yield stronger concentration results. To achieve these two
antipodal goals is a difficult task whose solution will differ from case to case.
In order to prove that a given interval indeed provides the main contribution to a sum, 
we bound the tail of the sum using various complementary methods including maximising 
techniques (e.g.\ \Cref{lem:SigmaC,lem:SigmaQ,nQandl:lower,lem:caseB}), Chernoff bounds (\Cref{lem:SigmaC,lem:SigmaD}), 
and approximations by integrals (\Cref{lem:SigmaQ:small,nQandl}).

Determining the main contributions to \eqref{eq:general}, \eqref{eq:complexcore}, and \eqref{eq:core}
will yield structural statements like the typical order of the \comppart, the core, and the kernel
of $G=\general_g(n,m)$. In order to derive the component structure of $G$, we further apply combinatorial techniques
like double counting (e.g.\ \Cref{main2} and \Cref{kernelpumping}) and probabilistic methods such as Markov's 
and Chebyshev's inequalities (\Cref{generalstructure}).

\section{Kernels, cores, and \comppart s}\label{sec:kernelcorecomplex}

For the remainder of the paper, let $n,m,n_\complex,n_\core,l,d\in\N$ be such that $m=m(n)\le(1+o(1))n$
and such that $n_\complex$, $n_\core$, $l$, and $d$ are admissible (in terms of \Cref{def:admissible}).
Furthermore, set $n_\outside=n-n_\complex$ and $m_\outside=m-n_\complex-l$.

The aim of this section is to determine the main contributions (in the sense of \Cref{def:sums}) to the sums in \eqref{eq:general},
\eqref{eq:complexcore}, and \eqref{eq:core}. In other words, we derive the typical orders of the
\comppart\ and the core of $G=\general_g(n,m)$, as well as the excess and the deficiency of $G$. These orders
will be the main ingredients for the proofs of Theorems~\ref{main1}--\ref{main3}.
For all results in this section, we defer the proofs to \Cref{sec:proofs:lemmas}.

\subsection{Kernels}\label{sec:kernels}

Throughout this section, we assume $l\ge1$.
As a basis of our analysis of \eqref{eq:general}, \eqref{eq:complexcore}, and \eqref{eq:core}, we
first have to determine the sum $\abs{\class{\kernel}_g(2l-d,3l-d)}$
of weights of the multigraphs in $\class{\kernel}_g(2l-d,3l-d)$. 
We start with the case when the kernel is cubic (or equivalently, $d=0$). The number of cubic kernels was determined in~\cite{enum} by Fang and 
the authors of the present paper.

\begin{thm}[\cite{enum}]\label{cubic:enum}
  The number of cubic multigraphs with $2l$ vertices and $3l$ edges embeddable on $\Sg$, weighted by
  their compensation factor, is given by
  \begin{equation*}
    \abs{\class{\kernel}_g(2l,3l)}=\br{1+O\br{l^{-1/4}}}e_g l^{5g/2-7/2}\gamma_{\kernel}^{2l}(2l)!\,,
  \end{equation*}
  where $\gamma_{\kernel}=\frac{79^{3/4}}{54^{1/2}}\approx 3.606$ and $e_g>0$ is a constant depending
  only on $g$.
\end{thm}

The number of \emph{connected} cubic kernels will be of interest as well.

\begin{thm}[\cite{enum}]\label{cubic:connected}
  The number of \emph{connected} multigraphs in $\class{\kernel}_g(2l,3l)$, weighted by
  their compensation factor, is
  \begin{equation*}
    \br{1+O\br{l^{-1/4}}}c_g l^{5g/2-7/2}\gamma_{\kernel}^{2l}(2l)!\,,
  \end{equation*}
  where $\gamma_\kernel$ is as in \Cref{cubic:enum} and $c_g>0$ is a constant depending only on $g$.
\end{thm}

In particular, \Cref{cubic:enum,cubic:connected} imply that $\kernel_g(2l,3l)$ is connected with
probability tending to $\frac{c_g}{e_g}>0$; in other words, the probability that a random cubic kernel
is connected is bounded away from zero.

Before we consider kernels with non-zero deficiency, we shall look at the structure of cubic kernels.
We aim to find the giant component of $\general_g(n,m)$ and prove that it is complex, hence finding the giant
component of the kernel would be a basis for a complex giant component in $\general_g(n,m)$.
Moreover, we would like this giant component to have genus $g$. The following result from~\cite{enum}
provides us with a component of genus $g$ in a cubic kernel.

\begin{lem}[\cite{enum}]\label{cubic:structure}
  If $g\ge 1$, then $\kernel_g(2l,3l)$ \whp\ has one component of genus $g$ and all its other
  components are planar.
\end{lem}

Intuitively, the non-planar component provided by \Cref{cubic:structure} should be the largest
component of the kernel, ideally even large enough to be the giant component. The following result shows that this
component indeed covers almost all vertices in the kernel.

\begin{lem}\label{cubic:lc-size}
  Let $g\ge 1$. Denote by $\pl(G)$ the subgraph of $G=\kernel_g(2l,3l)$ consisting of all planar
  components. Then $\abs{\pl(G)} = O_p(1)$. Furthermore, $\abs{\pl(G)}$ is even and there exist constants
  $c^+,c^-\in\R^+$ such that for every fixed integer $i\ge 1$ and sufficiently large $l$,
  \begin{equation}\label{cubic:planar:eq}
    c^-i^{-7/2}\br{1-\frac{i}{l}}^{5g/2-7/2} \le
    \prob{\abs{\pl(G)}=2i} \le
    c^+i^{-7/2}\br{1-\frac{i}{l}}^{5g/2-7/2}.
  \end{equation}
\end{lem}

For the case $g=0$, \cite[Lemma 2]{kang2012} provides an analogous statement to \eqref{cubic:planar:eq}
for the number of vertices outside the giant component of $\kernel_0(2l,3l)$.

Let us now look at general (not necessarily cubic) kernels. For such
kernels, we are not able to give a precise formula for their number, but we can bound their number
by comparing them to cubic kernels via a double counting argument.

\begin{lem}\label{kernelpumping}
  Let $k\in\N$ be fixed. For $\kernel\in\class{\kernel}_g$, denote by
  \begin{enumerate}
  \item $\mathcal{P}_1$ the property that $\kernel$ has precisely $k$ components;
  \item $\mathcal{P}_2$ the property that, if $g\ge1$, then each component of $\kernel$ has genus strictly 
    smaller than $g$.
  \end{enumerate}
  For $i=1,2$, denote by $\class{\kernel}_g(n_\kernel,m_\kernel;\mathcal{P}_i)$ the subclass of
  $\class{\kernel}_g(n_\kernel,m_\kernel)$ of kernels that have property $\mathcal{P}_i$. Then
  \begin{equation*}
    \frac{\abs{\class{\kernel}_g(2l-d,3l-d)}}{\abs{\class{\kernel}_g(2l,3l)}} \le \frac{6^d}{d!}
    \quad\text{and}\quad
    \frac{\abs{\class{\kernel}_g(2l-d,3l-d;\mathcal{P}_i)}}{\abs{\class{\kernel}_g(2l,3l;\mathcal{P}_i)}} \le \frac{6^d}{d!}
    \text{ for }i=1,2.
  \end{equation*}
  If in addition $d\le\frac{2}{7}l$, then also
  \begin{equation*}
    \frac{\abs{\class{\kernel}_g(2l-d,3l-d)}}{\abs{\class{\kernel}_g(2l,3l)}} \ge \frac{1}{216^dd!}
    \quad\text{and}\quad
    \frac{\abs{\class{\kernel}_g(2l-d,3l-d;\mathcal{P}_i)}}{\abs{\class{\kernel}_g(2l,3l;\mathcal{P}_i)}} \ge \frac{1}{216^dd!}
    \text{ for }i=1,2.
  \end{equation*}
\end{lem}

\Cref{kernelpumping} has two main applications. On one hand, together with \Cref{cubic:enum},
\Cref{kernelpumping} provides a way to bound the value $\abs{\class{\kernel}_g(2l-d,3l-d)}$ in
\eqref{eq:core}. On the other hand, \Cref{kernelpumping} will also
enable us to extend the structural results from \Cref{cubic:structure,cubic:lc-size} to kernels with
a fixed constant deficiency $d$ (see \Cref{generalstructure}).

\subsection{Core and deficiency}\label{sec:coredeficiency}

We first determine the main contributions to the sums in \eqref{eq:complexcore} and \eqref{eq:core}.
By definition, $\abs{\class{\complex}_g(0,0)}=1$. Thus, throughout this section we will assume that
both $n_\complex\ge1$ and $l\ge1$ (recall that $l=0$ if and only if $n_\complex=0$).
Observe that \eqref{eq:complexcore}, \eqref{eq:core}, and the identity
\begin{equation*}
  \binom{n_\complex}{n_\core}\binom{n_\core}{2l-d} =
  \frac{(n_\complex)_{n_\core}}{(2l-d)!(n_\core-2l+d)!}
\end{equation*}
imply that
\begin{equation}\label{eq:complex}
  \abs{\class{\complex}_g(n_\complex,n_\complex+l)} = 
  \sum_{n_\core,d}\frac{(n_\complex)_{n_\core}  \abs{\class{\kernel}_g(2l-d,3l-d)}  \varphi_{n_\core,l,d}\,  n_\core n_\complex^{n_\complex-n_\core-1}}{(2l-d)!(n_\core-2l+d)!}\,.
\end{equation}

The factor $\abs{\class{\kernel}_g(2l-d,3l-d)}$ in \eqref{eq:complex} can be bounded using \Cref{cubic:enum} 
and \Cref{kernelpumping}. The term $\varphi_{n_\core,l,d}$, however, is still unknown. Recall that this value
denotes the average number, over all $\kernel\in\class{\kernel}_g(2l-d,3l-d)$, of
different ways to subdivide the edges of $\kernel$ that result in a core with $n_\core$ vertices
and $n_\core+l$ edges.

\begin{lem}\label{binsandballs}
  There exists a function $\nu=\nu(n_\core,l,d)$ such that
  \begin{equation*}
    \varphi_{n_\core,l,d} = (n_\core-2l+d)!\binom{n_\core+\nu l -1}{3l-d-1}
  \end{equation*}
  and $-5\le\nu\le 1$.
\end{lem}

Let us now determine the value of the sum in \eqref{eq:complex} over $n_\core$, as well as its main contribution. 
To this end, we apply \Cref{kernelpumping,binsandballs} to \eqref{eq:complex}, gather all parts
of the equation that depend on $n_\core$, and denote the sum over these values by $\Sigma_\core$.

  \begin{lem}\label{complex:bounds}
    There exists a function $\tau=\tau(l,d)$ such that
    \begin{enumerate}
    \item $\frac{1}{216}\le\tau\leq6$ for all $0\le d\le \left\lfloor\frac{2l}{7}\right\rfloor$;
    \item $0\leq\tau\leq6$ for all $\left\lfloor\frac{2l}{7}\right\rfloor < d\leq2l$;
    \end{enumerate}
    and
    \begin{align}\label{complex:step1}
      \abs{\class{\complex}_g(n_\complex,n_\complex+l)}=n_\complex^{n_\complex-1}\frac{\abs{\class{\kernel}_g(2l,3l)}}{(2l)!}&\sum_{d=0}^{2l}\binom{2l}{d}\frac{\tau^d}{(3l-d-1)!}\Sigma_\core,
    \end{align}
    where
    \begin{equation}\label{sum:core}
      \Sigma_\core=\Sigma_\core(n_\complex,l,d):=\sum_{n_\core}\frac{(n_\complex)_{n_\core}}{n_\complex^{n_\core}}n_\core(n_\core+\nu l-1)_{3l-d-1}.
    \end{equation}
  \end{lem}

  The strategy to determine the main contribution to $\Sigma_\core$ is roughly
  as follows. Using inequalities from \Cref{sec:bounds}, we bound
  $\Sigma_\core(n_\complex,l,d)$ from above by a sum of the type
  \begin{equation*}
    \sum_{n_\core}\exp\br{A(n_\complex,n_\core,l,d)}.
  \end{equation*}
  The derivative of $A(n_\complex,n_\core,l,d)$ with respect to $n_\core$ will show to have a zero at $n_\core =
  (1+o(1))\overline{n}_\core$, where
  \begin{equation*}\label{core:max}
    \ol{n}_\core = \sqrt{n_\complex(3l-d)}.
  \end{equation*}
  We then substitute $n_\core = \overline{n}_\core + r$ and prove that the resulting sum---up to
  a scaling factor---corresponds to a normally distributed random variable to which the Chernoff
  bound \eqref{eq:Chernoff:N} applies. Finally, for $n_\core$ from the range of the main 
  contribution to the upper bound, we derive a similar lower bound, which will enable us to derive 
  the main contribution to $\Sigma_\core$.
  
  \begin{lem}\label{lem:SigmaC}
    Let $f_\core=f_\core(n_\complex,l,d)$ be such that
    \begin{equation}\label{eq:Sigmacore}
      \Sigma_\core(n_\complex,l,d) = \sqrt{n_\complex}\br{\frac{n_\complex(3l-d)}{e}}^{(3l-d)/2}\exp\br{f_\core}.
    \end{equation}
    \begin{enumerate}
    \item\label{ncore:upper}
      There exist constants $a_\core^+,b_\core^+\in\R$ such that
      \begin{equation*}
        f_\core \le a_\core^+ + b_\core^+\sqrt{\frac{l^3}{n_\complex}}.
      \end{equation*}
    \item\label{ncore:lower}
      For every function $\epsilon(n_\complex)=o(1)$, there exist constants $N_\complex\in\N$, $a_\core^-,b_\core^-\in\R$ such that
      whenever $n_\complex\ge N_\complex$ and $\frac{7}{2}d\le l\le \epsilon n_\complex$, then
      \begin{equation*}
        f_\core \ge a_\core^- + b_\core^-\sqrt{\frac{l^3}{n_\complex}}.
      \end{equation*}
    \item\label{ncore:main}
      For every $0<\delta<\frac12$, whenever
      $n_\complex,l\to\infty$ and $\frac{7}{2}d\le l\le \epsilon n_\complex$, where $\epsilon=\epsilon(n_\complex)=o(1)$ is given, the main contribution to
      $\Sigma_\core$ is provided by
      \begin{equation*}
        n_\core \in I_\core^\delta(n_\complex,l,d) := \left\{k\in\N \;\big\vert\; \abs{k-\ol{n}_\core} < \delta\ol{n}_\core\right\}.
      \end{equation*}
    \end{enumerate}
  \end{lem}

  Our next aim is to analyse the sum over $d$ in \eqref{complex:step1}. To this end, observe that for
  \begin{equation}\label{complex:sum:def}
    \Sigma_d = \Sigma_d(n_\complex,l) := \sum_d\binom{2l}{d}\frac{(3l-d)^{(3l-d+2)/2}e^{d/2}\tau^d}{(3l-d)!n_\complex^{d/2}}\exp(f_\core),
  \end{equation}
  \eqref{complex:step1} and \eqref{eq:Sigmacore} yield
    \begin{equation}\label{complex:step2}
      \abs{\class{\complex}_g(n_\complex,n_\complex+l)}=\frac{n_\complex^{n_\complex+3l/2-1/2}\abs{\class{\kernel}_g(2l,3l)}}{e^{3l/2}(2l)!}\Sigma_d.
    \end{equation}

We determine the value of $\Sigma_d$, as well as its main contribution, in a similar fashion as for $\Sigma_\core$.

\begin{lem}\label{lem:SigmaD}
  Let $f_d=f_d(n_\complex,l)$ be such that
  \begin{equation}\label{eq:Sigmad}
    \Sigma_d = (3l)^{-(3l-1)/2}e^{3l}\exp\br{f_d}.
  \end{equation}
  \begin{enumerate}
  \item\label{deficiency:upper}
    There exist constants $a_d^+\in\R$ and $b_d^+\in\R$ such that
    \begin{equation*}
      f_d\leq a_d^+ + b_d^+\sqrt{\frac{l^3}{n_\complex}}.
    \end{equation*} 
  \item\label{deficiency:lower}
    For every function $\epsilon(n_\complex)=o(1)$, there exist constants
    $N_\complex\in\N$ and $a_d^-,b_d^-\in\R$ such that
    \begin{equation*}
      f_d\geq a_d^- + b_d^-\sqrt{\frac{l^3}{n_\complex}},
    \end{equation*}
    whenever $n_\complex\ge N_\complex$ and $l\le \epsilon  n_\complex$.
  \item\label{deficiency:maincontribution}
    There exists a constant $\beta_d^+\in\R^+$ such that for $n_\complex,l\to\infty$ and $l=o(n_\complex)$, the main contribution to
    $\Sigma_d$ is provided by
      \begin{enumerate}
      \item $d\in I_d(n_\complex,l) := \{0\}$ if $l=o(n_\complex^{1/3})$;
      \item $d\in I^h_d(n_\complex,l) := \{ k\in\N \mid k\le h(n_\complex)\}$ for every fixed function $h=h(n_\complex)=\omega(1)$ if
        $l=\Theta(n_\complex^{1/3})$;
      \item $d\in I_d(n_\complex,l) := \left\{ k\in\N \mid k \le \beta_d^+\sqrt{\frac{l^3}{n_\complex}}\right\}$
         if $l=\omega(n_\complex^{1/3})$.
      \end{enumerate}
  \end{enumerate}
\end{lem}

  Interpreted in a probabilistic sense, \Cref{lem:SigmaC,lem:SigmaD} immediately
  yield the typical order of a core of a complex graph, as well as the typical deficiency.

\begin{coro}\label{coredeficiency}
    For every function $\epsilon(n_\complex)=o(1)$, if $n_\complex,l\to\infty$ and
    $l\le \epsilon  n_\complex$, then \whp\ $\complex=\complex_g(n_\complex,n_\complex+l)$ has a core with $\sqrt{3n_\complex l}(1+o(1))$ vertices.
    Furthermore, the deficiency of $\complex$ is given by
    \begin{equation*}
      \deficiency(\complex) =
      \begin{cases}
        0 \text{ \whp} & \text{if }l = o(n_\complex^{1/3}),\\
        O_p(1) & \text{if }l = \Theta(n_\complex^{1/3}),\\
        O\br{\sqrt{\frac{l^3}{n_\complex}}} \text{\whp} & \text{if }l = \omega(n_\complex^{1/3}).
      \end{cases}
    \end{equation*}
\end{coro}

Observe that \Cref{coredeficiency} requires $n_\complex$ and $l$ to be growing and $l$ to be of smaller
order than $n_\complex$. We shall later see that this will \whp\ be the case for the \comppart\
of $\general_g(n,m)$.

In addition to \Cref{coredeficiency}, which tells us the deficiency and the order of the core of
$\complex_g(n_\complex,n_\complex+l)$, \Cref{lem:SigmaD} also enables us to express
the number of complex graphs that are embeddable on $\Sg$.

\begin{coro}\label{complex:number}
  For all positive admissible values $n_\complex,l$, we have
  \begin{equation*}
    \abs{\class{\complex}_g(n_\complex,n_\complex+l)}=\frac{n_\complex^{n_\complex+3l/2-1/2}\abs{\class{\kernel}_g(2l,3l)}e^{3l/2}}{(3l)^{(3l-1)/2}(2l)!}\exp\br{f_d}.
  \end{equation*}
\end{coro}

This finalises our analysis of \eqref{eq:complexcore} and \eqref{eq:core}.

\subsection{Complex part and excess}\label{sec:complexexcess}

In this section we derive the main contributions with respect to $n_\complex$ and $l$ to the double 
sum \eqref{eq:general}.
In the previous section, we had to distinguish the cases $n_\complex=0$ and $n_\complex>0$ in order
to determine the number of complex graphs. Similarly, it will turn out that our asymptotic formulas will be quite different
depending on whether the number $m_\outside=m-n_\complex-l$ of edges \emph{outside} the complex part is zero or not.
In order to keep expressions simple, we will deal with the special cases $n_\complex=0$ and $m_\outside=0$
separately.

To this end, define $\class{\general}_g^*(n,m)$ to be the subclass of $\class{\general}_g(n,m)$ consisting
of all graphs for which the complex part is non-empty and the non-complex part has at least one edge.
After bounding $\abs{\class{\general}_g^*(n,m)}$, we shall see that the two special cases $n_\complex=0$ and $m_\outside=0$ are `rare' in the sense that almost all graphs in
$\class{\general}_g(n,m)$ are also in $\class{\general}_g^*(n,m)$.

\begin{prop}\label{lem:exceptional}
  For every $m=m(n)$ as in \Cref{main1}\ref{firstsup}, \Cref{main2}, or \Cref{main3} we have
  \begin{equation*}
    \abs{\class{\general}_g(n,m)\setminus\class{\general}_g^*(n,m)}=o\br{\abs{\class{\general}_g^*(n,m)}}.
  \end{equation*}
\end{prop}

By \Cref{lem:exceptional}, we can determine the main contributions to \eqref{eq:general} by deriving the
main contributions to the corresponding sum for $\abs{\class{\general}_g^*(n,m)}$, namely
\begin{equation}\label{eq:generalstar}
   \abs{\class{\general}_g^*(n,m)} =
  \sum_{n_\complex,l}\binom{n}{n_\complex}\abs{\class{\complex}_g(n_\complex,n_\complex+l)}\cdot\abs{\class{\outside}(n_\outside,m_\outside)},
\end{equation}
where $n_\complex$ and $l$ take all admissible values with $n_\complex>0$ and $m_\outside>0$.

In order to analyse \eqref{eq:generalstar}, we derive an upper bound for the sum over $n_\complex$ and
subsequently also for the sum over $l$. These upper bounds indicate which intervals $I_\complex(n)$ and $I_l(n)$
for $n_\complex$ and $l$, respectively, `should' provide the main contribution to \eqref{eq:generalstar}. For
$n_\complex$ and $l$ from these intervals, we then derive lower bounds and prove that the lower bound for
$n_\complex\in I_\complex(n)$ and $l\in I_l(n)$ is much larger than the tails of the upper bound, thus 
proving that the main contribution to \eqref{eq:generalstar} is indeed provided by $n_\complex\in I_\complex(n)$
and $l\in I_l(n)$.

Applying \eqref{binomial:refined}, \Cref{complex:number}, \Cref{britikov}, and \Cref{cubic:enum}
to \eqref{eq:generalstar}, we have
\begin{equation}\label{second:start}
  \abs{\class{\general}_g^*(n,m)} = \Theta(1)n^{n+\frac12}\br{\frac{e}{2}}^m\sum_l l^{\frac{5g}{2}-3-\frac{3l}{2}}\phi^l\sum_{n_\complex}\rho(n_\outside,m_\outside) \psi(n_\complex,l),
\end{equation}
where $\phi=2\sqrt{e}\gamma_\kernel^23^{-\frac32}$ and
\begin{equation}\label{psi}
  \psi(n_\complex,l) = \left(\frac{2}{e}\right)^{n_\complex}n_\complex^{\frac{3l}{2}-1}n_\outside^{2m_\outside-n_\outside-\frac12}m_\outside^{-m_\outside-\frac12}\exp\br{f_d}.
\end{equation}

Consider the sum
\begin{equation*}
  \Sigma_\complex=\Sigma_\complex(n,m,l):=\sum_{n_\complex}\rho(n_\outside,m_\outside) \psi(n_\complex,l),
\end{equation*}
where we sum over all values of $n_\complex$ that are admissible in $\class{\general}_g^*(n,m)$.
We shall see in \Cref{nQandl} that for fixed $l>0$, the main contribution to $\Sigma_\complex$ is
centred around
\begin{equation*}
  \ol{n}_\complex = 2m-n-2l.
\end{equation*}
The corresponding numbers of vertices and edges in the non-complex components are given by
\begin{equation*}
  \ol{n}_\outside = 2(n-m+l)
  \qquad\text{and}\qquad
  \ol{m}_\outside = n-m+l.
\end{equation*}
The bounds for $\Sigma_\complex$ will depend on whether $l$ is `small' or `large', more precisely, whether
\begin{equation}\label{eq:casesl}
  9\ol{m}_\outside^2\br{\frac{3l}{2}-1} \le \ol{n}_\complex^3
\end{equation}
is satisfied (if so, $l$ is considered small) or not (if so, $l$ is large). 

\begin{lem}\label{lem:SigmaQ}
  Define
  $M_\complex=M_\complex(n,m,l)$ by
  \begin{equation*}
    M_\complex =
    \begin{cases}
      \displaystyle\br{\frac{2}{e}}^{2m-n}\ol{n}_\complex^{\frac{3l}{2}-1}\ol{m}_\outside^{-\ol{m}_\outside-1} & \text{if }\eqref{eq:casesl}\text{ holds},\\
      \displaystyle\br{\frac{2}{e}}^{2m-n}l^{\frac{l}{2}-\frac13}\ol{m}_\outside^{-\ol{m}_\outside+l-\frac53} & \text{otherwise}.
    \end{cases}
  \end{equation*}
  Then
  \begin{equation}\label{eq:SigmaQ}
    \Sigma_\complex \le n^{\frac{3}{2}}\exp\br{O(l)}M_\complex.
  \end{equation}
  Furthermore, for every fixed positive valued function $\epsilon=\epsilon(n)=o(1)$ and every $\delta>0$,
  there exists $N\in\N$ such that for all $n\ge N$
  \begin{equation}\label{eq:SigmaQ:better}
    \Sigma_\complex \le \Theta(1)n^{\frac{3}{2}}\br{\frac{e}{2}}^{2l}(1+\delta)^lM_\complex,
  \end{equation}
  whenever
  \begin{equation}\label{eq:smalll}
    9\ol{m}_\outside^2\br{\frac{3l}{2}-1} \le \epsilon\ol{n}_\complex^3.
  \end{equation}
\end{lem}
For the case that $m$ is larger than $\frac{n}{2}$ by only a small margin, we prove a stronger bound
with the help of \Cref{britikov}\ref{britikov:sup} and a more careful analysis of the sums involved.
\begin{lem}\label{lem:SigmaQ:small}
  Let $m=(1+\edgesfirst n^{-1/3})\frac{n}{2}$ with $\edgesfirst=o(n^{1/12})$ and $\edgesfirst\to\infty$.
  Then we have
  \begin{equation}\label{eq:SigmaQ2}
    \Sigma_\complex \le \edgesfirst n^{\frac{2}{3}}\exp\br{O(l)}M_\complex.
  \end{equation}
\end{lem}

In \Cref{lem:SigmaQ,lem:SigmaQ:small}, the exact bound depends on whether \eqref{eq:casesl} is satisfied or violated. 
Correspondingly, we set
\begin{equation*}
  \Sigma_l := \sum_l l^{\frac{5g}{2}-3-\frac{3l}{2}}\phi^l\Sigma_\complex(n,m,l),
\end{equation*}
where $l$ takes all admissible values for which \eqref{eq:casesl} holds, and
\begin{equation*}
  \tilde\Sigma_l := \sum_l l^{\frac{5g}{2}-3-\frac{3l}{2}}\phi^l\Sigma_\complex(n,m,l),
\end{equation*}
where $l$ takes all admissible values for which \eqref{eq:casesl} is violated. Heuristically, $\Sigma_l$ should be the larger
of the two sums, because $l^{-\frac{3l}{2}}$ should be the dominating term and this term is small when $l$ is large (which is 
the case when \eqref{eq:casesl} is violated). We shall see in \Cref{lem:caseB} that $\tilde\Sigma_l$ is indeed negligible.

Accordingly, we focus on $\Sigma_l$ for the moment.
Applying the bound \eqref{eq:SigmaQ}, we have $\Sigma_l \le \Sigma_l^+$, where
\begin{equation*}
  \Sigma_l^+ = \br{\frac{2}{e}}^{2m-n}\sum_l l^{\frac{5g}{2}-3-\frac{3l}{2}}\phi^l\ol{n}_\complex^{\frac{3l}{2}-1}\ol{m}_\outside^{-\ol{m}_\outside}\exp(O(l)).
\end{equation*}
The main contribution to $\Sigma_l^+$ should be centred around its largest summand. We approximate the largest
summand by ignoring polynomial terms and replacing the term $\exp(O(l))$ by $(e/2)^{2l}$ (which we saw in \Cref{lem:SigmaQ} 
to be a good approximation when \eqref{eq:smalll} holds). The remaining terms attain their largest value at the unique solution $l_0$ of the equation
\begin{equation}\label{eq:l0}
  l_0 = \frac{\phi^{2/3}(2m-n-2l_0)}{e^{1/3}2^{4/3}(n-m+l_0)^{2/3}}\,,\qquad m-n<l_0<m-\frac{n}{2}\,.
\end{equation}
Before we proceed to prove that the main contribution to $\abs{\class{\general}_g^*(n,m)}$ is indeed provided by
$l$ `close to' $l_0$ (and thus the `typical excess' of a graph in $\class{\general}_g^*(n,m)$ is close
to $l_0$), let us take a closer look at the value $l_0$. We introduce
the following notation for the seven different cases of $m(n)$ from our main results.

\begin{enumerate}
\item[\firstsub:]
  $m(n)=(1+\edgesfirst n^{-1/3})\frac{n}{2}$ with $\edgesfirst=\edgesfirst(n)=o(n^{1/3})$ and $\edgesfirst\to-\infty$, the \emph{first subcritical regime};
\item[\firstcrit:]
  $m(n)=(1+\edgesfirst n^{-1/3})\frac{n}{2}$ with $\edgesfirst\to c_{\edgesfirst}\in\R$, the \emph{first critical regime};
\item[\firstsup:]
  $m(n)=(1+\edgesfirst n^{-1/3})\frac{n}{2}$ with $\edgesfirst=o(n^{1/3})$ and $\edgesfirst\to\infty$, the \emph{first supercritical regime};
\item[\intermediate:]
  $m(n)=\intdegree\frac{n}{2}$ with $\intdegree=\intdegree(n)\to c_{\intdegree}\in(1,2)$, the \emph{intermediate regime};
\item[\secondsub:]
  $m(n)=(2+\edgessecond n^{-2/5})\frac{n}{2}$ with $\edgessecond=\edgessecond(n)=o(n^{2/5})$ and $\edgessecond\to-\infty$, the \emph{second subcritical regime};
\item[\secondcrit:]
  $m(n)=(2+\edgessecond n^{-2/5})\frac{n}{2}$ with $\edgessecond\to c_{\edgessecond}\in\R$, the \emph{second critical regime};
\item[\secondsup:]
  $m(n)=(2+\edgessecond n^{-2/5})\frac{n}{2}$ with $\edgessecond=o((\log n)^{-2/3}n^{2/5})$ and $\edgessecond\to\infty$, the \emph{second supercritical regime}.
\end{enumerate}

The union of the first three cases will also be referred to as \emph{the first phase transition}, while the union
of the last three cases is called \emph{the second phase transition}. In \firstsub\ and \firstcrit, our main results
will follow from well-known results. Thus, for the rest of this section, we assume that we are in one of the other
five cases.

The definition of $l_0$ immediately yields its asymptotic order.
\begin{lem}\label{sizesl0}
  The value $l_0$ defined in \eqref{eq:l0} is positive and \whp\ satisfies
  \begin{equation*}
    l_0 =
    \begin{cases}
      \Theta(\edgesfirst) & \text{in \firstsup},\\
      \Theta(n^{1/3}) & \text{in \intermediate},\\
      \Theta\br{\abs{\edgessecond}^{-2/3}n^{3/5}} & \text{in \secondsub},\\
      \Theta(n^{3/5}) & \text{in \secondcrit},\\
      \frac12\edgessecond n^{3/5} + \Theta\br{\edgessecond^{-3/2}n^{3/5}} & \text{in \secondsup}.
    \end{cases}
  \end{equation*}
  Furthermore, in \secondcrit, we have $0<l_0-\frac12 \edgessecond n^{3/5}=\Theta\br{n^{3/5}}$.
\end{lem}

In general, $l_0$ will not be an integer and thus in particular not admissible. Set
\begin{equation*}
  l_1 := \left\lceil l_0 \right\rceil.
\end{equation*}
Now \eqref{eq:l0} and \Cref{sizesl0} yield
\begin{equation}\label{eq:l1}
  l_1 = (1+o(1))\frac{\phi^{2/3}(2m-n-2l_1)}{e^{1/3}2^{4/3}(n-m+l_1)^{2/3}}\,.
\end{equation}
From \Cref{sizesl0} we deduce that all $l$ 'close to' $l_1$ are admissible and use this fact to derive a lower bound on $\abs{\class{\general}_g^*(n,m)}$.
\begin{lem}\label{nQandl:lower}
  Let $c>1$ be given and suppose that $l\in\N$ with $\frac{l_0}{c} \le l \le c\,l_0$ and
  \begin{equation*}
    0 < \ol{m}_\outside =
    \begin{cases}
      \Theta(n^{3/5}) & \text{in \secondcrit},\\
      \Theta\br{\edgessecond^{-3/2}n^{3/5}} & \text{in \secondsup}.
    \end{cases}
  \end{equation*}
  Then $l$ is admissible. Furthermore, there exists 
  \begin{equation*}
    \tilde n_\complex = \ol{n}_\complex + O\br{\ol{m}_\outside^{2/3}}
  \end{equation*}
  such that
  \begin{equation*}
    \Sigma_\complex(n,m,l) \ge \Theta(1)\br{\frac{e}{2}}^{2l}\ol{m}_\outside^{\frac23}\exp\br{f_d(\tilde n_\complex,l)}M_\complex(n,m,l).
  \end{equation*}
  In particular, for every $\delta>0$ and $n$ large enough,
 \begin{equation*}
   \abs{\class{\general}_g^*(n,m)}\geq \Theta(1)n^{n+\frac12}\br{\frac{e}{2}}^{m+2l_1}l_1^{-\frac{3l_1}{2}}\phi^{l_1}(n-m+l_1)^{\frac23}(1-\delta)^{l_1}M_\complex(l_1,n,m). 
 \end{equation*}
\end{lem}
The bound in \Cref{nQandl:lower} enables us to show that $\tilde\Sigma_l$ is negligible.
\begin{lem}\label{lem:caseB}
  For $n\to\infty$, we have
  \begin{equation*}
    n^{n+\frac12}\br{\frac{e}{2}}^m \tilde\Sigma_l = o\br{\abs{\class{\general}_g^*(n,m)}}.
  \end{equation*}
\end{lem}

\Cref{lem:caseB} implies that the main contribution to $\abs{\class{\general}_g^*(n,m)}$ is provided by the same
intervals that provide the main contribution to $\Sigma_l$. After determining lower bounds for the summands
in \eqref{eq:generalstar}, our aim is to determine the `optimal' intervals in view of \Cref{def:sums}. In other
words, we are looking for intervals $I_\complex(n)$ and $I_l(n)$ such that a) the lower bound, summed over
$I_\complex(n)$ and $I_l(n)$, is much larger than the `tail' of the upper bound and b) $I_\complex(n)$ and $I_l(n)$
are as small as possible. To that end, in the second phase transition, we need an auxiliary result that tells us
that $f_d$ (defined in \Cref{lem:SigmaD}) does not change `too much' if we fix $l$ and change $n_\complex$ by a 
small fraction.

\begin{lem}\label{lem:boundfd}
  Suppose that $m(n)$ lies in \secondsub, \secondcrit, or \secondsup. Let positive valued functions $h=h(n)=\omega(1)$ and $\epsilon=\epsilon(n)=o(1)$
  satisfying $h\epsilon=\omega(1)$ be given. Then for all $\delta>0$, there exists 
  $N\in\N$ such that for all $n>N$, $n_\complex = (1+o(1))n$, and $h \le l \le \frac{n_\complex}{h}$, we have
  \begin{equation*}
    \abs{f_d((1-\epsilon)n_\complex,l)-f_d(n_\complex,l)}\leq \delta\epsilon l.
  \end{equation*}
\end{lem}

With this auxiliary result, we can now determine the desired intervals $I_\complex(n)$ and $I_l(n)$ that
provide the main contribution to $\abs{\class{\general}_g^*(n,m)}$.

\begin{lem}\label{nQandl}
  There exist constants $\beta_l^+,\beta_l^-\in\R^+$ and functions 
  $\eta_l^+,\eta_l^-: (1,2)\to\R^+$, and $\vartheta_l^+,\vartheta_l^-\colon \R\to\R^+$ 
   with 
   \begin{align*}
    &\beta_l^+>\beta_l^-,&
    &\eta_l^+(x)>\eta_l^-(x),&
    &\vartheta_l^+(x)>\vartheta_l^-(x)>\frac{x}{2}&
  \end{align*}
  for all $x\in\R$ such that the following holds.
 
  For every fixed function $h=h(n)=\omega(1)$, the main contribution to \eqref{second:start} is provided by
  $l\in I_l(n)$ and $n_\complex\in I_\complex^h(n,l)$, where
  \begin{equation*}
    I_l(n) =
    \begin{cases}
      \left\{ k\in\N \mid \beta_l^-\edgesfirst \le k \le \beta_l^+\edgesfirst \right\} & \text{in \firstsup},\\
      \left\{ k\in\N \mid \eta_l^-(c_{\intdegree})n^{1/3}
        \le k \le \eta_l^+(c_{\intdegree})n^{1/3} \right\} & \text{in \intermediate},\\
      \left\{ k\in\N \mid \beta_l^-\abs{\edgessecond}^{-2/3}n^{3/5} \le k \le \beta_l^+\abs{\edgessecond}^{-2/3}n^{3/5} \right\} & \text{in \secondsub},\\
      \left\{ k\in\N \mid \vartheta_l^-(c_{\edgessecond})n^{3/5} \le k \le \vartheta_l^+(c_{\edgessecond})n^{3/5} \right\} & \text{in \secondcrit},\\
      \left\{ k\in\N \mid \beta_l^-\edgessecond^{-3/2}n^{3/5} \le k - \frac12\edgessecond n^{3/5} \le \beta_l^+\edgessecond^{-3/2}n^{3/5} \right\} & \text{in \secondsup},
    \end{cases}
  \end{equation*}
  and
  \begin{equation*}
    I_\complex^h(n,l) = \left\{ k\in\N \mid \abs{k-\ol{n}_\complex} \le h \ol{m}_\outside^{2/3} \right\}.
  \end{equation*}
\end{lem}

\section{Internal structure}\label{sec:structure}

In the \Cref{sec:kernelcorecomplex}, we have determined the main contributions to $\abs{\class{\general}_g^*(n,m)}$ and
thus, by \Cref{lem:exceptional}, also the main contributions to $\abs{\class{\general}_g(n,m)}$. Interpreting
these results in a probabilistic sense, we deduce the typical orders $n_\complex,n_\core$ of the \comppart\ and the core of
$G=\general_g(n,m)$, respectively, as well as its typical excess $\ex(G)$ and deficiency $\deficiency(G)$. 
All results in this section are proved in \Cref{sec:proofs:main}.

The \comppart, for instance, grows from
order $\edgesfirst n^{2/3}$ in the first supercritical regime to linear order in the intermediate regime. The number $m_\outside$ of edges
\emph{outside} the \comppart\ is about half the number $n_\outside$ of \emph{vertices} outside the \comppart.

\begin{thm}\label{allsizesfirst}
  Let $G=\general_g(n,m)$. Then $n_\complex$, $n_\core$, $\ex(G)$, and $\deficiency(G)$ \whp\ lie in the following ranges.
  \begin{center}
  \begin{tabular}{c|cc}
    & \firstsup & \intermediate \\
    \hline
    $n_\complex$ & $\edgesfirst n^{2/3}+O_p(n^{2/3})$ & $\br{\intdegree-1}n+O_p\br{n^{2/3}}$\\
    $n_\core$ & $\Theta\br{\edgesfirst n^{1/3}}$ & $\Theta\br{n^{2/3}}$\\
    $\ex(G)$ & $\Theta\br{\edgesfirst}$ & $\Theta\br{n^{1/3}}$\\
    $\deficiency(G)$ & $0$ & $O_p(1)$
  \end{tabular}
  \end{center}
  Furthermore, 
  \begin{equation*}
    m_\outside = \frac{n_\outside}{2}+O_p(n_\outside^{2/3}).
  \end{equation*}
\end{thm}

In the second phase transition, the \comppart\ covers almost all vertices and thus, it is more
convenient to consider the number $n_\outside=n-n_\complex$ of vertices \emph{outside} the \comppart.

\begin{thm}\label{allsizessecond}
  Let $G=\general_g(n,m)$. Then $n_\outside$, $n_\core$, $\ex(G)$, and $\deficiency(G)$
   \whp\ lie in the following ranges.
  \begin{center}
    \begin{tabular}{c|ccc}
      & \secondsub & \secondcrit & \secondsup \\
      \hline
      $n_\outside$ & $\abs{\edgessecond}n^{3/5}+O_p\br{\abs{\edgessecond}^{2/3}n^{2/5}}$ & $\Theta\br{n^{3/5}}$ & $\Theta\br{\edgessecond^{-3/2}n^{3/5}}$ \\
      $n_\core$ & $\Theta\br{\abs{\edgessecond}^{-1/3}n^{4/5}}$ & $\Theta\br{n^{4/5}}$ & $\Theta\br{\edgessecond^{1/2}n^{4/5}}$ \\
      $\ex(G)$ & $\Theta\br{\abs{\edgessecond}^{-2/3}n^{3/5}}$ & $\Theta\br{n^{3/5}}$ & $\frac12\edgessecond n^{3/5}+\Theta\br{\edgessecond^{-3/2}n^{3/5}}$ \\
      $\deficiency(G)$ & $O\br{\abs{\edgessecond}^{-1}n^{2/5}}$ & $O\br{n^{2/5}}$ & $O\br{\edgessecond^{3/2}n^{2/5}}$
    \end{tabular}
  \end{center}
  Furthermore, we have
  \begin{equation*}
    n_\outside = 2\ex(G) - \edgessecond n^{3/5} + O_p\br{\br{2\ex(G) - \edgessecond n^{3/5}}^{2/3}}
  \end{equation*}
  and
  \begin{equation*}
    m_\outside = \frac{n_\outside}{2}+O_p(n_\outside^{2/3}).
  \end{equation*}
\end{thm}

As an immediate corollary of \Cref{allsizesfirst,allsizessecond}, we deduce the typical order of the
kernel of $G=\general_g(n,m)$.

\begin{coro}\label{kernelsize}
  The number $n_\kernel$ of vertices and $m_\kernel=\frac{3}{2}n_\kernel+\deficiency(G)$ of edges of
  the kernel of $G=\general_g(n,m)$ lie in the following ranges \whp.
  \begin{center}
  \begin{tabular}{c|ccccc}
    & \firstsup & \intermediate & \secondsub & \secondcrit & \secondsup \\
    \hline
    $n_\kernel$ & $\Theta(\edgesfirst)$ & $\Theta(n^{1/3})$ & $\Theta\br{\abs{\edgessecond}^{-2/3}n^{3/5}}$ & $\Theta\br{n^{3/5}}$ & $\frac12\edgessecond n^{3/5}+\Theta\br{\edgessecond^{-3/2}n^{3/5}}$\\
    $\deficiency(G)$ & $0$ & $O_p(1)$ & $O\br{\abs{\edgessecond}^{-1}n^{2/5}}$ & $O\br{n^{2/5}}$ & $O\br{\edgessecond^{3/2}n^{2/5}}$
  \end{tabular}
  \end{center}
\end{coro}

\Cref{allsizesfirst,allsizessecond} and \Cref{kernelsize} tell us the orders of the \comppart, the core, and the kernel.
What we are ultimately looking for, however,
are orders of components. \Cref{cubic:structure,cubic:lc-size} cover the case of cubic kernels, which
are precisely the kernels of $\general_g(n,m)$ in \firstsup. However, we are not interested in the
properties a kernel has if we pick it uniformly at random from the class of all kernels. We are
rather looking for properties of the \emph{kernel of $\general_g(n,m)$}, where the randomness lies
in $\general_g(n,m)$.
Clearly, we cannot expect the probability distribution on the class of kernels given by this
construction to be uniform.

However, by a double counting argument, we prove that the aforementioned probability distribution does not differ
`too much' from the uniform distribution if we are in \firstsup\ or \intermediate. From this, we use 
Markov's inequality \eqref{eq:Markov} to deduce that in these regimes, the kernel $\kernel_G$,
the core $\core_G$, and the complex part $\complex_G$ of $G=\general_g(n,m)$ have a component of genus $g$ 
that covers almost all vertices, while all other components are planar. 
Recall that $H_i(G')$ denotes the $i$-th largest component of a graph $G'$.
Denote by $R(G')$ the graph $G'\setminus H_1(G')$.

\begin{thm}\label{generalstructure}
  Let $G=\general_g(n,m)$, where $m=m(n)$ lies in \firstsup\ or \intermediate.
  \begin{enumerate}
  \item\label{fewcomponents}
    $\kernel_G$, $\core_G$, and $\complex_G$ have the same number $k=O_p(1)$ of components;
  \item\label{icomponents}
    for every $i\ge2$, the probability that $\kernel_G$, $\core_G$, and $\complex_G$
    have at least $i$ components is bounded away both from $0$ and $1$;
  \item\label{nonplanarcomp}
    \whp\ $H_1(\kernel_G)$, $H_1(\core_G)$, and $H_1(\complex_G)$ have genus $g$;
  \item\label{planarcomp}
    \whp\ $R(\kernel_G)$, $R(\core_G)$, and $R(\complex_G)$ are planar;
  \item\label{componentsizes}
    if $\kernel_G$, $\core_G$, and $\complex_G$ have at least $i\ge 2$ components, then
    \begin{equation*}
      \abs{H_i(\kernel_G)} = \Theta_p(1),
      \quad
      \abs{H_i(\core_G)} = \Theta_p\br{n^{1/3}},
      \quad
      \abs{H_i(\complex_G)} = \Theta_p\br{n^{2/3}};
    \end{equation*}
  \item\label{kernelplanar}
    $\abs{R(\kernel_G)}=O_p(1)$;
  \item\label{coreplanar}
    $\abs{R(\core_G)}=O_p(n^{1/3})$;
  \item\label{complexplanar}
    $\abs{R(\complex_G)}=O_p(n^{2/3})$.
  \end{enumerate}
\end{thm}

For the second phase transition, the proof method of \Cref{generalstructure} fails. For these cases, we prove 
the existence of the giant component by different means in \Cref{sec:proofs:main}.

From \Cref{generalstructure},
we deduce the typical order of the largest components of the \comppart, the core, and the kernel of $\general_g(n,m)$,
respectively.

\begin{prop}\label{compsizesfirst}
  For $G=\general_g(n,m)$, the largest components of the \comppart\ $\complex_G$, the core $\core_G$, 
  and the kernel $\kernel_G$, respectively, have the following order.

  \centering
  \begin{tabular}{c|cc}
    & \firstsup & \intermediate \\
    \hline
    $\abs{H_1(\complex_G)}$ & $\edgesfirst n^{2/3}+O_p(n^{2/3})$ & $\br{\intdegree-1}n+O_p\br{n^{2/3}}$\\
    $\abs{H_1(\core_G)}$ & $\Theta\br{\edgesfirst n^{1/3}}$ & $\Theta\br{n^{2/3}}$\\
    $\abs{H_1(\kernel_G)}$ & $\Theta\br{\edgesfirst}$ & $\Theta\br{n^{1/3}}$
  \end{tabular}
\end{prop}

\section{Proofs of main results}\label{sec:proofs:main}

In this section, we prove the main results (\Cref{main1,main2,main3}) of this paper,
as well as the structural results from \Cref{sec:structure}.

\subsection{Proof of \Cref{main1}}

In \firstsub, i.e.\ $m=(1+\edgesfirst n^{-1/3})\frac{n}2$ with $\edgesfirst=o(n^{1/3})$
and $\edgesfirst\to-\infty$, the \ER\ random graph $G(n,m)$ \whp\ is embeddable on $\Sg$ by \Cref{britikov}.
Thus, \Cref{main1}\ref{firstsub} follows immediately from \Cref{thm:ER}\ref{Gnm:sub}.

In \firstcrit, i.e.\ $\edgesfirst\to c_{\edgesfirst}\in\R$, \Cref{britikov}\ref{britikov:crit} implies that
$G(n,m)$ has no complex components with positive probability. Thus, \Cref{thm:ER}\ref{Gnm:crit} 
yields the second statement of \Cref{main1}\ref{firstcrit}. By \cite[Theorem 5]{luczakpittel}, the probability that $G(n,m)$ is 
planar, and thus in particular embeddable on $\Sg$, is larger than the probability that $G(n,m)$ has no complex components. 
Hence the first statement of \Cref{main1}\ref{firstcrit} follows as well.

In \firstsup, i.e.\ $\edgesfirst=o(n^{1/3})$ and $\edgesfirst\to\infty$, by 
\Cref{generalstructure}\ref{nonplanarcomp}--\ref{componentsizes} and \Cref{compsizesfirst}, the \comppart\ of $G=\general_g(n,m)$ \whp\ has one component that 
has genus $g$ and order $\edgesfirst n^{2/3}+O_p(n^{2/3})$, while all other components are planar and 
have order $\Theta_p(n^{2/3})$. By \Cref{generalstructure}\ref{fewcomponents} and \ref{icomponents},
it remains to show that for each $i\ge 1$, the \mbox{$i$-th} largest non-complex component has order $\Theta_p\br{n^{2/3}}$. By \Cref{britikov}\ref{britikov:crit}
and the fact that $m_\outside = \frac{n_\outside}{2}+O_p\br{n_\outside^{2/3}}$, there is a positive probability that
$G(n_\outside,m_\outside)$ has no complex component and therefore the claim follows from
\Cref{thm:ER}\ref{Gnm:crit}.\qed

\subsection{Proof of \Cref{main2}}

Let $m(n)$ be a function from the second phase transition, that is, $m(n)=(2+\edgessecond n^{-2/5})\frac{n}2$ 
with $\edgessecond = \edgessecond(n) = o(n^{2/5})$. Again, we denote the number $n-n_\complex$ of
vertices outside the \comppart\ of a given graph $G\in\class{\general}_g(n,m)$ by $n_\outside$ and the number 
of edges outside the \comppart\ by $m_\outside$. We claim that $n_\complex-\abs{H_1(\complex_G)} = O_p(n_\outside)$. 
In other words, for every $\delta>0$ we need to find a constant $c_\delta$ so 
that $n_\complex-\abs{H_1(\complex_G)} \le c_\delta n_\outside$ with probability greater than $1-\delta$ 
for sufficiently large $n$. Fix $\delta >0$ and denote by $\class{E}_g(n,m)$ the subclass of 
$\class{\general}_g(n,m)$ of those graphs $G$ for which $n_\complex-\abs{H_1(\complex_G)} > c_\delta n_\outside$ 
with $c_\delta := \frac5\delta$. We have to prove that $\abs{\class{E}_g(n,m)} < \delta \abs{\class{\general}_g(n,m)}$ 
for sufficiently large $n$.
 
Suppose that there exists an infinite set $I\subset\N$ such that 
$\abs{\class{E}_g(n,m)} \ge \delta\abs{\class{\general}_g(n,m)}$ for all $n\in I$. We use double counting 
in order to derive a contradiction from this assumption. Let $n\in I$ be fixed
and pick a graph $G\in\class{E}_g(n,m)$. \Cref{allsizessecond} together with the assumption
$\abs{\class{E}_g(n,m)} \ge \delta\abs{\class{\general}_g(n,m)}$ yields that
\begin{equation*}
  m_\outside = m-m_\complex  = \frac{n_\outside}{2}+O_p\br{n_\outside^{2/3}}.
\end{equation*}
By definition, $\abs{H_1(\complex_G)} < n_\complex-\frac5\delta n_\outside$. Thus, there is a partition 
$(A,B)$ of the vertices in $\complex_G$ such that each component is contained either in $A$ or in 
$B$ and that $\abs{A} \ge \abs{B} \ge \frac5\delta n_\outside$. Now we perform the following operation. 
We delete one edge from the non-complex components and instead add an edge between $A$ 
and $B$. The resulting graph is still embeddable on $\Sg$ and thus lies in $\class{\general}_g(n,m)$.
The number of choices for this operation is therefore
\begin{equation*}
  m_\outside\abs{A}\cdot\abs{B} %\ge (1+o(1))\frac{n_\outside}{2}\cdot\frac{n_\complex}{2}\cdot \frac5\delta n_\outside
  \ge (1+o(1))\frac{5}{4\delta}n_\outside^2n_\complex.
\end{equation*}
The reverse operation is to delete an edge $uv$ from the \comppart\ that separates $u$ and $v$ and add an edge outside the 
\comppart\ (not creating any new complex components). There are less than $n_\complex$ choices for $uv$, because any spanning 
tree of a component has to contain all edges of that component that are feasible for $uv$. Thus, there are less than 
$n_\outside^2n_\complex$ possibilities for the reverse operation, yielding
\begin{equation*}
  (1+o(1))\frac{5}{4\delta}n_\outside^2n_\complex  \abs{\class{E}_g(n,m)} < n_\outside^2n_\complex  \abs{\class{\general}_g(n,m)}
\end{equation*}
and thus
\begin{equation*}
  \abs{\class{E}_g(n,m)} < (1+o(1))\frac{4\delta}{5}  \abs{\class{\general}_g(n,m)}
  < \delta\abs{\class{\general}_g(n,m)}
\end{equation*}
for sufficiently large $n\in I$, a contradiction.

We have thus proved that 
$n_\complex-\abs{H_1(\complex_G)} = O_p(n_\outside)$, which in turn implies that 
$n-\abs{H_1(\complex_G)} = O_p(n_\outside)$. \Cref{main2} now follows from \Cref{allsizessecond} 
and the trivial fact that $n-\abs{H_1(\complex_G)} \ge n_\outside$.\qed

\subsection{Proof of \Cref{main3}}

Analogously to the proof of the case \firstsup\ of \Cref{main1}, \Cref{main3} follows from \Cref{thm:ER}\ref{Gnm:crit},
\Cref{britikov}\ref{britikov:crit}, \Cref{generalstructure}, and \Cref{compsizesfirst}.\qed

\proofof{main4}

In \ref{enum:first}, the cases \firstsub\ and \firstcrit\ follow directly from \Cref{britikov}\ref{britikov:sub}
and \ref{britikov:crit}, respectively, and Stirling's formula, applied to the binomial coefficient
$\binom{\binom{n}{2}}{m}$.

For all other regimes, \Cref{nQandl:lower} provides a lower bound. On the other hand, \Cref{lem:exceptional} and
\Cref{nQandl} tell us that the main contribution to $\abs{\class{\general}_g(n,m)}$ is provided by an interval $I_l(n)$ that is
centred around $l_1=\lceil l_0\rceil$. Moreover, by~\eqref{eq:l0}, there exist constants $0<c^-<c^+$ such that
\begin{equation}\label{eq:l}
  c^-\frac{2m-n-2l}{(n-m+l)^{2/3}} \le l \le c^+\frac{2m-n-2l}{(n-m+l)^{2/3}}
\end{equation}
for all $l\in I_l(n)$. The length of $I_l(n)$ is $O(l_0)$. Let $l_2\in I_l(n)$ be the index that maximises the 
summand. Applying \Cref{lem:SigmaQ,lem:SigmaQ:small} and \eqref{eq:l}, we observe that the resulting upper bound differs from the lower bound from 
\Cref{nQandl:lower} by a factor of the type $\exp(O(l_0))$. Now \Cref{main4} follows by inserting the values for $l_0$ from
\Cref{sizesl0} into the lower bound from \Cref{nQandl:lower}.\qed

\subsection{Proof of \Cref{allsizesfirst,allsizessecond}}

The results on the excess and the order of the complex part follow from \Cref{nQandl}. Observe that
$\ex(G) = o(n_\complex)$ for all regimes and thus \Cref{coredeficiency} is applicable, yielding
the order $n_\core$ of the core and the deficiency $\deficiency(G)$. Finally, by \Cref{nQandl} we
know that
\begin{equation*}
  n_\outside = 2(n-m+\ex(G))+ O_p\br{(n-m+\ex(G))^{2/3}}
\end{equation*}
and
\begin{equation*}
  m_\outside = n-m+\ex(G)+ O_p\br{(n-m+\ex(G))^{2/3}},
\end{equation*}
which yields the last statements of \Cref{allsizesfirst,allsizessecond}.\qed

%-------------------------------------------------------------------------------------------------------------------

\proofof{kernelsize}

\Cref{kernelsize} follows directly from \Cref{kernel-ld} and the values of $\ex(G)$ and $\deficiency(G)$
stated in \Cref{allsizesfirst,allsizessecond}.\qed

\proofof{generalstructure}

  Given a fixed kernel $\kernel$, call a subdivision of $\kernel$ \emph{good} if it is a simple
  graph (and thus a valid core). We first prove that the fraction of good subdivisions
  among all subdivisions of $\kernel$ is bounded away from zero.

  To this end, suppose that $\kernel$ is a kernel with $2l-d$ vertices and $3l-d$
  edges and that we want to subdivide its edges $k$ times (with $k\ge 6l-2d$) in order to construct a core $\core$ with
  $k+2l-d$ vertices. We subdivide $\kernel$ in the following way. First, decide which labels the
  vertices of $\kernel$ should have in $\core$; there are $\binom{k+2l-d}{2l-d}$ choices for this.
  Let $I$ be the set of the remaining $k$ labels. We recursively subdivide edges of $\kernel$
  and assign the smallest remaining label in $I$ to the new vertex. The number of choices increases by one in each recursion step and
  thus we have $(k+3l-d-1)_{k}$ choices in total. This way, we construct each subdivision precisely
  once. Hence the total number of subdivisions of $\kernel$ is
  \begin{equation*}
    \binom{2l-d+k}{2l-d}(k+3l-d-1)_{k}\,.
  \end{equation*}
  In order to give a lower bound on the number of \emph{good} subdivisions, we change our construction
  slightly by introducing a preliminary step. After choosing the labels for the vertices in $\kernel$,
  we subdivide each edge of $\kernel$ twice and then choose labels from $I$ for the new vertices;
  there are $(k)_{6l-2d}$ choices for this. After this step, we proceed as before, with the
  additional rule that an edge may only be subdivided if none of its end vertices is a vertex of
  $\kernel$. Similar to our first construction, there are $(k-3l+d-1)_{k-6l+2d}$ choices for this
  part of the construction. Every graph obtained by this type of subdivision is simple and no graph
  is constructed more than once. Thus, the total number of good subdivisions is at least
  \begin{equation*}
    \binom{2l-d+k}{2l-d}(k)_{6l-2d}(k-3l+d-1)_{k-6l+2d}\,.
  \end{equation*}
  The fraction of good subdivisions among all subdivisions of $\kernel$ is thus at least
  \begin{align*}
    \frac{(k)_{6l-2d}(k-3l+d-1)_{k-6l+2d}}{(k+3l-d-1)_{k}} 
    \ge \left(\frac{k-3l+d}{k-6l+2d}\right)^{-6l+2d}
    \stackrel{\eqref{explogu}}{\ge} \exp\br{-\frac{2(3l-d)^2}{k-6l+2d}}.
  \end{align*}
  Substituting $l=\ex(G)$, $d=\deficiency(G)$, and $k=n_\core-3l+d$ from \Cref{allsizesfirst}
  (and observing that these values satisfy $k\ge 6l-2d$ \whp) yields that the fraction of good subdivisions is bounded away from zero.

  To make this more precise, denote by
  $\Sub(\kernel_G)$ the proportion of subdivisions of $\kernel_G$ that lie in 
  $\class{\core}_g(n_\core,n_\core+l)$. We have shown that for every $\delta>0$ there exists an $\varepsilon>0$
  such that
  \begin{equation}\label{goodsubdivisions}
    \begin{split}
    1-\delta \le \Sub(\kernel_G) &\le 1 \text{ \whp\ in \firstsup,}\\
    \varepsilon \le \Sub(\kernel_G) &\le 1 \text{ with probability at least }1-\delta\text{ in \intermediate.}
    \end{split}
  \end{equation}
  Recall the construction steps \ref{construction:kerneltocore}--\ref{construction:complextogeneral}: the core $\core_G$ is 
  constructed from $\kernel_G$ by subdividing edges; the \comppart\ $\complex_G$ is obtained from $\core_G$ by attaching rooted trees to all vertices; 
  adding trees and unicyclic components to $\complex_G$ yields $G$.
  Let $X$ be an event that depends on the choice of $\kernel\in\class{\kernel}_g$.
  From the above construction, \eqref{goodsubdivisions}, and the fact that the kernel of $G=\general_g(n,m)$ has
  a growing number of vertices by \Cref{allsizesfirst}, we deduce that
  \begin{equation}\label{kernelprob}
    \varepsilon \le \frac{\prob{X \text{ holds for }\kernel=\kernel_G}}{\prob{X \text{ holds for }\kernel=\kernel_g(2l-d,3l-d)}} \le \frac{1}{\varepsilon},
  \end{equation}
  provided that the denominator is non-zero.

  To prove \ref{fewcomponents}, observe that the kernel, the core, and the \comppart\ of a graph have the 
  same number $k$ of components by construction. \Cref{cubic:structure,cubic:lc-size} (for $g\ge 1$) and
  \cite[Lemma 2]{kang2012} (for $g=0$) tell us that the cubic kernel $\kernel_g(2l,3l)$ has $O_p(1)$ components.
  Thus by \eqref{kernelprob}, we have $k=O_p(1)$ if
  the kernel is cubic, which is the case \whp\ in \firstsup\ by \Cref{allsizesfirst}. In \intermediate, we
  have $\deficiency(G)=O_p(1)$. Thus, we apply \Cref{kernelpumping} and deduce that $k=O_p(1)$. By analogous
  arguments, we deduce \ref{icomponents}, \ref{kernelplanar}, and the statements about $\kernel_G$ from
  \ref{nonplanarcomp}, \ref{planarcomp}, and \ref{componentsizes}.

  The observation that subdividing edges (when constructing $\core_G$) and attaching trees (constructing $\complex_G$) 
  does not change the genus of any component proves the remaining statements of \ref{nonplanarcomp} and \ref{planarcomp}.

  In order to prove \ref{componentsizes}, \ref{coreplanar} and \ref{complexplanar}, let $A_\kernel$ be any fixed component 
  of $\kernel_G$. Denote by $A_\core$ and $A_\complex$ the corresponding components of $\core_G$ and $\complex_G$, respectively.
  Observe that
  \begin{itemize}
  \item in a random (not necessarily good) subdivision of the kernel,
    the expected number of subdivisions of any given edge $e$ is $\frac{n_\core}{n_\kernel}-1$;
  \item if we attach a rooted forest to the core in order to construct the \comppart, the
    expected order of the tree attached to any given vertex $v$ is $\frac{n_\complex}{n_\core}$.
  \end{itemize}
  By \Cref{allsizesfirst}, we have $\frac{n_\core}{n_\kernel} = \Theta(n^{1/3})$ and $\frac{n_\complex}{n_\core} =
  \Theta(n^{1/3})$ \whp. Therefore, \eqref{goodsubdivisions} and Markov's inequality \eqref{eq:Markov}, applied to the 
  random variables $\abs{A_\core}$ and $\abs{A_\complex}$, imply that 
  \begin{equation*}
    \abs{A_\core}=O_p(n^{1/3})\abs{A_\kernel}
    \qquad\text{and}\qquad
    \abs{A_\complex}=O_p(n^{1/3})\abs{A_\core}
  \end{equation*}
  for every \emph{fixed} component $A_\kernel$. On the other hand, there are $O_p(1)$ components, which proves \ref{coreplanar} and \ref{complexplanar}.

  It remains to prove the lower bound for $\abs{\component_i(\core_G)}$ and $\abs{\component_i(\complex_G)}$ in \ref{componentsizes}.
  For an edge $e$ of $\kernel_G$, denote by $X_e$ the random variable of subdivisions of $e$. Both the expectation $\expec{X_e}$ and 
  the variance $\sigma^2$ have order $\Theta(n^{1/3})$. Therefore, Chebyshev's inequality~\eqref{eq:Chebyshev} implies that 
  \begin{equation*}
    \prob{X_e \le \frac12\expec{X_e}} = O(n^{-1/3}).
  \end{equation*}
  Thus, for a fixed component $A_\kernel\not=\component_1(\kernel_G)$, a union bound over all $O_p(1)$ edges in $A_\kernel$ 
  proves that $\abs{A_\core}=\Theta_p(n^{1/3})\abs{A_\kernel}$. By another union bound, this is true for all $O_p(1)$ components
  (apart from $\component_1(\kernel_G)$), proving $\abs{\component_i(\core_G)}=\Theta_p(n^{1/3})$ for all $i\ge 2$.

  Similarly, for a vertex $v$ of $\core_G$, denote by $Y_v$ the number of vertices in the tree attached to $v$ when we construct
  $\complex_G$. Again, both the expectation and the variance have order $\Theta(n^{1/3})$ and we deduce
  \begin{equation*}
    \prob{Y_v \le \frac12\expec{Y_v}} = O(n^{-1/3})
  \end{equation*}
  from Chebyshev's inequality~\eqref{eq:Chebyshev}. This implies that, for any given $\delta>0$, there exists an $\varepsilon>0$ such that with probability at
  least $1-\delta$, every component $A_\core$ contains at least $\varepsilon n^{1/3}$ vertices $v$ with $Y_v > \varepsilon n^{1/3}$, which yields
  $\abs{A_\complex}\ge \varepsilon^2n^{2/3}$. This proves \ref{componentsizes} and thus finishes the proof of \Cref{generalstructure}.\qed
%-------------------------------------------------------------------------------------------------------------------

\proofof{compsizesfirst}

\Cref{compsizesfirst} is an immediate consequence of \Cref{allsizesfirst,generalstructure}.\qed

\section{Proofs of auxiliary results}\label{sec:proofs:lemmas}

In this section we prove all results from \Cref{sec:strategy,sec:kernelcorecomplex}.

\subsection{Proof of \Cref{britikov}}

It remains to prove \ref{britikov:sup2}. From Lemma 3, (10.11), and (10.12) in \cite{jansonea},
we deduce that
\begin{equation*}
 \rho(n,m)= \frac{2^{2m-n}e^nm!n!}{(n-m)!n^{2m}2\pi i}\oint \sqrt{1-z}\exp\br{n k(z)}\frac{\dd z}{z},
\end{equation*}
where the contour of the integral is a closed curve around the origin with $\abs{z}\leq1$ and
\begin{equation*}
 k(z)=z-1-\frac{m}{n}\log(z)+\br{1-\frac{m}{n}}\log(2-z).
\end{equation*}
We use the contour consisting of a) the line segment from $1$ to $i$, b) the semicircle of radius one 
with negative real value, and c) the line segment from $-i$ to $1$.
Along this contour we have $\abs{\exp(k(z))}\leq 1$ and thus 
\begin{align*}
 \rho(n,m)&\leq \frac{2^{2m-n}e^nm!n!}{(n-m)!n^{2m}2\pi}\abs{\oint \frac{\sqrt{1-z}}{z}\dd z}\\
 &\stackrel{\eqref{binomial:refined}}{\leq}\frac{e^2(\pi+2\sqrt2)}{\sqrt2\pi^{\frac32}}\br{\frac{2}{e}}^{2m-n}\frac{m^{m+1/2}n^{n-2m+1/2}}{(n-m)^{n-m+1/2}}\,,
\end{align*}
proving the lemma.\qed

\proofof{cubic:lc-size}

We abbreviate the class of cubic kernels embeddable on $\Sg$ by $\class{A}_g$ and the subclass of 
$\class{A}_g$ of \emph{connected} cubic kernels by $\class{B}_g$. Clearly, every graph in $\class{A}_g$ has an even number of vertices.
We first prove \eqref{cubic:planar:eq}.

By \Cref{cubic:enum,cubic:connected} there exist positive constants $a_g^-,a_g^+,b_g^-,b_g^+$ such that for all $l$ 
\begin{equation*}
  a_g^-\leq \frac{\abs{\class{A}_g(2l)}}{(2l)^{5g/2-7/2}\gamma_{\kernel}^{2l}(2l)!}\leq a_g^+
  \qquad\text{and}\qquad
  b_g^-\leq \frac{\abs{\class{B}_g(2l)}}{(2l)^{5g/2-7/2}\gamma_{\kernel}^{2l}(2l)!}\leq b_g^+.
\end{equation*}

By \Cref{cubic:structure}, the elements of $\class{A}_g(2l)$ \whp\ have a unique non-planar component. 
Therefore the probability that $\pl(G)$ has exactly $2i$ vertices is given by 
\begin{equation*}
 \prob{\abs{\pl(G)}= 2i}=(1+o(1))\binom{2l}{2i}\frac{\abs{\class{B}_g\br{2l-2i}}\cdot \abs{\class{A}_0\br{2i}}}{\abs{\class{A}_g\br{2l}}}
\end{equation*}
and we can therefore conclude that \eqref{cubic:planar:eq} holds.

It remains to show that for every $\delta>0$ there exists a constant $c_\delta$ such that 
$\prob{\abs{\pl(G)}>2c_\delta}<\delta$ for sufficiently large $l$.
By \Cref{cubic:structure}, \eqref{cubic:planar:eq}, and the fact that $g\geq1$, we have for any $c_\delta\in\N_{>0}$
\begin{equation*}
 \prob{\abs{\pl(G)}>2c_\delta}\le (1+o(1))\sum_{i=c_\delta+1}^{l-3} c_g^+ i^{-7/2}\br{1-\frac{i}{l}}^{-1}.
\end{equation*}
The summand (as a function in $i$) has a unique minimum at $i=\frac{7l}{9}$. Therefore,
\begin{align*}
 \prob{\abs{\pl(G)}\ge 2c_\delta} &\le (1+o(1))c_g^+\int_{c_\delta}^{l-2}x^{-7/2}\br{1-\frac{x}{l}}^{-1}\dd x\\
  &= \br{\frac{2}{5}+o(1)}c_g^+c_\delta^{-5/2}(1+O(l^{-1/2})) < \delta
\end{align*}
for $c_\delta$ and $l$ large enough, as desired.\qed

%--------------------------------------------------------------------------------------------------------------

\proofof{kernelpumping}

For $\kernel\in\class{\kernel}_g\br{2l,3l}$ and $\ol{\kernel}\in\class{\kernel}_g\br{2l-d,3l-d}$, we say that 
$\kernel$ \emph{contracts} to $\ol{\kernel}$ if for each vertex in $\kernel$ with label 
$i\in\{2l-d+1,\ldots,2l\}$ we can choose an edge $e_i=\{i,v_i\}$ so that contracting these edges results in 
$\ol{\kernel}$ (the contracted vertices obtain the smaller of the two labels). We say that 
$e_{2l-d+1},\dotsc,e_{2l}$ are the \emph{contracted edges}. Denote by $\class{\kernel}_g^{\Delta=4}(2l-d,3l-d)$ the subclass
of $\class{\kernel}_g(2l-d,3l-d)$ consisting of multigraphs with maximum degree four. 
We say that a contraction of $\kernel$ to $\ol{\kernel}$ has \emph{degree four} if $\ol{\kernel}\in\class{\kernel}_g^{\Delta=4}(2l-d,3l-d)$.

If $\kernel$ contracts to $\ol{\kernel}$, then
the compensation factor defined in \eqref{eq:compensationfactor} satisfies
\begin{equation}\label{eq:changedweight}
  w(\ol{\kernel}) \le w(\kernel) \le 6^d w(\ol{\kernel}).
\end{equation}

Each $\kernel\in\class{\kernel}_g(2l,3l)$ contracts in at most $3^d$ ways, because $\kernel$ is 
cubic and hence there are at most $3^d$ choices for the edges $e_{2l-d+1},\dotsc,e_{2l}$.
Vice versa, we claim that every fixed $\ol{\kernel}\in\class{\kernel}_g\br{2l-d,3l-d}$
is obtained by at least $d!2^{-d}$ different contractions from graphs in $\class{\kernel}_g\br{2l,3l}$. 
By recursively splitting vertices of $\ol{\kernel}$ of degree at least four into two new adjacent vertices of 
degree at least three each, not increasing the genus throughout the process, we obtain a weighted multigraph 
$\kernel\in\class{\kernel}_g\br{2l,3l}$ that contracts to $\ol{\kernel}$. The new vertices can be labelled in $d!$ ways,
of which at least $2^{-d}d!$ result in distinct multigraphs in $\class{\kernel}_g\br{2l,3l}$.
Together with \eqref{eq:changedweight}, this proves the upper bound 
\begin{equation*}
  \frac{\abs{\class{\kernel}_g(2l-d,3l-d)}}{\abs{\class{\kernel}_g(2l,3l)}} \le \frac{6^d}{d!}.
\end{equation*}
The corresponding bound for $\abs{\class{\kernel}_g(2l-d,3l-d;\class{P}_i)}$ follows analogously observing that
the two constructions above do neither change the number of components nor increase the genus of any component.

For the lower bound, we claim that the elements of $\class{\kernel}_g(2l,3l)$ have at least $6^{-d}$ contractions of degree four \emph{on average}.
Indeed, first observe that $\kernel\in\class{\kernel}_g\br{2l,3l}$ contracts to $\ol{\kernel}\in\class{\kernel}_g^{\Delta=4}(2l-d,3l-d)$ 
if and only if the contracted edges form a matching in $\kernel$. By choosing the edges of the matching
recursively, we see that $\kernel$ contains at least $\br{2^dd!}^{-1}\prod_{j=0}^{d-1}(2l-6j)$ matchings of size $d$.

Denote by $\class{A}(\kernel)$ the class of all weighted multigraphs that are isomorphic to $\kernel$. 
If we choose $A\in\class{A}(\kernel)$ and a matching $M$ of size $d$ in $A$ uniformly at random, then the 
probability that every edge in $M$ has precisely one end vertex with label in $\{2l-d+1,\dotsc,2l\}$ is 
$\frac{2^d}{\binom{2l}{d}}$. Therefore, the average number of contractions of degree four of graphs in $\class{A}(\kernel)$ is at least
\begin{equation*}
  \frac{\prod_{j=0}^{d-1}(2l-6j)}{2^dd!}\cdot\frac{2^d}{\binom{2l}{d}}
  \ge \br{\frac{2l-6d}{2l-d}}^d \ge 6^{-d},
\end{equation*}
where the last inequality uses the fact that $d\le\frac{2l-d}{6}$. The fact that the classes $\class{A}(\kernel)$ partition
$\class{\kernel}_g(2l,3l)$ proves that $\kernel\in\class{\kernel}_g(2l,3l)$ has at least $6^{-d}$ contractions of degree four
on average.

Vice versa, let $\ol{\kernel}\in\class{\kernel}_g^{\Delta=4}\br{2l-d,3l-d}$. By recursively splitting the $d$ vertices of 
degree four in $\ol{\kernel}$, we see that $\ol{\kernel}$ can be obtained by at most $\binom{4}{2}^dd!=6^dd!$ contractions 
of degree four. Together with \eqref{eq:changedweight}, we deduce that
\begin{equation*}
  \frac{\abs{\class{\kernel}_g(2l-d,3l-d)}}{\abs{\class{\kernel}_g(2l,3l)}} \ge \frac{1}{216^dd!}.
\end{equation*}
The corresponding bound for $\abs{\class{\kernel}_g(2l-d,3l-d;\class{P}_i)}$ follows analogously.\qed

%-------------------------------------------------------------------------------------------------------------------

\subsection{Remark}\label{sec:remark:pumping}

Observe that the proof of \Cref{kernelpumping} applies to any class $\class{F}$ of (multi-) graphs that is a) closed under taking minors
and b) \emph{weakly addable}, that is, if $G$ is obtained by adding an edge between two distinct components of $F\in\class{F}$,
then also $G\in\class{F}$. For more details, see \Cref{sec:discussion}.

\proofof{binsandballs}

  Let $\kernel\in\class{\kernel}_g(2l-d,3l-d)$.
  We subdivide the edges of $\kernel$ by inserting $n_\core-2l+d$ vertices and then assign labels to these new vertices
  in one of $(n_\core-2l+d)!$ possible ways so as to obtain a core with $n_\core$ vertices.

  Call a distribution of $n_\core-2l+d$ new vertices to the edges of $\kernel$ \emph{feasible} if
  the resulting graph has no loops or multiple edges. The number
  $\binom{n_\core+l-1}{3l-d-1}$ of \emph{all} distributions is clearly an upper
  bound for the number of feasible distributions. On the other 
  hand, a distribution is feasible if and only if each loop is subdivided at least twice and for every multiple
  edge, at most one of its edges is not subdivided. Denote by $s_\kernel$ the minimal number of times that
  we need to subdivide the edges of $\kernel$ in order to obtain a simple graph. Then
  $\binom{n_\core+l-s_\kernel-1}{3l-d-1}$ is a lower bound on the number of feasible
  distributions.
  
  By construction, $s_\kernel\le 2(3l-d) \le 6l$ and we thus deduce that
  \begin{equation*}
    \min_{-5\le\nu\le1}\br{(n_\core-2l+d)!\binom{n_\core+\nu l -1}{3l-d-1}} \le
    \varphi_{n_\core,l,d} \le
    (n_\core-2l+d)!\binom{n_\core+l -1}{3l-d-1}.
  \end{equation*}
  Now the lemma follows from the intermediate value theorem and the fact that the function
  $\binom{x}{k}$ for fixed $k\in\N$ is continuous for $x\in\R$.\qed

%-------------------------------------------------------------------------------------------------------------------

\proofof{complex:bounds}

\Cref{complex:bounds} follows directly from \eqref{eq:complex}, \Cref{kernelpumping,binsandballs},
the intermediate value theorem, and the fact that $x^d$ is continuous.\qed

%-------------------------------------------------------------------------------------------------------------------

\proofof{lem:SigmaC}

We first derive an upper bound for $\Sigma_\core$, as well as the main contribution to this upper bound.
We substitute $n_\core=\ol{n}_\core+r$ (recall that $\ol{n}_\core=\sqrt{n_\complex(3l-d)}$).
Applying \eqref{falling:bound} to \eqref{sum:core}, and then using \eqref{explogu} and \eqref{bound:frac}
we deduce that
\begin{equation}\label{eq:core:2}
 \Sigma_\core \le \Sigma_\core^+ := \sum_{r}\exp\left(-\frac{r^2}{2n_\complex}+rA_1+A_2\right),
\end{equation}
where
\begin{align*}
 A_1=&\frac{1-2\ol{n}_\core}{2n_\complex}+\frac{3l-d-1}{\ol{n}_\core}+\frac{(3l-d-1)(3l-d-2)}{2(\ol{n}_\core+ l-1)^2}\,,\\
 A_2=&(3l-d)\br{\log\br{\ol{n}_\core}-\frac{1}{2}}+\sqrt{\frac{3l-d}{4n_\complex}}+(3l-d-1)\br{\frac{l-1}{\ol{n}_\core}-\frac{3l-d-2}{2(\ol{n}_\core+l-1)}}.
\end{align*}

Evaluating the `Gaussian' sum in \eqref{eq:core:2} we obtain
\begin{equation*}
\Sigma_\core^+ \leq \sqrt{2\pi n_\complex}\exp\br{A_2+\frac{n_\complex A_1^2}{2}}.
\end{equation*}
The existence of the constants $a_\core^+,b_\core^+$ from \ref{ncore:upper} now follows from
\begin{equation*}
  \exp(A_2) \leq \br{\frac{n_\complex(3l-d)}{e}}^{(3l-d)/2}\exp\br{O\br{\sqrt{\frac{l^3}{n_\complex}}}}
\end{equation*}
and the observation that $n_\complex A_1^2 = O\br{\frac{l^2}{n_\complex}}$, which is $O\br{\sqrt{\frac{l^3}{n_\complex}}}$, because $l=O(n_\complex)$.

In order to prove \ref{ncore:lower}, suppose that $\frac{7}{2}d\le l\le \epsilon n_\complex = o(n_\complex)$; then also $l=o(\ol{n}_\core)$. In \eqref{sum:core},
we set $n_\core=\ol{n}_\core-\nu l+1+s$. If we let the parameter $s=r+\nu l-1$ take only values for which $n_\core\in I_\core^\delta(n_\complex,l,d)$
with fixed $0<\delta<\frac12$, then
\begin{equation*}
  \Sigma_\core \ge \sum_s \frac{(n_\complex)_{n_\core}}{n_\complex^{n_\core}}n_\core\br{\ol{n}_\core+s}_{3l-d-1}.
\end{equation*}
The interval $I_\core^\delta(n_\complex,l,d)$ has length $2\delta\ol{n}_\core > 2\delta\sqrt{n_\complex}$ and hence
we can choose for $s$ an interval $I_s$ of length $\delta\sqrt{n_\complex}$ in which $\abs{s}<\delta\ol{n}_\core$ holds. 

We use \eqref{falling:bound} for both falling factorials 
and obtain
\begin{equation}\label{SigmaC:lower:start}
  \Sigma_\core \ge \ol{n}_\core^{3l-d}\sum_s \br{1+\frac{s}{\ol{n}_\core}}^{3l-d-1}\br{1+\frac{s+1+\nu l}{\ol{n}_\core}}\exp\br{B_1},
\end{equation}
where
\begin{equation*}
  B_1 = -\frac{\br{\ol{n}_\core-\nu l+1+s}^2}{2(n_\complex-\ol{n}_\core+\nu l-1-s)}-\frac{(3l-d-1)^2}{2(\ol{n}_\core-3l+d+1+s)}.
\end{equation*}
Observe that $1+\frac{s}{\ol{n}_\core}=\Theta(1)$ and $1+\frac{s+1+\nu l}{\ol{n}_\core}=\Theta(1)$. Using \eqref{explogl}, we deduce that
\begin{equation}\label{SigmaC:lower:2}
  \br{1+\frac{s}{\ol{n}_\core}}^{3l-d}\exp\br{B_1} \geq \exp\br{-\frac{3l-d}{2} + O\br{\sqrt{\frac{l^3}{n_\complex}}} + O(1)}.
\end{equation}
Now \eqref{SigmaC:lower:start} and \eqref{SigmaC:lower:2}, together with $\abs{I_s}=\delta\sqrt{n_\complex}$ prove \ref{ncore:lower}.

\medskip

It remains to prove \ref{ncore:main}. First observe that 
if we take the sum \eqref{eq:core:2} over all $r\in\Z$ and normalise, we obtain a normally distributed random
variable $X$ with mean $n_\complex A_1 = O(l)$ and variance $n_\complex$. Applying the Chernoff bound \eqref{eq:Chernoff:N} to $X$, we deduce
\begin{equation*}
  \sum_{\abs{r-n_\complex A_1} > \frac{\delta}{2} \ol{n}_\core}\exp\br{-\frac{r^2}{2n_\complex}+rA_1+A_2} \le
  2\exp\br{-\delta^2\frac{3l-d}{8}}\Sigma_\core^+.
\end{equation*}
Note that $\abs{r-n_\complex A_1}<\frac{\delta}{2} \ol{n}_\core$ implies that $n_\core\in I_\core^\delta(n_\complex,l,d)$
for sufficiently large $n_\complex$ and $l=o(n_\complex)$, because then $n_\complex A_1=o(n_\complex)$. Therefore,
\begin{equation*}
  \frac{\displaystyle\sum_{n_\core\notin I_\core^\delta}\frac{(n_\complex)_{n_\core}}{n_\complex^{n_\core}}n_\core(n_\core+\nu l-1)_{3l-d-1}}
  {\displaystyle\sum_{n_\core\in I_\core^\delta}\frac{(n_\complex)_{n_\core}}{n_\complex^{n_\core}}n_\core(n_\core+\nu l-1)_{3l-d-1}}
  \le \exp\br{-\delta^2\frac{3l-d}{8}+\Theta(1)+\Theta\br{\sqrt{\frac{l^3}{n_\complex}}}}.
\end{equation*}
Now $\delta^2\frac{3l-d}{8} = \Theta(l)$, $\sqrt{\frac{l^3}{n_\complex}}=o(l)$, and the fact that $l\to\infty$ finish the proof of \ref{ncore:main}.\qed

%-------------------------------------------------------------------------------------------------------------------

\proofof{lem:SigmaD}

We start by proving \ref{deficiency:upper}. We apply
\begin{equation*}
  \frac{(3l-d)^{(3l-d+2)/2}}{(3l-d)!} \stackrel{\eqref{stirling:bound}}{\le} \frac{e^{3l-d}}{\sqrt{2\pi}}(3l-d)^{-\frac{3l-d-1}{2}} \stackrel{\eqref{explogu}}{\le}
  \frac{e^{3l-\frac{d}{2}}}{\sqrt{2\pi}}(3l)^{-\frac{3l-d-1}{2}}
\end{equation*}
and \Cref{complex:bounds,lem:SigmaC} to deduce that
\begin{align*}
  \Sigma_d &\leq \frac{\exp\br{a_\core^++b_\core^+\sqrt{\frac{l^3}{n_\complex}}}}{\sqrt{2\pi}}e^{3l}(3l)^{-\frac{3l-1}{2}}
  \sum_{d=0}^{2l}\binom{2l}{d}\br{\frac{108l}{n_\complex}}^{\frac{d}{2}},
\end{align*}
proving \ref{deficiency:upper} with $a_d^+=a_\core^+-\frac12\log(2\pi)$ and $b_d^+=b_\core^++2\sqrt{108}$.

For \ref{deficiency:lower}, first note that we have a lower bound for $\Sigma_d$ if we restrict the sum \eqref{complex:sum:def} to 
$0\leq d\leq\left\lfloor\frac{2l}{7}\right\rfloor$. By analogous arguments as for the upper bound, we deduce that
\begin{equation*}
 \Sigma_d\geq \frac{\exp\br{a_\core^-+b_\core^-\sqrt{\frac{l^3}{n_\complex}}}}{e}e^{3l}(3l)^{-\frac{3l-1}{2}}
 \sum_{d=0}^{\frac{2l}{7}}\binom{2l}{d}\br{\frac{3l}{216^2en_\complex}}^{\frac{d}{2}}.
\end{equation*}
The sum above can be extended to a sum $Y=\sum_{d=0}^{2l}\binom{2l}{d}y^d$ with $y=o(1)$. Normalising this sum results in
a binomially distributed random variable $X=\text{Bi}(2l,p)$ with $p=\frac{y}{1+y}$ and $\expec{X}=\Theta\br{\sqrt{l^3/n_\complex}}$. 
If $\expec{X}\to0$, then the main contribution to $Y$ is provided by the term with $d=0$. Otherwise, the Chernoff bound 
\eqref{eq:Chernoff:B} yields that the main contribution to $X$---and thus also to $Y$---is provided by an interval contained in the range $0\le d\le \frac{2l}{7}$.
Thus, with \eqref{explogl} we deduce that
\begin{align*}
  \Sigma_d &\geq \frac{\exp\br{a_\core^-+b_\core^-\sqrt{\frac{l^3}{n_\complex}}}}{e(1+o(1))}e^{3l}(3l)^{-\frac{3l-1}{2}}
 \exp\br{\frac{\sqrt{3}}{108\sqrt{e}}\sqrt{\frac{l^3}{n_\complex}}-\frac{\sqrt{3}l^2}{216\sqrt{e}n_\complex}}.
\end{align*}
Observing that $l^2/n_\complex=o(\sqrt{l^3/n_\complex})$, we have thus proved \ref{deficiency:lower} for any choice of
$a_d^-<a_\core^--1$ and $b_d^-<b_\core^-+\frac{\sqrt{3}}{108\sqrt{e}}$.

In order to prove \ref{deficiency:maincontribution}, it remains to show that the tail of $\Sigma_d$ has smaller order than its total value, that is
\begin{equation}\label{Sigmad:main}
 e^{b_d^+\sqrt{\frac{l^3}{n_\complex}}}\sum_{d\not\in I_d}\binom{2l}{d}\br{\frac{6 \sqrt{3l}}{\sqrt{n_\complex}}}^d=o\br{e^{b_d^-\sqrt{\frac{l^3}{n_\complex}}}}.
\end{equation}
Write
\begin{equation*}
  Z = \sum_{d=0}^{2l}\binom{2l}{d}\br{\frac{6 \sqrt{3l}}{\sqrt{n_\complex}}}^d.
\end{equation*}
For $\sqrt{\frac{l^3}{n_\complex}}\to 0$, the exponential terms in \eqref{Sigmad:main} are both $1+o(1)$ and the sum on the left hand side is
$o(1)$, because its range does not include the main contribution of the binomial sum $Z$, which is located at $d=0$.

If $\sqrt{\frac{l^3}{n_\complex}}\to c\in\R^+$, then both exponential terms in \eqref{Sigmad:main} are $\Theta(1)$.
For any fixed $h=h(n_\complex)=\omega(1)$, we deduce from \eqref{eq:Chernoff:B}, applied to the normalised sum $Z$,
\begin{equation*}
  \sum_{d>h}\binom{2l}{d}\br{\frac{6 \sqrt{3l}}{\sqrt{n_\complex}}}^d \le \exp\br{-ch}
\end{equation*}
for some constant $c>0$, which proves \eqref{Sigmad:main}.

Finally, if $\sqrt{\frac{l^3}{n_\complex}}\to\infty$, we can choose $\beta_d^+$ sufficiently large so that 
\eqref{eq:Chernoff:B} yields
\begin{equation*}
  \sum_{d>\beta_d^+\sqrt{\frac{l^3}{n_\complex}}}\binom{2l}{d}\br{\frac{6 \sqrt{3l}}{\sqrt{n_\complex}}}^d \le \exp\br{-(b_d^+-b_d^-+1)\sqrt{\frac{l^3}{n_\complex}}},
\end{equation*}
which proves \eqref{Sigmad:main} also in this last case.\qed

%-------------------------------------------------------------------------------------------------------------------

\proofof{coredeficiency}

The typical range for $d(G)$ follows directly from \Cref{lem:SigmaD}\ref{deficiency:maincontribution}. Substituting
this deficiency in the formulas for the main contribution for $n_\core$ from \Cref{lem:SigmaC} yields the 
typical order of the core.\qed

%-------------------------------------------------------------------------------------------------------------------

\proofof{complex:number}

This follows directly from \eqref{complex:step2} and \eqref{eq:Sigmad}.\qed

%-------------------------------------------------------------------------------------------------------------------

\proofof{lem:exceptional}

We prove \Cref{lem:exceptional} using the lower bound on $\abs{\class{\general}_g^*(n,m)}$ from \Cref{nQandl:lower}.
It is important to note that vice versa, the proof of \Cref{nQandl:lower} does \emph{not} rely on \Cref{lem:exceptional}.

Suppose first $n_\complex=0$, i.e.~the complex part is empty and the graph only consists of trees and unicyclic 
components. In this case \Cref{britikov}\ref{britikov:sup2} implies that the number of such graphs satisfies
\begin{equation*}
  \abs{\class{\outside}(n,m)} \le \Theta(1)n^m2^{m-n}e^{n-\frac{m^2}{n}}.
\end{equation*}
Comparing this to the lower bound from \Cref{nQandl:lower} shows that
\begin{equation*}
  \frac{\abs{\class{\outside}(n,m)}}{\abs{\class{\general}_g^*(n,m)}} \le e^{-l_1} = o(1).
\end{equation*}

The remaining case is $m_\outside=0$, i.e.\
$m=n_\complex+l \ge n_\complex+1$ (recall that $n_\complex>0$ implies $l>0$). The number of such graphs is given by
\begin{equation*}
 \sum_{n_\complex\le m-1}\binom{n}{n_\complex}\abs{\class{\complex}_g(n_\complex,m)}.
\end{equation*}
The case $n_\complex=m-1$ in the sum above is of smaller order than the lower bound for $\abs{\class{\general}_g^*(n,m)}$
from \Cref{nQandl:lower}. For every $n_\complex<m-1$, \Cref{complex:number} implies that
\begin{equation}\label{quotient}
 \frac{\binom{n}{n_\complex}\abs{\class{\complex}_g(n_\complex,m)}}{\binom{n}{n_\complex}\abs{\class{\complex}_g(n_\complex,m-1)}\binom{n-n_\complex}{2}}
  = \Theta(1)n_\complex^{\frac32}(m-n_\complex)^{-\frac32}(n-n_\complex)^{-2}.
\end{equation}
In \firstsup\ and \intermediate, the right hand side of \eqref{quotient} is $O(n^{-1/2})$.
Observing that the denominator is a summand of
$\abs{\class{\general}_g^*(n,m)}$, we deduce that
\begin{equation*}
  \sum_{n_\complex\le m-1}\binom{n}{n_\complex}\abs{\class{\complex}_g(n_\complex,m)} = o\br{\abs{\class{\general}_g^*(n,m)}}
  \qquad\text{in \firstsup\ and \intermediate.}
\end{equation*}
Suppose now that we are in the second phase transition and write $I_l=[p_l,q_l]$. For $n_\complex< m-p_l$, the right hand side of \eqref{quotient} is $o(1)$
and thus
\begin{equation*}
  \sum_{n_\complex< m-p_l}\binom{n}{n_\complex}\abs{\class{\complex}_g(n_\complex,m)} = o\br{\abs{\class{\general}_g^*(n,m)}}.
\end{equation*}
For $n_\complex \ge m-p_l$, or equivalently $l \le p_l$, we have
\begin{equation*}
  \sum_{n_\complex = m-p_l}^{m-2}\binom{n}{n_\complex}\abs{\class{\complex}_g(n_\complex,m-1)}\binom{n-n_\complex}{2}
  \le \exp(-f(n))\abs{\class{\general}_g^*(n,m)},
\end{equation*}
where $f=\omega(\log n)$ is a positive valued function. From this, we deduce that
\begin{equation*}
  \sum_{n_\complex = m-p_l}^{m-2}\binom{n}{n_\complex}\abs{\class{\complex}_g(n_\complex,m)}
  \stackrel{\eqref{quotient}}{\le} \Theta\br{n^{\frac32}}\exp(-f(n))\abs{\class{\general}_g^*(n,m)} = o\br{\abs{\class{\general}_g^*(n,m)}}.
\end{equation*}
This concludes the proof of \Cref{lem:exceptional}.
\qed

%-------------------------------------------------------------------------------------------------------------------

\proofof{lem:SigmaQ}

In $\Sigma_\complex=\sum_{n_\complex}\rho\psi$ (see \eqref{psi} for the definition of $\psi$), we substitute 
$n_\complex = \ol{n}_\complex+r$.
We then have $n_\outside = n-n_\complex = \ol{n}_\outside-r$ and
$m_\outside = m-n_\complex-l = \ol{m}_\outside-r$.

With this substitution, we obtain
\begin{align*}
 \psi=\br{\frac{2}{e}}^{\ol{n}_\complex+r}(\ol{n}_\complex+r)^{\frac{3l}{2}-1}(\ol{n}_\outside-r)^{-r-\frac{1}{2}}(\ol{m}_\outside-r)^{-\ol{m}_\outside+r-\frac12}\exp\br{f_d}.
\end{align*}
Because $n_\complex,l$ are admissible, we have $l = O(n_\complex)$ and thus
\begin{equation}\label{complex:f}
f_d \leq a_d^++b_d^+\sqrt{\frac{l^3}{n_\complex}} = O(l).
\end{equation}
If in addition \eqref{eq:smalll} holds, then $l=o\br{\ol{n}_\complex}$ and thus, for every fixed
$h(n)=\omega(1)$,
\begin{equation}\label{complex:f:smalll}
  f_d \le a_d^+ + o(1) l,
\end{equation}
whenever $r \ge -\ol{n}_\complex+hl$.
In either case, we distinguish whether $r>0$ or $r\le0$.

Let $\Sigma_{r>0}$ be the part of $\Sigma_\complex$ consisting of the summands with $r>0$.
We bound $\rho(n_\outside,m_\outside)$ from above by $1$. Additionally we claim that
\begin{equation}\label{psi:rpositive}
 \br{\frac{2}{e}}^{r}(\ol{n}_\outside-r)^{-r}(\ol{m}_\outside-r)^{-\ol{m}_\outside+r} < \ol{m}_\outside^{-\ol{m}_\outside}\exp\br{-\frac{r^3}{24\ol{m}_\outside^2}}.
\end{equation}
Indeed, for $r\ge0$, the quotient of the two sides in \eqref{psi:rpositive} has a unique maximum at $r=0$, where we have equality.
Furthermore, there exists a constant $c>0$ with
\begin{equation}\label{psi:rootterms}
  (\ol{n}_\outside-r)^{-\frac{1}{2}}(\ol{m}_\outside-r)^{-\frac12} \le c\ol{m}_{\outside}^{-1}\exp\br{\frac{r^3}{216\ol{m}_\outside^2}}.
\end{equation}
Now \eqref{complex:f}, \eqref{psi:rpositive}, and \eqref{psi:rootterms} yield
\begin{equation}\label{SigmaQ:positive}
 \Sigma_{r>0}\le \br{\frac{2}{e}}^{\ol{n}_\complex}\ol{m}_\outside^{-\ol{m}_\outside-1}\exp\br{O(l)}\sum_r (\ol{n}_\complex+r)^{\frac{3l}{2}-1}\exp\br{-\frac{r^3}{27\ol{m}_\outside^2}}.
\end{equation}
If in addition \eqref{eq:smalll} holds, we can replace $\exp(O(l))$ by $(1+o(1))^l$.
The summand above is maximised at the (not necessarily integral) unique positive solution $r_0$ of
\begin{equation*}
 r_0^3+r_0^2\ol{n}_\complex=9\ol{m}_\outside^2\br{\frac{3l}{2}-1}.
\end{equation*}

Suppose first that \eqref{eq:casesl} holds, that is,
$\ol{n}_\complex^3\geq 9\ol{m}_\outside^2\br{\frac{3l}{2}-1}$. Then 
\begin{equation}\label{range:A}
 \frac{1}{2}\sqrt{\frac{9\ol{m}_\outside^2\br{\frac{3l}{2}-1}}{\ol{n}_\complex}}\leq r_0\leq \sqrt{\frac{9\ol{m}_\outside^2\br{\frac{3l}{2}-1}}{\ol{n}_\complex}}
\end{equation}
and thus
\begin{align*}
 (\ol{n}_\complex+r_0)^{\frac{3l}{2}-1}\exp\br{-\frac{r_0^3}{27\ol{m}_\outside^2}}
  &\stackrel{\eqref{explogu}}{\le}
 \ol{n}_\complex^{\frac{3l}{2}-1} \exp\br{\frac{r_0\br{\frac{3l}{2}-1}}{\ol{n}_\complex}-\frac{r_0^3}{27\ol{m}_\outside^2}}\\
  &\stackrel{\eqref{eq:casesl},\eqref{range:A}}{\le} \ol{n}_\complex^{\frac{3l}{2}-1} \exp\br{O(l)}.
\end{align*}
Summing over $1\le r\le \ol{m}_\outside-1$, we deduce that
\begin{equation*}
  \Sigma_{r>0}\le \br{\frac{2}{e}}^{\ol{n}_\complex}\ol{n}_\complex^{\frac{3l}{2}-1}\ol{m}_\outside^{-\ol{m}_\outside}
  \exp\br{O(l)},
\end{equation*}
which proves \eqref{eq:SigmaQ} for $\Sigma_{r>0}$ if \eqref{eq:casesl} holds.
If the stronger condition \eqref{eq:smalll} is satisfied, the factor $\exp\br{O(l)}$ improves to $\exp\br{O\br{\sqrt{\epsilon} l}}=\exp\br{o(1) l}$,
proving \eqref{eq:SigmaQ:better} for $\Sigma_{r>0}$.

Now consider the case $\ol{n}_\complex^3<9\ol{m}_\outside^2\br{\frac{3l}{2}-1}$. Then
\begin{equation}\label{range:B}
 \frac{1}{2}\sqrt[3]{9\ol{m}_\outside^2\br{\frac{3l}{2}-1}}\leq r_0\leq 2\sqrt[3]{9\ol{m}_\outside^2\br{\frac{3l}{2}-1}}
\end{equation}
and hence
\begin{equation*}
  (\ol{n}_\complex+r_0)^{\frac{3l}{2}-1}\exp\br{-\frac{r_0^3}{27\ol{m}_\outside^2}}
  \le
  (3r_0)^{\frac{3l}{2}-1}\exp\br{O(l)}.
\end{equation*}
Summing over less than $\ol{m}_\outside$ values for $r$, we deduce that
\begin{equation*}
 \Sigma_{r>0} \le \br{\frac{2}{e}}^{\ol{n}_\complex}r_0^{\frac{3l}{2}-1}\ol{m}_\outside^{-\ol{m}_\outside}
  \exp\br{O(l)}.
\end{equation*}
Together with \eqref{range:B}, this proves \eqref{eq:SigmaQ} for $\Sigma_{r>0}$ in the case that \eqref{eq:casesl} is violated. 

Finally, consider the part $\Sigma_{r\le0}$ of $\Sigma_\complex$ consisting of the summands with $r\le0$.
Observe that $-\ol{n}_\complex+1\le r\le0$; in particular, the case $r\le 0$ only occurs if $\ol{n}_\complex>0$. We use \Cref{britikov}\ref{britikov:sup2} as 
an upper bound for $\rho=\rho(\ol{n}_\outside-r,\ol{m}_\outside-r)$ to deduce
\begin{align}\label{eq:negativer:start}
  \rho  \psi \le c\br{\frac{2}{e}}^{\ol{n}_\complex}(\ol{n}_\complex+r)^{\frac{3l}{2}-1}\ol{m}_\outside^{-\ol{m}_\outside-\frac12}\exp(f_d).
\end{align}
We bound the factor $\exp(f_d)$ by \eqref{complex:f}.
Furthermore, $(\ol{n}_\complex+r)^{\frac{3l}{2}-1}\le \ol{n}_\complex^{\frac{3l}{2}-1}$, because $r\le 0$. 
Summing over $r$, we deduce that
\begin{equation*}
  \Sigma_{r\le0} \le c\br{\frac{2}{e}}^{\ol{n}_\complex}\ol{n}_\complex^{\frac{3l}{2}}\ol{m}_\outside^{-\ol{m}_\outside-\frac12}\exp(O(l)).
\end{equation*}
This proves \eqref{eq:SigmaQ} for $\Sigma_{r\le0}$, independent of whether \eqref{eq:casesl} is satisfied.

Finally, suppose that \eqref{eq:smalll} holds. Then in \eqref{eq:negativer:start}, we bound the factor $\exp(f_d)$ by
\eqref{complex:f:smalll} for $r\ge r_1 := -\ol{n}_\complex+hl$ and deduce by analogous arguments as above that
\begin{equation*}
  \sum_{r=r_1}^{0}\rho  \psi \le \Theta(1)\br{\frac{2}{e}}^{\ol{n}_\complex}\ol{n}_\complex^{\frac{3l}{2}}\ol{m}_\outside^{-\ol{m}_\outside-\frac12}(1+o(1))^l.
\end{equation*}
For $r<r_1$, observe that Euler's formula yields $r \ge r_2 := -\ol{n}_\complex+\Theta(l)$. In this range, the summand
$\rho\psi$ is maximised at the upper bound $r=r_1-1$; this yields
\begin{equation*}
  \sum_{r=r_2}^{r_1-1}\rho  \psi \le \Theta(1)\br{\frac{2}{e}}^{\ol{n}_\complex}\br{hl}^{\frac{3l}{2}}\ol{m}_\outside^{-\ol{m}_\outside-\frac12}(1+o(1))^l.
\end{equation*}
If we choose $h$ to be growing slowly enough so that $hl = o(\ol{n}_\complex)$, then this proves \eqref{eq:SigmaQ:better} for
$\Sigma_{r<0}$.

The trivial observation $\Sigma_\complex = \Sigma_{r>0}+\Sigma_{r\le0}$ finishes the proof.\qed

%-------------------------------------------------------------------------------------------------------------------

\proofof{lem:SigmaQ:small}

Like in the proof of \Cref{lem:SigmaQ}, we distinguish the cases $r>0$ and $r\le0$ as well as
whether \eqref{eq:casesl} holds or not.

First consider $\Sigma_{r>0}$ when \eqref{eq:casesl} holds. Then \eqref{range:A} implies $r_0\le\ol{n}_\complex$, which yields
\begin{align*}
  \sum_{r=1}^{\ol{n}_\complex}\br{\ol{n}_\complex+r}^{\frac{3l}{2}-1}\exp\br{-\frac{r^3}{27\ol{m}_\outside^2}}
  &\le \ol{n}_\complex^{\frac{3l}{2}}\br{1+\frac{r_0}{\ol{n}_\complex}}^{\frac{3l}{2}-1}\exp\br{-\frac{r_0^3}{27\ol{m}_\outside^2}}\\
  &\le \ol{n}_\complex^{\frac{3l}{2}}\exp(O(l)).
\end{align*}
The sum over the remaining values for $r$ is bounded by the integral
\begin{equation*}
  \int_{\ol{n}_\complex}^{\infty}\br{2r}^{\frac{3l}{2}-1}\exp\br{-\frac{r^3}{27\ol{m}_\outside^2}}\dd r
  \le \ol{m}_\outside^l\Gamma\br{\frac{l}{2}}\exp(O(l))
  = \ol{m}_\outside^ll^{\frac{l}{2}}\exp(O(l)).
\end{equation*}
Now \eqref{eq:casesl}, \eqref{SigmaQ:positive}, and the fact
that $\ol{n}_\complex = 2m-n-2l < \edgesfirst n^{2/3}$ prove \eqref{eq:SigmaQ2} for $\Sigma_{r>0}$.

If \eqref{eq:casesl} is violated, we split $\Sigma_{r>0}$ into the sums for $1\le r\le r_0$ and $r_0<r$.
Observe that \eqref{range:B} implies $\ol{n}_\complex<2r_0$. Thus, the sum for $1\le r\le r_0$ is smaller than
$\ol{m}_\outside^ll^{\frac{l}{2}}\exp(O(l))$, while the sum for $r_0<r$ is bounded by the integral
\begin{equation*}
  \int_{r_0}^{\infty}\br{3r}^{\frac{3l}{2}-1}\exp\br{-\frac{r^3}{27\ol{m}_\outside^2}}\dd r
  \le \ol{m}_\outside^ll^{\frac{l}{2}}\exp(O(l)) < \ol{m}_\outside^ll^{\frac{l}{2}-\frac13}\exp(O(l)).
\end{equation*}
Now \eqref{eq:SigmaQ2} for $\Sigma_{r>0}$ follows from \eqref{SigmaQ:positive}
and the trivial fact that $\ol{m}_\outside = O(n)$.

For $r\le0$, observe that $m_\outside = \frac{n_\outside}{2}-\frac{r}{2}$. Furthermore, we have $n_\complex\le
\ol{n}_\complex = O(\edgesfirst n^{2/3})$ and thus
\begin{equation*}
  r = O(\edgesfirst n^{2/3})
  \qquad\text{and}\qquad
  n_\outside = (1+o(1))n.
\end{equation*}
By the assumption $\edgesfirst = o(n^{1/12})$, \Cref{britikov}\ref{britikov:sup} applies to
$\rho(n_\outside,m_\outside)$ and summing over $-\ol{n}_\complex+1\le r\le 0$ yields
\begin{align*}
  \Sigma_{r\le0} \le c\br{\frac{2}{e}}^{\ol{n}_\complex}\ol{n}_\complex^{\frac{3l}{2}}\ol{n}_\outside^{-\frac12}\ol{m}_\outside^{-\ol{m}_\outside-\frac12}\exp(O(l)).
\end{align*}
Now \eqref{eq:SigmaQ2} follows for $\Sigma_{r\le0}$ analogously to the proof of \Cref{lem:SigmaQ}, with
the additional fact $\ol{n}_\complex = O(\edgesfirst n^{2/3})$.\qed

%-------------------------------------------------------------------------------------------------------------------

\proofof{sizesl0}

By \eqref{eq:l0}, $l_0$ is positive.
We prove the order of $l_0$ separately for each of the five regimes.

\firstsup: In this regime, we have
\begin{align*}
 l_0=\frac{\phi^{2/3}(\edgesfirst n^{2/3}-2l_0)}{e^{1/3}2^{4/3}\br{\frac{n}{2}-\edgesfirst n^{2/3}+l_0}^{2/3}}.
\end{align*}
The denominator is of order $\Theta(n^{2/3})$. Thus, in order for the equality to be true, the numerator must 
be of order $\edgesfirst n^{2/3}$ and thus $l_0=\Theta(\edgesfirst)$.

\intermediate: Here, the denominator is still of order $n^{2/3}$ and 
the numerator is of order $\Theta(n)$ and thus $l_0=\Theta(n^{1/3})$.

\secondsub: The numerator is of order $\Theta(n)$ and thus
\begin{equation*}
 l_0=\frac{\Theta(n)}{(l_0-\frac12\edgessecond n^{3/5})^{2/3}}.
\end{equation*}
If $l_0=\Omega(\abs{\edgessecond} n^{3/5})$, then we have
$l_0=\Theta\br{\frac{n}{l_0^{2/3}}}$ and thus $l_0=\Theta(n^{3/5})=o(\abs{\edgessecond} n^{3/5})$, a contradiction.
Therefore, $l_0=o(\abs{\edgessecond} n^{3/5})$ and 
\begin{equation*}
 l_0=\Theta\br{\frac{n}{(\abs{\edgessecond} n^{3/5})^{2/3}}}=\Theta\br{\abs{\edgessecond}^{-2/3}n^{3/5}}.
\end{equation*}

\secondcrit: The numerator has order $\Theta(n)$. For the 
 denominator we have a contradiction similar to the previous case if $l_0$ is not $\Theta(n^{3/5})$. 
 Furthermore, the denominator has order $\Theta\br{n^{3/5}}$.

\secondsup: The numerator is $\Theta(n)$ and we obtain a contradiction if there is no cancellation 
 in the denominator. Thus we set $l_0=\frac12\edgessecond n^{3/5}+r$ with $r=o(\edgessecond n^{3/5})$ and deduce that 
$r=\Theta\br{\edgessecond^{-3/2}n^{3/5}}$.
\qed

%-------------------------------------------------------------------------------------------------------------------
\proofof{nQandl:lower}

By \Cref{sizesl0}, we have $0<l=o(\ol{n}_\complex)$ and $0<\ol{n}_\complex<n$. Thus, $\ol{n}_\outside$ and
$\ol{m}_\outside$ are also positive. Therefore, we have
$\class{\complex}_g(\ol{n}_\complex,\ol{n}_\complex+l) \not= \emptyset$ and $\class{\outside}(\ol{n}_\outside,\ol{m}_\outside)
\not= \emptyset$, showing that the given value $l$ and $n_\complex = \ol{n}_\complex$ are admissible. Recall that
\begin{equation*}
  \Sigma_\complex = \sum_{n_\complex}\rho(n_\outside,m_\outside)\psi(n_\complex,l).
\end{equation*}
Observe that (at least) all $n_\complex$ with $\ol{n}_\complex \le n_\complex \le \ol{n}_\complex + \ol{m}_\outside -1$
are admissible in this sum. For each such $n_\complex$, we have $m_\outside \le \frac{n_\outside}{2}$ and thus
\Cref{britikov}\ref{britikov:crit} yields
\begin{equation*}
  \Sigma_\complex \ge \Theta(1)\sum_{n_\complex=\ol{n}_\complex}^{\ol{n}_\complex+\ol{m}_\outside-1}\psi(n_\complex,l).
\end{equation*}
Set $n_\complex = \ol{n}_\complex + r$. There exists a $c>0$ such that
\begin{equation*}
  \psi(\ol{n}_\complex+r,l) \ge \br{\frac{2}{e}}^{\ol{n}_\complex}\br{\ol{n}_\complex+r}^{\frac{3l}{2}-1}\ol{m}_\outside^{-\ol{m}_\outside-1}\exp\br{f_d-\frac{r^3}{12\ol{m}_\outside^2}}
\end{equation*}
holds for $0\le r\le c\ol{m}_\outside$. The factor $\br{\ol{n}_\complex+r}^{\frac{3l}{2}-1}\exp\br{-\frac{r^3}{12\ol{m}_\outside^2}}$
is increasing until the unique positive solution $r_0$ of
\begin{equation*}
  r_0^3+r_0^2\ol{n}_\complex=4\ol{m}_\outside^2\br{\frac{3l}{2}-1}.
\end{equation*}
The assumptions on the size of $l$ imply that $l$
satisfies \eqref{eq:casesl}, which in turn yields $r_0=\Theta\br{\ol{m}_\outside^{2/3}}$. Therefore, for $1\le r\le r_0$, we have
\begin{equation*}
  \psi(\ol{n}_\complex+r,l) \ge \br{\frac{2}{e}}^{\ol{n}_\complex}\ol{n}_\complex^{\frac{3l}{2}-1}\ol{m}_\outside^{-\ol{m}_\outside-1}\exp\br{f_d(\ol{n}_\complex+r,l)-\frac{1}{12\ol{m}_\outside^2}}.
\end{equation*}
Let $\tilde{n}_\complex = \ol{n}_\complex+r$ be the value that minimises
$f_d(\ol{n}_\complex+r,l)$ for $1\le r\le r_0$; then $\tilde{n}_\complex = \ol{n}_\complex+O(\ol{m}_\outside^{2/3})$, since
$r \le r_0$. This proves the 
lower bound for $\Sigma_\complex$. The lower bound for $\abs{\class{\general}_g^*(n,m)}$ follows directly 
from \eqref{second:start}, the bound for $\Sigma_\complex$, and the fact that $l_1$ is admissible.\qed

%-------------------------------------------------------------------------------------------------------------------
\proofof{lem:caseB}

First observe that there exists $l_b>0$ such that \eqref{eq:casesl} is violated precisely when $l\ge l_b$. In the
first supercritical regime, we have $l_b=\Theta\br{\edgesfirst^3}$, in all other regimes $l_b=\Theta(n)$.

By \Cref{lem:SigmaQ}, we have
\begin{equation*}
  \tilde{\Sigma}_l\leq \sum_l n^{\frac{3}{2}} l^{-l+\frac{5g}{2}-\frac{10}{3}}(n-m+l)^{m-n-\frac{5}{3}}\exp\br{O(l)},
\end{equation*}
where the sum is taken over all $l\ge l_b$. The sum on the right hand side is bounded from above by
a geometric sum $\sum_l\exp\br{-cl}$ with $c>0$ and thus
\begin{equation*}
  \tilde{\Sigma}_l\leq (1+o(1)) n^{\frac{3}{2}} l_b^{-l_b+\frac{5g}{2}-\frac{10}{3}}(n-m+l_b)^{m-n-\frac{5}{3}}\exp\br{O(l_b)}.
\end{equation*}
Comparing this with the lower bound for $\abs{\class{\general}_g^*(n,m)}$ from \Cref{nQandl:lower} and implementing
\eqref{eq:l1}, we deduce that
\begin{equation*}
  \frac{n^{n+1/2}\br{\frac{e}{2}}^m\tilde{\Sigma}_l}{\abs{\class{\general}_g^*(n,m)}} \le \br{2m-n-2l_1}n^{\frac16}l_b^{-l_b+\frac{5g}{2}-\frac{10}{3}}\br{\frac{n-m+l_b}{n-m+l_1}}^{m-n}\exp\br{O(l_b)}.
\end{equation*}
The right hand side is $o(1)$, unless we are in the first supercritical regime and $\edgesfirst$ (and thus also
$l_b$) is too small for the term $l_b^{-l_b}$ to compensate the polynomial terms in $n$. For this to be the case,
we would in particular have $\edgesfirst=o\br{n^{1/12}}$. For such $\edgesfirst$, we have the stronger upper bound 
for $\tilde{\Sigma}_l$ provided by \Cref{lem:SigmaQ:small}, which is smaller than the one from \Cref{lem:SigmaQ} by a
factor of $\edgesfirst^{-1}n^{5/6}$. Thus, for these $\edgesfirst$, we have
\begin{align*}
  \frac{n^{n+1/2}\br{\frac{e}{2}}^m\tilde{\Sigma}_l}{\abs{\class{\general}_g^*(n,m)}} &\le \edgesfirst\br{2m-n-2l_1}n^{-\frac23}l_b^{-l_b+\frac{5g}{2}-\frac{10}{3}}\br{\frac{n-m+l_b}{n-m+l_1}}^{m-n}\exp\br{O(l_b)}\\
  &\le \edgesfirst^2l_b^{-l_b}\exp\br{O(l_b)},
\end{align*}
which is $o(1)$, because $l_b=\Theta\br{\edgesfirst^3}$.\qed

%-------------------------------------------------------------------------------------------------------------------

\proofof{lem:boundfd}

We first show that for $d=o(l)$ we have
\begin{equation}\label{eq:fcore}
 \abs{f_\core(d,(1-\epsilon)n_\complex,l)-f_\core(n_\complex,l,d)} = o(\epsilon l).
\end{equation}
By \eqref{eq:Sigmacore}, we have
\begin{equation}\label{eq:cond1:fcore}
 \frac{\Sigma_\core(d,(1-\epsilon)n_\complex,l)}{\Sigma_\core(n_\complex,l,d)}=(1-\epsilon)^{\frac{3l-d+1}2}\exp\br{f_\core(d,(1-\epsilon)n_\complex,l)-f_\core(n_\complex,l,d)}.
\end{equation}

We can also compare the summands of the two terms $\Sigma_\core(d,(1-\epsilon)n_\complex,l)$ and 
$\Sigma_\core(n_\complex,l,d)$ separately. Denote the summands by
\begin{equation*}
  s(n_\core,d,n_\complex,l) = \frac{(n_\complex)_{n_\core}}{n_\complex^{n_\core}}n_\core(n_\core+\nu l-1)_{3l-d-1}\,.
\end{equation*}
Then we have for $1\le n_\core = o(n_\complex)$
\begin{align*}
  \frac{s(n_\core,d,(1-\epsilon)n_\complex,l)}{s(n_\core,d,n_\complex,l)}%=\frac{((1-\epsilon)n_\complex)_{n_\core}}{(1-\epsilon)^{n_\core}(n_\complex)_{n_\core}}
   &=\frac{\binom{(1-\epsilon)n_\complex}{n_\core}}{(1-\epsilon)^{n_\core}\binom{n_\complex}{n_\core}}\\
   &\stackrel{\eqref{binomial:refined}}{=} \Theta(1)\br{1+\frac{\epsilon n_\core}{(1-\epsilon)n_\complex-n_\core}}^{(1-\epsilon)n_\complex-n_\core}\br{1-\frac{n_\core}{n_\complex}}^{\epsilon n_\complex}\\
   &\stackrel{\eqref{explog}}{=} \Theta(1)\exp\br{-(1+o(1))\frac{\epsilon n_\core^2}{2n_\complex}}.
\end{align*}
There exists an interval $I$ that contains the ranges of the main contribution to both $\Sigma_\core(n_\complex,l,d)$ 
and $\Sigma_\core(d,(1-\epsilon)n_\complex,l)$, such that $n_\core=(1+o(1))\sqrt{n_\complex(3l-d)}$, and thus in particular $n_\core = o(n_\complex)$, 
for all $n_\core\in I$. Then for $n_\core\in I$ and $d=o(l)$,
\begin{equation*}
  \frac{s(n_\core,d,(1-\epsilon)n_\complex,l)}{s(n_\core,d,n_\complex,l)} = \Theta(1)\exp\br{-\br{\frac32+o(1)} \epsilon l}.
\end{equation*}
Summing over $n_\core\in I$, we deduce that
\begin{align*}
  \Sigma_\core(d,(1-\epsilon)n_\complex,l)
  = \Theta(1)\exp\br{-\br{\frac32+o(1)} \epsilon l}\Sigma_\core(n_\complex,l,d).
\end{align*}
Combining this with \eqref{eq:cond1:fcore} and the condition $\epsilon l=\omega(1)$ yields \eqref{eq:fcore}.

\Cref{lem:SigmaD} yields
\begin{equation}\label{eq:fd}
 \frac{\Sigma_d((1-\epsilon)n_\complex,l)}{\Sigma_d(n_\complex,l)}=\exp\br{f_d((1-\epsilon)n_\complex,l)-f_d(n_\complex,l)}.
\end{equation}
Suppose that $J$ is an interval that contains the ranges of the main contributions to both $\Sigma_d(n_\complex,l)$ and
$\Sigma_d((1-\epsilon)n_\complex,l)$, such that $d\le d_0=o(l)$ for all $d\in J$. Denote the summands of $\Sigma_d(n_\complex,l)$ by
\begin{equation*}
  s_d(n_\complex,l) = \binom{2l}{d}\frac{(3l-d)^{(3l-d+2)/2}e^{d/2}\tau^d}{(3l-d)!n_\complex^{d/2}}\exp(f_\core(n_\complex,l,d)).
\end{equation*}
Recall that $\tau=\tau(d,l)$ does \emph{not} depend on $n_\complex$. With 
\eqref{eq:fcore}, we have
\begin{equation*}
 \frac{s_d((1-\epsilon)n_\complex,l)}{s_d(n_\complex,l)}=(1-\epsilon)^{-\frac{d}{2}} \exp\br{o(\epsilon l)}
\end{equation*}
for $d\in J$. Summing over $J$ and comparing with \eqref{eq:fd} proves the lemma.\qed

%-------------------------------------------------------------------------------------------------------------------
\proofof{nQandl}

Let us write $I_l(n)=[p_l(n),q_l(n)]$ and $I_\complex^h(n,l) = [p_\complex(n,l),q_\complex(n,l)]$. Without loss of
generality $p_l<l_1<q_l$.
We first prove that the main contribution with respect to $l$ is provided by $l\in I_l(n)$. To that
end, we bound the tail of the sum (the part with $l\notin I_l(n)$) from above and prove that this upper bound
has smaller order than the lower bound from \Cref{nQandl:lower}.

Observe that for $l\in I_l(n)$, we have
\begin{equation}\label{eq:orderl}
  \frac{l\ol{m}_\outside^{2/3}}{\ol{n}_\complex} = \Theta(1).
\end{equation}

For this proof, let $s_l(n,m)$ be the summand of the sum $\Sigma_l$, i.e.~$\Sigma_l=\sum_l s_l(n,m)$, and 
$s_\complex(n,m,l)$ be the summand of $\Sigma_\complex=\sum_{n_\complex}s_\complex(n,m,l)$. We need to show
that
\begin{equation*}
  T := n^{n+\frac12}\br{\frac{e}2}^m\sum_{l\not\in I_l(n)}s_l(n,m) = o\br{\abs{\class{\general}_g^*(n,m)}}.
\end{equation*}

By \Cref{lem:caseB}, we may take our sum only over $l$ that satisfy \eqref{eq:casesl}. 
If $l_2$ denotes the index where $s_l(n,m)$ takes its maximal value outside $I_l(n)$, then
\begin{equation}\label{sum:l}
 \sum_{l\not\in I_l(n)}s_l(n,m) \le n  s_{l_2}(n,m).
\end{equation}
In \firstsup, \eqref{eq:casesl} is violated for all $l\ge l_b = \Theta(\edgesfirst^3)$ and thus we have the stronger bound
\begin{equation}\label{sum:l:better}
 \sum_{l\not\in I_l(n)}s_l(n,m) \le \Theta(\edgesfirst^3) s_{l_2}(n,m).
\end{equation}
By \Cref{lem:SigmaQ}, there exists a constant $\alpha>1$ such that
\begin{equation}\label{sl:upper}
  s_l(n,m) \le n^{\frac{3}{2}}\br{\frac{2}{e}}^{2m-n}M(l,n,m;\alpha),
\end{equation}
where
\begin{equation*}
  M(l,n,m;\alpha) = l^{-\frac{3l}{2}}\br{\frac{e^2\phi}{4}}^l\br{2m-n+2l}^{\frac{3l}{2}-1}\br{n-m+l}^{m-n-l-1}\alpha^l.
\end{equation*}
By choosing $\beta_l^-$ (respectively $\eta_l^-$ or $\vartheta_l^-$) small enough and $\beta_l^+$ (respectively $\eta_l^+$ or 
$\vartheta_l^+$) large enough, we may assume that $M(l,n,m;\alpha)$ is strictly increasing (with respect to $l$) for $l\le p_l$ and
strictly decreasing for $l\ge q_l$. For $l\le p_l$, \eqref{eq:smalll} is satisfied and thus \eqref{sl:upper} holds for every
$\alpha=1+\delta$, where $\delta>0$ is any given constant. Thus,
\begin{equation}\label{sl2}
  s_{l_2}(n,m) \le n^{\frac{3}{2}}\br{\frac{2}{e}}^{2m-n}  \max\{M(p_l,n,m;1+\delta),M(q_l,n,m;\alpha)\}.
\end{equation}
In \firstsup, when $\edgesfirst=o\br{n^{1/12}}$, \Cref{lem:SigmaQ:small} together with analogous arguments gives us an upper bound
\begin{equation}\label{sl2:better}
  s_{l_2}(n,m) \le \edgesfirst n^{\frac{2}{3}}\br{\frac{2}{e}}^{2m-n}  \max\{M(p_l,n,m;\alpha),M(q_l,n,m;\alpha)\}.
\end{equation}

If $m$ is such that \eqref{sl2:better} applies and if the maximum in \eqref{sl2:better} is $M(q_l,n,m;\alpha)$, 
then \eqref{sum:l:better} and \eqref{sl2:better} yield (for large enough $\beta_l^+$)
\begin{equation*}
  \frac{T}{\abs{\class{\general}_g^*(n,m)}}
  \le
  \edgesfirst^4 e^{-l_1},
\end{equation*}
which is $o(1)$ by \Cref{sizesl0} and the fact that $\edgesfirst\to\infty$. If \eqref{sl2:better} does \emph{not} apply and
the maximum in \eqref{sl2} is $M(q_l,n,m,\alpha)$, then \eqref{sum:l}, \eqref{sl2}, and \Cref{nQandl:lower} imply that
if we choose $\beta_l^+$, $\eta_l^+$, or $\vartheta_l^+$ large enough, respectively, then
\begin{equation*}
  \frac{T}{\abs{\class{\general}_g^*(n,m)}}
  \le
  n^{\frac52} e^{-l_1},
\end{equation*}
which is $o(1)$.

If the maximum in \eqref{sl2} or \eqref{sl2:better} is $M(p_l,n,m,1+\delta)$ or $M(p_l,n,m,\alpha)$, respectively, then analogous 
considerations show that we can choose $\beta_l^-$, $\eta_l^-$, and $\vartheta_l^-$ so that for every $m=m(n)$ there exists a constant $c>0$ such that
\begin{equation*}
  \frac{T}{\abs{\class{\general}_g^*(n,m)}}
  \le
  \begin{cases}
    \edgesfirst^4 \exp\br{-cl_1} & \text{in \firstsup\ for } \edgesfirst=o(n^{1/12}),\\
    n^{\frac52} \exp\br{-c\edgessecond^{-3/2}n^{3/5}} & \text{in \secondsup},\\
    n^{\frac52} \exp\br{-cl_1} & \text{otherwise}.
  \end{cases}
\end{equation*}
In \secondcrit\ and \secondsup, the fact that we have $\alpha=1+\delta$ is essential for deducing the above bound. In all regimes---using
that $\edgessecond=o\br{(\log n)^{-2/3}n^{3/5}}$ in \secondsup---we deduce that
this upper bound is $o(1)$.
This proves that the main contribution to $\Sigma_l$ is indeed provided by $l\in I_l(n)$.

It remains to prove that for each $l\in I_l(n)$, the main contribution to $\Sigma_\complex$ is provided by $n_\complex \in
I_\complex^h(n,m,l)$. We substitute $n_\complex=\ol{n}_\complex+r$.

First consider the case $n_\complex<p_\complex = \ol{n}_\complex - h\ol{m}_\outside^{2/3}$, i.e. $r< -h\ol{m}_\outside^{2/3}$.
We shall split the sum into the three parts $-v\ol{m}_\outside^{2/3} \le r \le -h\ol{m}_\outside^{2/3}$, $-w\ol{m}_\outside^{2/3} \le
r \le -v\ol{m}_\outside^{2/3}$, and $r\le-w\ol{m}_\outside^{2/3}$, where
\begin{equation*}
  v := \ol{m}_\outside^{1/24}
  \qquad\text{and}\qquad
  w :=
  \begin{cases}
    \edgesfirst^{1/2} & \text{in \firstsup},\\
    l\ol{m}_\outside^{-2/9} & \text{otherwise}.
  \end{cases}
\end{equation*}
Observe that the interval $-w\ol{m}_\outside^{2/3} \le r \le -v\ol{m}_\outside^{2/3}$ is empty in \firstsup\ if
$\edgesfirst < \ol{m}_\outside^{1/12}$. Furthermore,
\begin{equation}\label{eq:orderw}
  w = \omega\br{\sqrt{\frac{l^3}{\ol{n}_\complex}}}
  \qquad\text{and}\qquad
  w = o\br{\frac{\ol{n}_\complex}{\ol{m}_\outside^{2/3}}} \stackrel{\eqref{eq:orderl}}{=} o(l).
\end{equation}

By \eqref{eq:negativer:start} and \Cref{nQandl:lower}, in each of the three intervals,
\begin{equation}\label{tailr:negative}
  \frac{\sum\rho\psi}{\Sigma_\complex(n,m,l)}
  \le \Theta(1)\ol{m}_\outside^{-\frac16}\sum s_r(n,m,l)
\end{equation}
with
\begin{equation*}
  s_r(n,m,l) = \br{1+\frac{r}{\ol{n}_\complex}}^{\frac{3l}{2}-1}
  \exp\br{f_d(\ol{n}_\complex+r,l)-f_d(\tilde n_\complex,l)}.
\end{equation*}

Recall that for \eqref{eq:negativer:start}, \Cref{britikov}\ref{britikov:sup2} was used to bound $\rho$.
Observe that for $-v\ol{m}_\outside^{2/3} \le r \le -h\ol{m}_\outside^{2/3}$, \Cref{britikov}\ref{britikov:sup} is 
applicable and thus \eqref{tailr:negative} holds with a factor of $\ol{m}_\outside^{-\frac23}$ instead of $\ol{m}_\outside^{-\frac16}$.
Furthermore, we claim that $f_d(\ol{n}_\complex+r,l)-f_d(\tilde n_\complex,l) = o\br{\frac{rl}{\ol{n}_\complex}}$. Indeed, in 
\firstsup\ and \intermediate, the left hand side is $O(1)$ and the claim follows by observing that $\frac{rl}{\ol{n}_\complex}=\Omega(h)$ by \eqref{eq:orderl}. 
In the second phase transition, 
such $r$ satisfy the conditions of \Cref{lem:boundfd} with $\epsilon=\Theta\br{\frac{r}{\ol{n}_\complex}}$ and thus the claim follows. 
Therefore, there exists a constant $c>0$ such that
\begin{align*}
  \frac{\sum_{r=-v\ol{m}_\outside^{2/3}}^{-h\ol{m}_\outside^{2/3}}\rho\psi}{\Sigma_\complex(n,m,l)} &\le 
  \Theta(1)\ol{m}_\outside^{-\frac23}\sum_{r=-v\ol{m}_\outside^{2/3}}^{-h\ol{m}_\outside^{2/3}}\exp\br{\br{\frac{3}{2}-o(1)}\frac{rl}{\ol{n}_\complex}}\\
  &\stackrel{\eqref{eq:orderl}}{\le} \Theta(1)\int_{h}^{\infty}e^{-cx}\dd x = \Theta(1)\exp\br{-ch} = o(1).
\end{align*}
Observe that in \firstsup, if $\edgesfirst=o(n^{1/24})$, then $r>-\ol{n}_\complex>-v\ol{m}_\outside^{2/3}$ and thus the interval
$-v\ol{m}_\outside^{2/3} \le r \le -h\ol{m}_\outside^{2/3}$ covers all cases for negative $r$. From now on, we may thus assume that
$\edgesfirst=\Omega(n^{1/24})$, which implies $w=\Omega(n^{1/48})$.

Now consider the interval $-w\ol{m}_\outside^{2/3} \le r \le -v\ol{m}_\outside^{2/3}$. In this regime, we still have
$f_d(\ol{n}_\complex+r,l)-f_d(\tilde n_\complex,l) = o\br{\frac{rl}{\ol{n}_\complex}}$ and thus
\begin{equation*}
  \frac{\sum_{r=-w\ol{m}_\outside^{2/3}}^{-v\ol{m}_\outside^{2/3}}\rho\psi}{\Sigma_\complex(n,m,l)} \le  
  \Theta(1)\ol{m}_\outside^{1/2}\exp\br{-cv} = o(1).
\end{equation*}

Finally, suppose that $r \le -w\ol{m}_\outside^{2/3}$. In this regime,
\begin{equation*}
  s_r \leq 
  \exp\br{\br{\frac{3l}{2}-1}\log\br{1+\frac{r}{\ol{n}_\complex}}+c_1\sqrt{\frac{l^3}{\ol{n}_\complex+r}}-c_2\sqrt{\frac{l^3}{\tilde{n}_\complex}}}.
\end{equation*}
The right hand side has its maximum (with respect to $r$) at $r=-w\ol{m}_\outside^{2/3}$. For this $r$, the first summand 
is negative and has order $w$ by \eqref{eq:orderl}. The other two summands are $o(w)$ by \eqref{eq:orderw}. Thus, there exists
a constant $c>0$ such that
\begin{equation*}
  \frac{\sum_{r\le-w\ol{m}_\outside^{2/3}}\rho\psi}{\Sigma_\complex(n,m,l)} \le  
  \Theta(1)n\exp\br{-cw},
\end{equation*}
which is $o(1)$, because $w=\Omega(n^{1/48})$. This finishes the proof for $r<0$.

Suppose now that $n_\complex>q_\complex = \ol{n_\complex} + h\ol{m}_\outside^{2/3}$, 
i.e. $r>h\ol{m}_\outside^{2/3}$. By \eqref{psi:rpositive}, \eqref{psi:rootterms}, and \Cref{nQandl:lower} we conclude that
\begin{equation*}
  \frac{\sum\rho\psi}{\Sigma_\complex(n,m,l)}
  \le \Theta(1)\ol{m}_\outside^{-\frac23}\sum_{r>h\ol{m}_\outside^{2/3}}
  \exp\br{\frac{3lr}{2\ol{n}_\complex}+f_d(\ol{n}_\complex+r,l)-f_d(\tilde n_\complex,l)-\frac{r^3}{27\ol{m}_\outside^2}}.
\end{equation*}
Note that for all $r$ in this sum, $\frac{rl}{\ol{n}_\complex}=o\br{\frac{r^3}{\ol{m}_\outside^2}}$. We claim 
that additionally
\begin{equation*}
  f_d(\ol{n}_\complex+r,l)-f_d(\tilde n_\complex,l)=o\br{\frac{r^3}{\ol{m}_\outside^2}}.
\end{equation*}
Indeed, this difference is $O(1)$ in \firstsup\ and \intermediate, while $\frac{r^3}{\ol{m}_\outside^2}\ge h^3=\omega(1)$.
In the second phase transition, the claim follows immediately if $r\geq\sqrt{l}\ol{m}_\outside^{2/3}$. If
$h\ol{m}_\outside^{2/3}<r<\sqrt{l}\ol{m}_\outside^{2/3}$, the conditions of \Cref{lem:boundfd} are satisfied with $\epsilon=\Theta\br{\frac{r}{\ol{n}_\complex}}$ and 
thus $f_d(\ol{n}_\complex+r,l)-f_d(\tilde n_\complex,l)=o\br{\frac{rl}{\ol{n}_\complex}}=o\br{\frac{r^3}{\ol{m}_\outside^2}}$.
Therefore, we deduce that
\begin{align*}
  \frac{\sum\rho\psi}{\Sigma_\complex(n,m,l)}
  &\le \Theta(1)\ol{m}_\outside^{-\frac23}\sum_{r>h\ol{m}_\outside^{2/3}}
  \exp\br{-\frac{r^3}{30\ol{m}_\outside^2}}\\
  &\le \Theta(1)\int_{h}^{\infty}\exp\br{-x^3}\dd x
  \le \Theta(1)\exp(-h).
\end{align*}
This finishes the proof also for $r>0$.\qed

\section{Discussion and open problems}\label{sec:discussion}

Comparing the range for $m$ that we cover in Theorems \ref{main1}--\ref{main4} with the `dense' regime
$m=\left\lfloor\mu n\right\rfloor$ for $1<\mu<3$ considered in~\cite{Chapuy2011-enumeration-graphs-on-surfaces,gimenez2009},
a gap of order $(\log n)^{2/3}$ becomes apparent---a significant improvement of~\cite{kang2012}, where the gap
had order $n^{1/3}$. The order term $\edgessecond^{-3/2}n^{3/5}$ in \Cref{main2,main4} becomes constant when
$\edgessecond=\Theta(n^{2/5})$, which matches the results from~\cite{Chapuy2011-enumeration-graphs-on-surfaces,gimenez2009}
that the giant component covers all but finitely many vertices in the dense regime. Therefore, we expect \Cref{main2,main4}
to hold for all $m=(1+o(1))n$.

The gap of order $(\log n)^{2/3}$ originates from the fact that we can only determine the number of kernels up to an
exponential error term (see \Cref{kernelpumping}) in the second phase transition. We thus believe that the key to closing
the gap would be to determine the number of kernels more exactly.

\begin{question}\label{problem:kernels}
  What is the exact value of $\abs{\class{\kernel}_g(2l-d,3l-d)}$ for any admissible $l,d\in\N$?
\end{question}

Solving \Cref{problem:kernels} would pave the way to prove \Cref{main2} for all $m=(1+o(1))n$. Moreover, it might open
the possibility to prove an analogous version of \Cref{generalstructure} in the second phase transition, thus rendering
the additional double counting argument in the proof of \Cref{main2} unnecessary; observe that this double counting argument
is responsible for the fact that the upper and lower bound on the order of $n-\abs{\component_1}$ are not quite the same.
We believe that these bounds should actually be of the same order.

\begin{conjecture}\label{conjecture:second:order}
  Let $m=\br{2+\edgessecond n^{-2/5}}\frac{n}{2}$, where $\edgessecond=\edgessecond(n)=o(n^{2/5})$.
  Then \whp\ the largest component $\component_1$ of $\general_g(n,m)$ is complex and satisfies
  \begin{equation*}
    n-\abs{\component_1} =
    \begin{cases}
      \Theta\br{\abs{\edgessecond}n^{3/5}} & \text{if }\edgessecond\to-\infty,\\
      \Theta\br{n^{3/5}} & \text{if }\edgessecond\to c\in\R,\\
      \Theta\br{\edgessecond^{-3/2}n^{3/5}} & \text{if }\edgessecond\to\infty.\\
    \end{cases}
  \end{equation*}
\end{conjecture}

Observe that in contrast to \Cref{main1,main3}, \Cref{main2} does neither provide a statement about the genus of
the largest component nor does it state the order of the $i$-th largest component for $i\ge2$. 
By~\cite{Chapuy2011-enumeration-graphs-on-surfaces}, the largest
component of $\general_g(n,\left\lfloor\mu n\right\rfloor)$ has genus $g$,
thus it is to be expected that this also holds throughout the second phase transition.

\begin{conjecture}\label{conjecture:second:genus}
  Let $m=(1+o(1))n$. Then the giant component of $\class{\general}_g(n,m)$ \whp\ is not embeddable on $\torus{g-1}$ 
  and all other components are planar.
\end{conjecture}

In view of \Cref{conjecture:second:genus}, an analogous statement to \Cref{cubic:structure} for general kernels would be
necessary. Similarly, proving \Cref{cubic:lc-size} for arbitrary kernels would open the possibility to determine the order
of the $i$-th largest component for $i\ge2$.

\begin{question}\label{problem:secondlargest:order}
  For $i\ge2$, what is the order of the $i$-th largest component in the second phase transition?
\end{question}

In view of enumeration of graphs embeddable on $\Sg$, \Cref{main4} provides an asymptotic result. Observe that the error terms
in \Cref{main4} become larger the bigger $m$ is. In particular, if we increase $\edgessecond$ to $\Theta(n^{2/5})$ in
\Cref{main4}\ref{enum:second}, then the main term of $\abs{\class{\general}_g(n,m)}$ becomes $n^n$---which matches the results
from \cite{Chapuy2011-enumeration-graphs-on-surfaces,gimenez2009} for the dense regime $m=\lfloor\mu n\rfloor$ with $\mu\in(1,3)$---but the
error term has order $\exp\br{O(n)}$. It should be possible to improve the error terms to being smaller than $(1\pm\delta)^{l_0}$ for 
every $\delta>0$ (with $l_0$ defined as in \eqref{eq:l0}) by a careful analysis of 
\Cref{lem:SigmaQ,lem:SigmaQ:small}, yet even better bounds would still be desirable.

\begin{problem}\label{problem:enum}
  Find asymptotic expressions for $\abs{\class{\general}_g(n,m)}$ with better error terms than in \Cref{main4}.
\end{problem}

It is important to note that the results in this paper apply to more general graph classes than $\class{\general}_g(n,m)$.
Indeed, the constructive decomposition that yields \eqref{eq:general}, \eqref{eq:complexcore}, and \eqref{eq:core} relies on
the fact that a graph is in $\class{\general}_g$ if and only if its kernel is in the corresponding class $\class{\kernel}_g$
of multigraphs. The only other ingredients of the proof that are specifically tailored for graphs on $\Sg$ are 
\Cref{cubic:enum,cubic:connected}, and \Cref{cubic:structure,cubic:lc-size,kernelpumping}. Recall that we saw in
\Cref{sec:remark:pumping} that \Cref{kernelpumping} holds for any class of multigraphs that is weakly addable (that is, closed
under adding an edge between two components) and closed under taking minors.

\begin{remark}\label{remark:general}
  Let $\class{X}$ be a graph class and $\class{Y}$ be a class of (weighted) multigraphs of minimum degree at least three.
  Suppose that
  \begin{enumerate}
  \item\label{class:kernel}
    a graph lies in $\class{X}$ if and only if its kernel is in $\class{Y}$;
  \item\label{class:enum}
    there are constants $c,\gamma>0$ and $k\in\R$ such that
    \begin{equation*}
      \abs{\class{Y}(2l,3l)} = (1+o(1))c\,l^{k}\gamma^{2l}(2l)!;
    \end{equation*}
  \item\label{class:connected}
    there is a constant $0<q\le1$ with
    \begin{equation*}
      \prob{Y(2l,3l) \text{ is connected}\,} \stackrel{l\to\infty}{\longrightarrow} q;
    \end{equation*}
  \item\label{class:giant}
    $\abs{\component_1(Y(2l,3l))}=2l-O_p(1)$ and for each fixed $i\in\N\setminus\{0\}$, the probability
    that $\abs{\component_1(Y(2l,3l))}=2l-2i$ is bounded away from both $0$ and $1$;
  \item\label{class:pumping}
    $\class{Y}$ is weakly addable and closed under taking minors.
  \end{enumerate}
  Then analogous statements to Theorems \ref{main1}--\ref{main4} hold for $\class{X}$.
\end{remark}

Obvious candidates for the classes $\class{X}$ and $\class{Y}$ would be (multi)graphs on non-orientable surfaces.
For such classes, \ref{class:kernel} and \ref{class:pumping} in \Cref{remark:general} are automatically satisfied, 
\ref{class:enum} and \ref{class:connected} would follow if \Cref{cubic:enum,cubic:connected} also hold for
non-orientable surfaces, and \ref{class:giant} holds if \Cref{cubic:structure} is true for non-orientable surfaces.

\begin{problem}\label{problem:nonorientable}
  Prove analogous versions of \Cref{cubic:enum,cubic:connected} and \Cref{cubic:structure} for non-orientable
  surfaces.
\end{problem}

One striking difference between $\general_g(n,m)$ and $G(n,m)$ is the order and the structure of the $i$-th largest 
component for $i\ge2$ in \firstsup\ and \intermediate. In $\general_g(n,m)$, the second largest component is much 
larger than in $G(n,m)$; $\Theta_p(n^{2/3})$ versus $o(n^{2/3})$. 
Moreover, the $i$-th largest component of $G(n,m)$ is a tree \whp. In contrast, $\general_g(n,m)$ with
positive probability has both tree components and complex components of order $\Theta_p(n^{2/3})$. It would thus 
be interesting to know whether there is a hierarchy in the size of the largest tree component and the second largest 
complex component.

\begin{question}\label{problem:secondlargest:structure}
  Given $i\ge2$, what is the probability that the $i$-th largest component of $\general_g(n,m)$ is a tree?
\end{question}

For $G(n,m)$, the giant component is in fact far better understood than it is stated in \Cref{thm:ER}. Central limit theorems
and local limit theorems provide much stronger concentration results about the order (i.e.\ the number of vertices) and the size (i.e.\ the number of edges) of the giant 
component~\cite{BehrCOKang10,coja-oghlan2014,BollobasRiordan2012A,BollobasRiordan2012B,PittelWormald05,Stepanov70} and give
more insight into the global and local structure of the giant component and its core.

\begin{problem}\label{problem:limittheorem}
  Derive central and local limit theorems for the giant component of $\general_g(n,m)$.
\end{problem}

As mentioned in \Cref{sec:intro}, the component structure of $G(n,m)$ is closely related to a Galton-Watson branching 
process. More precisely, the local structure of $G(n,\alpha\frac{n}{2})$ converges to that of a Galton-Watson tree with offspring 
distribution $\mathrm{Po}(\alpha)$ in the sense of Benjamini-Schramm local weak convergence~\cite{BenjSchramm,Karp90}.
For $\general_g(n,m)$, the additional constraint of the graph being embeddable on $\Sg$, exploration via a simple Galton-Watson
type process is not possible. This naturally raises the question if the local structure of $\general_g(n,m)$ can be described
in terms of the Benjamini-Schramm local weak convergence.

\begin{question}\label{problem:BenjSchramm}
  What is the limit of the local structure of $\general_g(n,m)$ in the sense of the Benjamini-Schramm local weak convergence?
\end{question}

The core, which plays a central role in our constructive decomposition, is also known as the \emph{$2$-core}. More generally, 
given $k\ge 2$, the \emph{$k$-core} of a graph $G$ is the largest subgraph of $G$ of minimum degree at least $k$. Like the core,
the $k$-core can be constructed by a \emph{peeling process} that recursively removes vertices of degree less than $k$. The order 
and size of the $k$-core of $G(n,m)$ has been determined in a seminal paper by Pittel, Spencer, and Wormald~\cite{PittelSpencerWormald}.
Following Pittel, Spencer, and Wormald, the $k$-core has been extensively studied~\cite{CoreForging,CoreMantle,JansonLuczakM08,Kim06,Luczakkcore,Riordan08}.
The most striking results in this area are the astonishing theorem by \Luczak~\cite{Luczakkcore} that the $k$-core for $k\ge 3$ jumps to linear 
order at the very moment it becomes non-empty, the central limit theorem by Janson and Luczak~\cite{JansonLuczakM08}, and the local
limit theorem by Coja-Oghlan, Cooley, Kang, and Skubch~\cite{CoreForging} that described---in addition to the order and size---several
other parameters of the $k$-core of $G(n,m)$. In~\cite{CoreMantle}, the same authors used a 5-type branching process in order to determine
the local structure of the $k$-core. In terms of \emph{global} structure, \cite{CoreForging} provides a randomised algorithm that constructs
a random graph with given order and size of the $k$-core.

\begin{question}\label{problem:kcore}
  What are the local and global structure of the $k$-core of $\general_g(n,m)$?
\end{question}

One of the main difficulties regarding $\general_g(n,m)$ is that while graph properties such as having a component of
a certain order are monotone for $G(n,m)$ (that is, for every fixed $n$, the probability that $G(n,m)$ has this property
is monotone for $0\le m\le\binom{n}{2}$), this is not necessarily the case for $\general_g(n,m)$. Indeed, monotonicity
of graph properties in $G(n,m)$ usually follows immediately from the equivalence between $G(n,m)$ and the \emph{random graph
process}, where we add one random edge at a time. For graphs on surfaces, however, not all edges are allowed to be added 
in the corresponding process. Thus, the process is fundamentally different from $\general_g(n,m)$. For instance, in the 
dense regime $m=\left\lfloor\mu n\right\rfloor$ with $\mu>1$, we know by~\cite{gimenez2009} that the probability that 
$P(n,m)$ is connected is bounded away from both $0$ and $1$. The planar graph \emph{process}, however, is connected \whp\
in that regime~\cite{gerke-process}. Knowing which graph properties are monotone for $\general_g(n,m)$ would yield a significant improvement
to the complexity of the arguments.

\begin{question}\label{problem:monotone}
  Which graph properties are monotone for $\general_g(n,m)$?
\end{question}

The constructive decomposition and generating functions of cubic planar graphs and their relation to 
the core of sparse planar graphs by Kang and {\L}uczak~\cite{kang2012} have been strengthened by Noy, 
Ravelomanana, and Ru{\'e} \cite{nrr} to yield an answer to a challenging open question of \Erdos\ and 
\Renyi~\cite{ErdosRenyi60} about the limiting probability of $G(n,m)$ being planar at the critical phase
\firstcrit, that is, for every constant $\edgesfirst\in\R$, the limit
$p(\edgesfirst)$ of the probability that $G\br{n,\br{1+\edgesfirst n^{-1/3}}\frac{n}{2}}$ is
planar. For graphs embeddable on a surface of positive genus, they gave a general strategy of how
to determine the corresponding probability. However, determining the \emph{exact} limiting probability
for $g\ge 1$ is still an open problem.

Furthermore, for $m$ beyond \firstcrit, we know that $G(n,m)$ \whp\ is not embeddable on any surface of fixed genus. 
This immediately raises the question what genus $g$ we need in order to embed $G(n,m)$ on $\Sg$.

\begin{question}\label{problem:Gnm:genus}
  Let $m=m(n)$ and $g=g(n)$ be given.
  \begin{enumerate}
  \item\label{genus:limit}
    When is the limiting probability of $G(n,m)$ being embeddable on $\Sg$ positive?
  \item\label{genus:whp}
    When is $G(n,m)$ embeddable on $\Sg$ \whp?
  \item\label{genus:expected}
    What is the expected genus of $G(n,m)$?
  \end{enumerate}
\end{question}

Another interesting direction, which might provide insight into the answer of \Cref{problem:Gnm:genus}, is to
consider $\general_g(n,m)$ for genus $g=g(n)$ that tends to infinity with $n$. If $g$ grows `fast enough' (e.g.\
as $\binom{n}{2}$), then $\general_g(n,m)$ will coincide with $G(n,m)$ and will thus exhibit the emergence of
the giant component, but not the second phase transition described in \Cref{main2}. For `slowly' growing $g$, on
the other hand, it is to be expected that the second phase transition does take place.

\begin{question}\label{problem:growinggenus}
  For which functions $g=g(n)$ does $\general_g(n,m)$ feature two phase transitions analogous to \Cref{main1,main2}?
\end{question}

\bibliographystyle{plain}
\bibliography{biblio}

\end{document}